\documentclass[english,american]{article}
\usepackage[T1]{fontenc}
\usepackage[latin9]{inputenc}
\usepackage{geometry}
\usepackage{color}
\usepackage{babel}
\usepackage{array}
\usepackage{float}
\usepackage{units}
\usepackage{mathtools}
\usepackage{amsmath}
\usepackage{amsthm}
\usepackage{amssymb}
\usepackage{stmaryrd}
\usepackage{graphicx}
\usepackage[unicode=true,pdfusetitle,
 bookmarks=true,bookmarksnumbered=false,bookmarksopen=false,
 breaklinks=false,pdfborder={0 0 1},backref=false,colorlinks=false]
 {hyperref}

\makeatletter

\newcommand{\noun}[1]{\textsc{#1}}
\providecommand{\tabularnewline}{\\}

\numberwithin{equation}{section}
\theoremstyle{plain}
\newtheorem{thm}{\protect\theoremname}[section]
\theoremstyle{plain}
\newtheorem{lem}[thm]{\protect\lemmaname}
\theoremstyle{definition}
\newtheorem{example}[thm]{\protect\examplename}
\theoremstyle{definition}
\newtheorem{defn}[thm]{\protect\definitionname}
\theoremstyle{remark}
\newtheorem*{acknowledgement*}{\protect\acknowledgementname}

\usepackage{algorithmic}
\usepackage{hyperref}

\makeatother

\addto\captionsamerican{\renewcommand{\acknowledgementname}{Acknowledgement}}
\addto\captionsamerican{\renewcommand{\definitionname}{Definition}}
\addto\captionsamerican{\renewcommand{\examplename}{Example}}
\addto\captionsamerican{\renewcommand{\lemmaname}{Lemma}}
\addto\captionsamerican{\renewcommand{\theoremname}{Theorem}}
\addto\captionsenglish{\renewcommand{\acknowledgementname}{Acknowledgement}}
\addto\captionsenglish{\renewcommand{\definitionname}{Definition}}
\addto\captionsenglish{\renewcommand{\examplename}{Example}}
\addto\captionsenglish{\renewcommand{\lemmaname}{Lemma}}
\addto\captionsenglish{\renewcommand{\theoremname}{Theorem}}
\providecommand{\acknowledgementname}{Acknowledgement}
\providecommand{\definitionname}{Definition}
\providecommand{\examplename}{Example}
\providecommand{\lemmaname}{Lemma}
\providecommand{\theoremname}{Theorem}

\begin{document}
\selectlanguage{english}%
\global\long\def\R{\mathbb{R}}%

\global\long\def\e{{x}}%

\global\long\def\et#1{{\e(#1)}}%

\global\long\def\ef{{x(\cdot)}}%

\global\long\def\x{{x}}%

\global\long\def\w{{w}}%

\global\long\def\m{{m}}%

\global\long\def\xt#1{{\x(#1)}}%

\global\long\def\xf{{x(\cdot)}}%

\global\long\def\d{{d}}%

\global\long\def\b{{b}}%

\global\long\def\k{{k}}%

\global\long\def\a{{a}}%

\global\long\def\u{{u}}%

\global\long\def\y{{y}}%

\global\long\def\yt#1{{\y(#1)}}%

\global\long\def\yf{{y(\cdot)}}%

\global\long\def\z{{z}}%

\global\long\def\v{{v}}%

\global\long\def\h{{h}}%

\global\long\def\s{{s}}%

\global\long\def\c{{c}}%

\global\long\def\p{{p}}%

\global\long\def\f{{f}}%

\global\long\def\rb{{r}}%

\global\long\def\rt#1{{\rb(#1)}}%

\global\long\def\rf{{r(\cdot)}}%

\global\long\def\mat#1{{#1}}%

\global\long\def\matN{\ensuremath{{N}}}%

\global\long\def\matA{\ensuremath{{A}}}%

\global\long\def\matB{\ensuremath{{B}}}%

\global\long\def\matC{\ensuremath{{C}}}%

\global\long\def\matD{\ensuremath{{D}}}%

\global\long\def\matP{\ensuremath{{P}}}%

\global\long\def\matU{\ensuremath{{U}}}%

\global\long\def\matV{\ensuremath{{V}}}%

\global\long\def\matW{\ensuremath{{W}}}%

\global\long\def\matM{\ensuremath{{M}}}%

\global\long\def\matZ{\ensuremath{{Z}}}%

\global\long\def\matR{R}%

\global\long\def\matQ{Q}%

\global\long\def\matS{S}%

\global\long\def\matSb{S_{B}}%

\global\long\def\matSw{S_{w}}%

\global\long\def\matY{Y}%

\global\long\def\matX{X}%

\global\long\def\matYhat{\hat{Y}}%

\global\long\def\matXhat{\hat{X}}%

\global\long\def\matI{I}%

\global\long\def\matzero{0}%

\global\long\def\matJ{J}%

\global\long\def\matL{L}%

\global\long\def\S#1{{\mathbb{S}_{N}[#1]}}%

\global\long\def\IS#1{{\mathbb{S}_{N}^{-1}[#1]}}%

\global\long\def\PN{\mathbb{P}_{N}}%

\global\long\def\TNormS#1{\|#1\|_{2}^{2}}%

\global\long\def\TNorm#1{\|#1\|_{2}}%

\global\long\def\InfNorm#1{\|#1\|_{\infty}}%

\global\long\def\FNorm#1{\|#1\|_{F}}%

\global\long\def\FNormS#1{\|#1\|_{F}^{2}}%

\global\long\def\UNorm#1{\|#1\|_{\matU}}%

\global\long\def\UNormS#1{\|#1\|_{\matU}^{2}}%

\global\long\def\UINormS#1{\|#1\|_{\matU^{-1}}^{2}}%

\global\long\def\ANorm#1{\|#1\|_{\matA}}%

\global\long\def\BNorm#1{\|#1\|_{B}}%

\global\long\def\ANormS#1{\|#1\|_{\matA}^{2}}%

\global\long\def\AINormS#1{\|#1\|_{\matA^{-1}}^{2}}%

\global\long\def\BINormS#1{\|#1\|_{\matB^{-1}}^{2}}%

\global\long\def\BINorm#1{\|#1\|_{\matB^{-1}}}%

\global\long\def\ONorm#1#2{\|#1\|_{#2}}%

\global\long\def\T{\textsc{T}}%

\global\long\def\pinv{\textsc{+}}%

\global\long\def\Expect#1{{\mathbb{E}}\left[#1\right]}%

\global\long\def\ExpectC#1#2{{\mathbb{E}}_{#1}\left[#2\right]}%

\global\long\def\dotprod#1#2{\left\langle #1,#2\right\rangle }%

\global\long\def\dotprodX#1#2#3{\left\langle #1,#2\right\rangle {}_{#3}}%

\global\long\def\dotprodM#1#2{\left\langle #1,#2\right\rangle _{\matM_{\matX}}}%

\global\long\def\dotprodsqr#1#2{\left\langle #1,#2\right\rangle ^{2}}%

\global\long\def\Trace#1{{\bf Tr}\left(#1\right)}%

\global\long\def\dist#1{{\bf dist}\left(#1\right)}%

\global\long\def\vectorization#1{{\bf vec}\left(#1\right)}%

\global\long\def\vecskew#1{{\bf vec_{skew}}\left(#1\right)}%

\global\long\def\nnz#1{{\bf nnz}\left(#1\right)}%

\global\long\def\blockdiag#1{{\bf blkdiag}\left(#1\right)}%

\global\long\def\vol#1{{\bf vol}\left(#1\right)}%

\global\long\def\rank#1{{\bf rank}\left(#1\right)}%

\global\long\def\diag#1{{\bf diag}\left(#1\right)}%

\global\long\def\grad#1{{\bf grad}#1}%

\global\long\def\span#1{{\bf span}#1}%

\global\long\def\hess#1{{\bf Hess}#1}%

\global\long\def\sym#1{{\bf sym}#1}%

\global\long\def\skew#1{{\bf skew}#1}%

\global\long\def\st{\,\,\,\text{s.t.}\,\,\,}%

\global\long\def\elp{\mathbb{S}^{\matB}}%

\global\long\def\elpCCA{\mathbb{S}_{\x\y}}%

\global\long\def\elpa{\mathbb{S}^{A}}%

\global\long\def\elplad{\mathbb{S}^{S_{w}+\lambda\matI_{d}}}%

\global\long\def\Stiefellda{{\bf St}_{(S_{w}+\lambda\matI_{d})}(p,d)}%

\global\long\def\ldaB{S_{w}+\lambda\matI_{d}}%

\global\long\def\elpsigx{\mathbb{S}^{\Sigma_{xx}}}%

\global\long\def\elpsigy{\mathbb{S}^{\Sigma_{\y\y}}}%

\global\long\def\elpparam#1{\mathbb{S}^{#1}}%

\global\long\def\poly#1{{\bf poly}\left(#1\right)}%

\global\long\def\id{\text{id}}%

\global\long\def\stiefel{{\bf St}}%

\global\long\def\stiefelB{{\bf St}_{\matB}}%

\global\long\def\qf#1{{\bf qf}\left(#1\right)}%

\global\long\def\qfm#1#2{{\bf qf}_{#2}\left(#1\right)}%

\global\long\def\qfmsmall#1#2{{\bf qf}_{#2}(#1)}%

\global\long\def\fcca{f_{{\bf CCA}}(\matZ)}%

\global\long\def\justfcca{f_{{\bf CCA}}}%

\global\long\def\sigmacca{\Sigma_{\nabla^{2}f_{{\bf CCA}}}}%

\global\long\def\flda{f_{{\bf LDA}}(\matW)}%

\global\long\def\justflda{f_{{\bf LDA}}}%

\global\long\def\nicehalf{\nicefrac{1}{2}}%
\selectlanguage{american}%

\title{\foreignlanguage{english}{{Riemannian optimization
with a preconditioning scheme on the generalized Stiefel manifold}}}
\author{Boris Shustin and Haim Avron\\
Tel Aviv University}
\maketitle
\begin{abstract}
Optimization problems on the generalized Stiefel manifold (and products
of it) are prevalent across science and engineering. For example,
in computational science they arise in symmetric (generalized) eigenvalue
problems, in nonlinear eigenvalue problems, and in electronic structures
computations, to name a few problems. In statistics and machine learning,
they arise, for example, in various dimensionality reduction techniques
such as canonical correlation analysis. In deep learning, regularization
and improved stability can be obtained by constraining some layers
to have parameter matrices that belong to the Stiefel manifold. Solving
problems on the generalized Stiefel manifold can be approached via
the tools of Riemannian optimization. However, using the standard
geometric components for the generalized Stiefel manifold has two
possible shortcomings: computing some of the geometric components
can be too expensive and convergence can be rather slow in certain
cases. Both shortcomings can be addressed using a technique called
Riemannian preconditioning, which amounts to using geometric components
derived by a precoditioner that defines a Riemannian metric on the
constraint manifold. In this paper we develop the geometric components
required to perform Riemannian optimization on the generalized Stiefel
manifold equipped with a non-standard metric, and illustrate theoretically
and numerically the use of those components and the effect of Riemannian
preconditioning for solving optimization problems on the generalized
Stiefel manifold.
\end{abstract}

\section{Introduction}

In this paper we consider large-scale optimization problems on the
\emph{generalized Stiefel manifold} (and products of it), i.e. optimization
with constraint spaces defined via generalized orthogonality constraints.
One well known example of a problem with a generalized orthogonality
constraints is the problem of finding the dominant generalized eigenspace
of a symmetric positive-definite (SPD) matrix pencil. Indeed, given
a pair of SPD matrices $\matA,\matB\in\R^{d\times d}$, minimizers
of $-\Trace{\matX^{\T}\matA\matX}$ subject to $\matX^{\T}\matB\matX=\matI_{p}$
(where $\matX\in\R^{d\times p}$) are bases for the subspace spanned
by the $p$ generalized eigenvectors that correspond to the $p$ largest
generalized eigenvalues of the pencil $(\matA,\matB$)~(this is a
consequence of the Courant--Fisher characterization of generalized
eigenvalues). More generally, problems with (generalized) orthogonality
constraints are prevalent across science and engineering. Examples
include, the Trust-Region Subproblem, Canonical Correlation Analysis
(CCA) \cite{hotelling1936relations}, and Fisher Linear Discriminant
Analysis \cite{fisher1936use}.

Some optimization problems with generalized orthogonality constraints
can be reformulated as (generalized) eigenvalue problems or (weighted)
Singular Value Decomposition (SVD) problems. This is true for some
of the cases mentioned in the previous paragraph. For example, CCA
on a pair of matrices $(\matX,\matY)$ amounts to computing the SVD
of $\matP^{\T}\matQ$ where $\matP$ and $\matQ$ are orthonormal
matrices whose column space spans the column space of $\matX$ and
$\matY$ (respectively)~\cite{bjorck1973numerical}. This allows
one to use direct methods, but that is unrealistic for large scale
problems.

Using iterative method in lieu of direct methods is a common modus
operandi for handling large scale problems. A natural framework for
solving optimization problems with generalized orthogonality constraints
is \emph{Riemannian optimization}~\cite{EAS98,AMS09,boumal2020intromanifolds}.
Indeed, when we have a single generalized orthogonality constraint
of the form $\matX^{\T}\matB\matX=\matI_{p}$, e.g., we want to minimize
$f(\matX)$ s.t. $\matX^{\T}\matB\matX=\matI_{p}$, one can impose
the structure of a smooth manifold on the constraint set, thereby
obtaining the generalized Stiefel manifold
\begin{equation}
\stiefelB(p,d)\coloneqq\left\{ \matX\in\R^{d\times p}:\ \matX^{\T}\matB\matX=\matI_{p}\right\} \,.\label{eq:general_Stiefel_def}
\end{equation}
(see \cite[Propositions 3.3.3 and 3.3.4]{AMS09}), and use Riemannian
optimization to minimize $f(\matX)$ s.t. $\matX\in\stiefelB(p,d)$.
If we have $k>1$ generalized orthogonality constraints , e.g., minimizing
$f(\matX_{1},\dots,\matX_{k})$ s.t. $\matX_{i}\in\stiefelB(p_{i},d_{i})\,(i=1,\dots,k)$,
as is the case in CCA (for $k=2$), then each of the constraints constrain
a disjoint set of variables, and the constraints are separable, so
they define a product of generalized Stiefel manifolds, which is a
smooth manifold as well, so Riemannian optimization can again be used.

In order to use Riemannian optimization on the generalized Stiefel
manifold $\stiefelB(p,d)$ we must further impose a Riemannian metric
on the tangent bundle of $\stiefelB(p,d)$. We refer to the Riemannian
metric naturally inherited by the scaled inner product $\dotprodX UVB=\Trace{\matU^{\T}\matB\matV}$
on $\R^{d\times p}$ as the \emph{standard metric} for (see~\cite[Section  3.6]{AMS09}
for explanation on how a Riemannian metric is inherited from an ambient
space in a natural way). Indeed, for the Stiefel manifold, i.e., when
$\matB=\matI_{d}$, reference to the last metric as the \emph{standard
metric} appears in the seminal work of Edelman, Arias and Smith~\cite{EAS98},
and this is also the metric used in the implementation of the generalized
Stiefel manifold in \noun{Manopt}~\cite{manopt}. Some of the geometric
components for working with $\stiefelB(p,d)$ equipped with the standard
metric appear in \cite[ Section 4.5]{EAS98}, while \noun{Manopt}
implements all the geometric components, but without providing a reference.

This paper is motivated by the observation that using the standard
metric in the context of Riemannian optimization with generalized
orthogonality constraints has one severe shortcoming: the computations
of some of the geometric components necessary for Riemaniann optimization
on the generalized Stiefel manifold, e.g., the Riemannian gradient
and Hessian, require taking products with the inverse of $\matB$.
Oftentimes, computing $\matB$ and its inverse is as expensive as
the direct method. In such cases there is no reason to use Riemannian
optimization as long as the standard metric, $\dotprodX UVB=\Trace{\matU^{\T}\matB\matV}$,
is used. Another issue with using the standard metric is that in some
cases it is suboptimal and using it will lead to slow convergence.

In this paper we propose to endow $\stiefelB(p,d)$ with a metric
inherited by the inner product $\dotprodM{\matU}{\matV}=\Trace{\matU^{\T}\matM_{\matX}\matV}$
on $\R^{d\times p}$ for some smooth mapping $\matX\mapsto\matM_{\matX}$
that maps a $\matX\in\stiefelB(p,d)$ to an SPD matrix $\matM_{\matX}$.
Using such a mapping is an instance of so-called \emph{Riemannian
preconditioning}~\cite{MS16}\emph{, }so we call the mapping $\matX\mapsto\matM_{\matX}$
a \emph{preconditioning scheme. }Indeed, using the metric defined
by the mapping $\matX\mapsto\matM_{\matX}$ still requires computing
$\matM_{\matX}$ in every iteration, and taking products with its
inverse, however one is free to design the mapping so that $\matM_{\matX}$
can always be cheaply decomposed. On flip side, as we discuss later,
one would like $\matM_{\matX}$ to well approximate $\matB$, or some
other matrix for which we can ensure well conditioning of the Riemannian
Hessian at the optimum. Thus in designing the mapping $\matX\mapsto\matM_{\matX}$
we have the same tradeoffs as when designing a preconditioner for
solving linear systems using a Krylov method.

In order to use Riemannian optimization with a preconditioning scheme,
one needs to implement all the necessary geometric components for
Riemannian optimization on $\stiefelB(p,d)$ endowed with the metric
defined by $\matX\mapsto\matM_{\matX}$. The majority of this paper
is devoted to developing these geometric components. We complement
these developments by considering the use of our approach on a couple
of simple theoretical examples, and on the problem of finding the
top canonical correlation between two datasets (which we explore both
theoretically and numerically).

\subsection{Related Work}

\paragraph{Riemannian Optimization.}

Riemannian optimization is an approach for solving constrained optimization
problems in which the constraints form a smooth manifold (e.g., nonlinear
differentiable equality constraints). It is based on extending classical
algorithms for unconstrained optimization on $\R^{n}$ (or any other
vector space equipped with an inner product), by generalizing the
main components needed to apply these algorithms to search spaces
that form smooth manifolds. Some early works are \cite{luenberger1972gradient,gabay1982minimizing,smith1994optimization}.
A more recent and detailed introduction can be found in \cite{AMS09}
and in \cite{boumal2020intromanifolds}.

\paragraph{Riemannian Optimization on the (Generalized) Stiefel Manifold.}

Optimization with orthogonality constraints are prevalent in many
applications across science, naturally giving rise to Riemannian optimization
on the (generalized) Stiefel manifold. Using Riemannian optimization
to solve problems with orthogonality constraints was considered in
the seminal work of Edelman et al.~\cite{EAS98}, and in particular
the components of the Stiefel manifold were developed with the \emph{standard}
(and also the \emph{canonical}) metric. Some recent works include
\cite{wen2013feasible,zhu2017riemannian,li2020efficient}, where the
Cayley transform is used to define a retraction map which leads to
more efficient algorithms. Another improved retraction computation
is proposed in \cite{sato2019cholesky}, where Sato and Aihara proposed
a Cholesky QR-based retraction on the generalized Stiefel manifold.
In \cite{kaneko2012empirical}, Kaneko et al. presented algorithms
to compute inverses of several retractions on the Stiefel manifold
in order to solve empirical arithmetic averaging problems over the
Stiefel manifold. Also, several optimization algorithms for non-smooth
optimization were developed on the Stiefel manifold such as a proximal
gradient method and a fast iterative shrinkage-thresholding algorithm
(FISTA \cite{beck2009fast}), see for example \cite{chen2018proximal,chen2019alternating,huang2019extending}.
Also in the context of this paper, a Riemannian optimization approach
for adaptive CCA on a product manifold of two generalized Stifel manifolds
was proposed in \cite{yger2012adaptive}. In addition, components
for the complex Stiefel manifold with the standard metric were developed
in several works , e.g., \cite{pechen2008control,oza2009optimization,sato2013complex,sato2014Riemannian}.
Unlike in our work, all the aforementioned works only use either the
standard or the canonical \cite{EAS98} metrics when optimizing on
the (generalized) Stiefel manifold.

\paragraph{Riemannian Preconditioning.}

In the context of Riemannian optimization, it is well-known that the
condition number of the Riemannain Hessian at the optimum is highly
indicative of the asymptotic convergence rate of Riemannian optimization
(e.g., \cite[Theorem 4.5.6, Theorem 7.4.11 and Eq. (7.50)]{AMS09}).
If the objective function is convex (in the Riemannian sense \cite[Chapter 3.2]{udriste2013convex})
then there also exist global convergence results depending on the
condition number of the Riemannian Hessian at all the points on the
manifold (e.g., \cite[Chapter 7, Theorem 4.2]{udriste2013convex}),
however these results are not applicable to optimization on the generalized
Stiefel manifold, since every continuous and convex function (in the
Riemannian sense) on the Stiefel manifold is constant.

The relation between convergence rate and condition number of the
Riemannian Hessian at the optimum motivates adjusting the metric based
on the cost or constraints, and this approach to preconditioning was
presented in several works, see e.g., \cite{ngo2012scaled,mishra2014r3mc,shi2016low,zhou2016riemannian}.
Most of the aforementioned works attempt to lower the condition number
of the Riemannian Hessian at the optimum by approximating the Euclidean
Hessian of the cost function. However, it is possible for the Riemannian
Hessian and the Euclidean Hessian to be very far from each other even
for simple examples (see Section \ref{sec:importance}). In \cite{MS16},
Mishra and Sepulchre showed that carefully selecting the metric based
on both the cost and the constraint (inspired by the Lagrangian) used
in Riemannian optimization affects convergence~\cite{MS16} of Riemannian
steepest-descent (the iterations become a version of Riemannian quasi-Newton
close to the optimum). They demonstrated this technique on a quotient
manifold (generalized Grassmann manifold) and on the fixed-rank manifold.
Unlike \cite{MS16}, we do not commit to a specific structure of the
metric, as long as it is inherited from the ambient space. Our framework
is suitable for the use of the metrics presented in \cite{MS16},
but also allows to use easier to compute metrics. Moreover, we develop
explicit components of Riemannian optimization on the generalized
Stiefel manifold with non-standard metric and consider their costs
with respect to the choice of metric (see Section \ref{sec:precond-geometry}).
This allows the use of various algorithms for smooth Riemannian optimization,
e.g, conjugate-gradient, trust-region, etc. We also motivate the choice
of metric by the condition number of the Riemannian Hessian at the
optimum.

Another similar view of Riemannian preconditioning in the sense of
Riemannian metric selection, which is specific for the Riemannian
trust-region algorithm, is to precondition the solver used to solve
the Trust-Region Subproblem~\cite{vandereycken2010riemannian}. The
aforementioned preconditioning approach generalizes the preconditioning
strategy for the unconstrained trust-region problem. Another example
of using Riemannian preconditioning for the Trust-Region Subproblem
can be found in \cite{mor2020solving}.

A different approach for preconditioning of Riemannian methods can
be found in \cite{kressner2016preconditioned} where linear systems
with tensor product structure are considered. That paper proposed
a Riemannian analogue to the preconditioned Richardson method for
Euclidean optimization based on the truncated Richardson iteration.
Similarly to Euclidean Preconditioned Richardson, in each iteration
the search direction is multiplied by an inverse of an SPD preconditioner
(and then projected to the tangent space). Another method proposed
in \cite{kressner2016preconditioned} is an approximate Riemannian
Newton method where the search direction is determined by an equation
involving an approximation to the Riemannian Hessian (known as constrained
Gauss--Newton, see e.g., \cite{bock1987randwertproblemmethoden}),
and a preconditioning term replacing a component in that equation.

\section{\label{sec:Preliminaries}Preliminaries}

\subsection{Notation and Basic Definitions}

We denote scalars using lower case Greek letters or using lower case
English letters $x,y,\dots$. Vectors are explicitly defined and also
denoted by $\x,\y,\dots$.  Matrices are denoted by $\matA,\matB,\dots$
or upper case Greek letters. Tangent vectors (of a manifold) are denoted
using lower case Greek letters with a subscript for the point on the
manifold to which they correspond (e.g., $\eta_{\x}$). Normal vectors
(of a manifold) are denoted using  lower and upper case English letters
with a subscript for the point on the manifold to which they correspond
(e.g., $\u_{\x}$). Vector fields on a manifold are denoted using
lower case Greek letters with brackets indicating the point on the
manifold to which they correspond (e.g., $\eta(\x)$). Normal vector
fields on a manifold are denoted using lower and upper case English
letters with brackets indicating the point on the manifold to which
they correspond (e.g., $\u(\x)$). We use the convention that vectors
are column-vectors.

We denote by $\dotprod{\cdot}{\cdot}_{\matC}$ the inner product with
respect to a matrix $\matC$: for vectors $\u$ and $\v$, $\dotprod{\u}{\v}_{\matC}\coloneqq\u^{\T}\matC\v$,
and for matrices $\matU$ and $\matV$, $\dotprod{\matU}{\matV}_{\matC}\coloneqq\Trace{\matU^{\T}\matC\matV}$
where $\Trace{\cdot}$ denotes the trace operator. The $s\times s$
identity matrix is denoted $\matI_{s}$. The $s\times s$ zero matrix
is denoted $\matzero_{s}$. We denote by ${\cal {\cal S}_{\text{sym}}}(p)$
and ${\cal {\cal S}_{\text{skew}}}(p)$ the set of all symmetric and
skew-symmetric matrices (respectively) in $\R^{p\times p}$.

Given a $d\times d$ matrix $\matA$ we denote by $\sym{(\matA)}\coloneqq\left(\matA+\matA^{\T}\right)/2$
and by $\skew{(\matA)}\coloneqq\left(\matA-\matA^{\T}\right)/2$ the
symmetric and skew-symmetric (respectively) components of $\matA$.
We describe a diagonal matrix using $\diag{\cdot}$ where the diagonal
components appear in the parenthesis, and similarly block diagonal
matrices are described using $\blockdiag{\cdot}$. For an SPD matrix
$\matB\in\R^{d\times d}$, we denote by $\matB^{\nicehalf}$ the unique
SPD matrix such that $\matB=\matB^{\nicehalf}\matB^{\nicehalf}$.
This matrix is obtained by keeping the same eigenvectors and taking
the square root of the eigenvalues. We denote the inverse of $\matB^{\nicehalf}$
by $\matB^{-\nicehalf}$.

Let $\matA$ be a symmetric $d\times d$ matrix. We use $\lambda_{1}(\matA)\geq\lambda_{2}(\matA)\geq\dots\geq\lambda_{d}(\matA)$
to denote the eigenvalues of $\matA$, and use $\kappa(\matA)$ to
denote the condition number of $\mat A$, which is the ratio between
the largest and smallest eigenvalues in absolute value. Let $\matB\in\R^{d\times d}$
be another symmetric positive semi-definite matrix, and assume that
$\ker(\matB)\subseteq\ker(\matA)$. If for $\lambda\in\R$ and $\v\notin\ker(\matB)$
it holds that $\matA\v=\lambda\matB\v$ then $\lambda$ is a generalized
eigenvalue and $\v$ is a generalized eigenvector of the matrix pencil
$(\matA,\matB)$. We use the notation $\lambda_{1}(\matA,\matB)\geq\lambda_{2}(\matA,\matB)\geq\dots\geq\lambda_{\rank{\matB}}(\matA,\matB)$
to denote the generalized eigenvalues of $(\matA,\matB)$. The (generalized)
condition number $\kappa(\matA,\matB)$ of the pencil $(\matA,\matB)$
is the ratio between the largest and smallest generalized eigenvalues
in absolute value. If $\matB$ is also non-singular, that is $\matB$
is an SPD matrix, then it holds that $\kappa(\matA,\matB)=\kappa(\matB^{-\nicehalf}\matA\matB^{-\nicehalf})$.

We denote by $\stiefelB(p,d)$ the generalized Stiefel manifold defined
by  (\ref{eq:general_Stiefel_def}). $\stiefelB(p,d)$ is a submanifold
of $\R^{d\times p}$. Given a function or vector field defined on
$\stiefelB(p,d)$, we use a bar decorator to denote a smooth extension
of that object to the entire $\R^{d\times p}$, either by committing
to a specific extension, or making sure that any statement made afterwards
holds for any such smooth extension. For example, given a smooth objective
function $f:\stiefelB(p,d)\to\R$, we use $\bar{f}:\R^{d\times p}\to\R$
to denote a smooth real-valued function defined on $\R^{d\times p}$
whose restriction to $\stiefelB(p,d)$ is $f$.

For $p=1$, we denote by $\elp$ the $d-1$ dimensional ellipsoid
defined by 
\[
\elp\coloneqq\left\{ \x\in\R^{d}\,:\,\x^{\T}\matB\x=1\right\} \,.
\]
In the special case $\matB=\matI_{d}$, we denote by $\stiefel(p,n)$
the Stiefel manifold defined by 
\[
\stiefel(p,d)\coloneqq\left\{ \matX\in\R^{d\times p}:\ \matX^{\T}\matX=\matI_{p}\right\} \,.
\]

Given an SPD matrix $\matB\in\R^{d\times d}$, we say that a decomposition
$\matA=\matQ\matR$ of $\matA\in\R^{d\times p}$ where $\matQ\in\R^{d\times p}$
and $\matR\in\R^{p\times p}$ is a thin $\matB$-QR decomposition
of $\matA$ if $\matQ\in\stiefelB(p,d)$ and $\mat R\in\R^{p\times p}$
is an upper triangular matrix. Note that the standard thin QR decomposition
(\cite[Chapter 5 and Lecture 7]{Golub:1996:MC:248979,trefethen1997numerical})
is a thin $\matI_{d}$-QR decomposition. Moreover, the thin $\matB$-QR
decomposition can be obtained using a standard thin QR decomposition
of the matrix $\matB^{\nicehalf}\matA$. Indeed, if $\matB^{\nicehalf}\matA=\matQ\matR$
with $\matQ\in\stiefel(p,d)$ then $\matA=(\matB^{-\nicehalf}\matQ)\matR$
and $\matB^{-\nicehalf}\matQ\in\stiefelB(p,d)$. The thin QR decomposition
is unique if $\matA$ is full rank and we require $\matR$ to have
strictly positive diagonal elements (\cite[Theorem 5.2.2 and Theorem 7.2]{Golub:1996:MC:248979,trefethen1997numerical}).
Consequently, we also have that the thin $\matB$-QR decomposition
is unique if $\matA$ is full rank and we require $\matR$ to have
strictly positive diagonal elements. In that case, we denote by $\qfm{\matA}{\matB}$
the unique $\matQ$ factor of the thin $\matB$-QR decomposition.
For the thin $\matI_{d}$-QR decomposition we abbreviate $\qf{\matA}\coloneqq\qfm{\matA}{\matI_{d}}$.
Using this notation we have the following relation \cite{sato2019cholesky}:
\[
\qfm{\matA}{\matB}=\matB^{-\nicehalf}\qf{\matB^{\nicehalf}\matA}.
\]

\subsection{Riemannian Optimization}

In this section we recall some basic definitions of Riemannian optimization,
and establish corresponding notations. A Riemannian manifold ${\cal M}$
is a real differentiable manifold ${\cal M}$ with a smoothly varying
inner product $g_{\x}$ on tangent spaces $T_{\x}{\cal M}$ (where
$\x\in{\cal M}$). A Riemannian manifold $({\cal M},g)$ is a Riemannian
submanifold of another Riemannian manifold $(\bar{{\cal M}},\bar{g})$,
if ${\cal M}$ is a submanifold of ${\cal \bar{{\cal M}}}$ and it
inherits the metric in a natural way: $g_{\x}(\eta_{\x},)=\bar{g}_{\x}(\eta_{\x},)$
for $\eta_{\x},\in T_{\x}{\cal M}$ where in the right-side $\eta_{\x}$
and $ $ are viewed as elements in $T_{\x}\bar{{\cal M}}$ (this is
possible since ${\cal M}$ is a submanifold of $\bar{{\cal M}}$).
The former notion is useful when the search space is embedded in a
larger space and the objective function is given in the coordinates
of the embedding space.

The fundamental idea in Riemannian optimization algorithms is to locally
approximate the constraint manifold by its tangent space at every
iteration. Each iterate on the tangent space minimizes some model
of the cost function, and then (possibly after several steps on the
tangent space, e.g., \cite{lezcano2019trivializations,criscitiello2020accelerated})
translates to the manifold using the retraction mapping $R_{\x}:T_{\x}{\cal M}\to{\cal M}$
\cite[Section 4.1]{AMS09}. Manipulation of tangent vectors from different
tangent spaces is performed via the vector transport $\tau_{\eta_{\x}}\in T_{R_{\x}(\eta_{\x})}{\cal M}$
\cite[Section 8.1]{AMS09}. In particular, the exponential mapping
\cite[Section 5.4]{AMS09} and parallel translation \cite[Section 5.4]{AMS09}
are examples of retraction and vector transport based on geodesics.
Note that computing them is costly, since it requires solving a system
of differential equations which might be solvable only numerically.

The notions of Riemannian gradient and Riemannian Hessian \cite[Section 3.6 and 5.5]{AMS09}
generalize the corresponding concepts from the Euclidean setting.
The Riemannian gradient is used for finding critical points, while
the Riemannian Hessian classifies them. Moreover, (asymptotic) convergence
of Riemannian methods is governed by the condition number of the Riemannian
Hessian at the optimal point.

For a smooth (objective) function defined on the manifold, $f:{\cal M}\to\R$,
denote the Riemannian gradient and Riemannian Hessian at $\x\in{\cal M}$
by $\grad{f(\x)\in T_{\x}{\cal M}}$ and $\hess{f(\x)}:T_{\x}{\cal M}\to T_{\x}{\cal M}$
respectively. Roughly speaking, the Levi-Civita (Riemannian connection)
$\nabla$ of $({\cal M},g)$ generalizes the notion of directional
derivative of vector fields.

With these components, various optimization algorithms are naturally
generalized from the Euclidean setting to the Riemannian setting (e.g.,
\cite{AMS09,boumal2020intromanifolds} for an extensive overview of
smooth techniques, and \cite{bento2017iteration,ferreira2002proximal,chen2018proximal,chen2019alternating,huang2019extending}
for some examples of non-smooth algorithms). For example, a variant
of Riemannian gradient descent is
\[
\x_{k+1}=R_{\x_{k}}(-\alpha_{k}\grad{f(\x_{k})})
\]
where $\alpha_{k}$ is the step size (possibly chosen by the Armijo's
backtracking procedure; see \cite[Algorithm 1]{AMS09}).

\section{\label{sec:precond-geometry}Preconditioned Geometric Components
for the Generalized Stiefel Manifold}

In this section we describe the necessary geometric components required
for Riemannian optimization on $\stiefelB(p,d)$ with a preconditioned
Riemannian metric. In the following, $\matB\in\R^{d\times d}$ is
an SPD matrix and we treat $\stiefelB(p,d)$ as an embedded submanifold
of $\R^{d\times p}$. Unlike previous articles in the literature,
we allow for a wider array of Riemannain metrics on $\stiefelB(p,d)$,
i.e., the metric is defined via a preconditioning scheme $\matX\mapsto\matM_{\matX}$.
We refer to the components we develop as \emph{preconditioned geometric
components} for $\stiefelB(p,d)$. In the following, we refer to $\R^{d\times p}$
as the ambient space. It is important to stress that all our formulas
are given in ambient space coordinates, and not in some local coordinates
of the manifold $\stiefelB(p,d)$.

In terms of computational costs of the geometric components, we remark
that an important feature of the components we develop is that they
access $\matB$ only via matrix-matrix products. In particular, the
formulas do not involve $\matB^{-1}$ but rather $\matM_{\matX}^{-1}$.
In quite a few problems involving generalized orthogonality constraints,
the matrix $\matB$ is given in a (semi-)implicit form, and it is
desirable to avoid computing it. In many applications, just forming
$\matB$ is as expensive as using a direct method. However, to use
the preconditioned geometric components one can avoid computing $\matB$.

\subsection{\label{subsec:metric-independent}Metric Independent Notions}

We first describe notions that are independent of the metric. This
is not our main contribution, as most of the following definitions
and formulas are well known (see e.g., \cite{yger2012adaptive,sato2019cholesky,zhu2017riemannian,kaneko2012empirical,li2020efficient,huang2019extending});
we include these definitions and formulas, and their derivations (which
appear in the appendix), for completeness. Additionally, most of the
formulas in this section can be derived via the known components of
the Stiefel manifold \cite{AMS09} via the change of variables $\hat{\matX}=\matB^{\nicehalf}\matX$.

We remark that the formulas for the inverses of various retractions
do not appear in the previous literature. However, for the most part
they too are simple generalizations of formulas for the Stiefel manifold,
which are derived in~\cite{EAS98,AMS09,zhu2017riemannian,li2020efficient}.
The inverse retraction is used in several recent algorithms proposed
in the literature: Riemannian CG with inverse retractions \cite{zhu2020riemannian},
Riemannian FISTA \cite{huang2019extending}, and empirical arithmetic
averaging over the Stiefel manifold \cite{kaneko2012empirical}.

The tangent space of $\text{St}_{\matB}(p,d)$ at $\matX\in\text{St}_{\matB}(p,d)$,
viewed as a subspace of $T_{\matX}\R^{d\times p}\simeq\R^{d\times p}$,
is 
\begin{equation}
T_{\matX}\stiefel_{\matB}(p,d)=\left\{ \matZ\in\R^{d\times p}\,:\,\matZ^{\T}\matB\matX+\matX^{\T}\matB\matZ=0_{p}\right\} .\label{eq:tangent_st_B}
\end{equation}
To explain ~(\ref{eq:tangent_st_B}), note that $\text{St}_{\matB}(p,d)$
is the kernel of $F(\matX)=\matX^{\T}\matB\matX-\matI_{p}$ which
is a submersion \cite[Proposition 3.3.3]{AMS09} (see further details
in Appendix~\ref{subsec:Metric-Independent-Notions-1}). $F$ is
a symmetric matrix valued function, so the dimension of the tangent
space (and, as such, the manifold itself) is $dp-p(p+1)/2$.

Obviously, if $\matZ\in T_{\matX}\stiefel_{\matB}(p,d)$ then the
matrix $\matX^{\T}\matB\matZ$ is skew-symmetric. Thus, a different
characterization of $T_{\matX}\stiefel_{\matB}(p,d)$ is as a decomposition
of every tangent vector into a sum of a product of a skew-symmetric
matrix with $\matX$, and a term whose columns are $\matB$-orthogonal
to the columns of $\matX$:
\begin{equation}
T_{\matX}\stiefel_{\matB}(p,d)=\left\{ \matZ=\matX\Omega+\matX_{\matB\perp}K\in\R^{d\times p}\,:\,\Omega\in{\cal {\cal S}_{\text{skew}}}(p),\ K\in\R^{(d-p)\times p}\right\} ,\label{eq:tangentstiefel_st_B_2}
\end{equation}
where $\Omega$ is a skew-symmetric matrix (i.e., $\Omega^{\T}=-\Omega$),
$K$ is arbitrary, and $\matX_{\matB\perp}\in\R^{d\times(d-p)}$ satisfies
that its columns are an orthonormal basis for the orthogonal complement
of the column space of $\matX$ with respect to the matrix $\matB$,
i.e., $\matX_{\matB\perp}^{\T}\matB\matX_{\matB\perp}=\matI_{d-p}$,
and $\matX_{\matB\perp}^{\T}\matB\matX=0_{(d-p)\times p}$.

There are several known retraction mappings suitable for the generalized
Stiefel manifold. We mention three of them (not including the exponential
map which is presented later and is also a retraction). The first
retraction mapping is based on the polar decomposition of the matrix
$\matX+\xi_{\matX}$ with respect to the inner product defined by
the matrix $\matB$ (i.e., decomposition of a matrix $\matA=\matQ\matP$
where $\matQ\in\text{St}_{\matB}(p,d)$ and $\matP$ is an SPD matrix
of the size $p\times p$; one such decomposition is $\matP=\left(\matA^{\T}\matB\matA\right)^{\nicehalf}$
and $\matQ=\matA\left(\matA^{\T}\matB\matA\right)^{-\nicehalf}$):
\begin{equation}
R_{\matX}^{\text{polar}}(\xi_{\matX})\coloneqq(\matX+\xi_{\matX})(I_{p}+\xi_{\matX}^{\T}B\xi_{\matX})^{-\nicehalf},\label{eq:stiefelpolarretraction}
\end{equation}
where $\xi_{\matX}\in T_{\matX}\text{St}_{\matB}(p,d)$. As for the
arithmetic complexity, once $\matB\xi_{\matX}$ has been computed,
we can compute $R_{\matX}^{\text{polar}}(\xi_{\matX})$ in $O(dp^{2}$)
operations.

Given $\matY\in\text{St}_{\matB}(p,d)$ close enough to $\matX$,
the inverse of the polar retraction is
\begin{equation}
R_{\matX}^{\text{polar}^{-1}}(\matY)\coloneqq\matY\matZ-\matX,\label{eq:stiefelpolarretractioninv}
\end{equation}
where $\matZ$ is the unique SPD solution of the following Lyapunov
equation
\begin{equation}
2\matI_{p}=\matX^{\T}\matB\matY\matZ+\matZ\matY^{\T}\matB\matX.\label{eq:polar-ret-lyapunov}
\end{equation}
Thus, once $\matB\matX$ is computed we can compute $R_{\matX}^{\text{polar}^{-1}}(\matY)$
using $O(dp^{2})$ operations. The expression for the inverse retraction
in ~(\ref{eq:stiefelpolarretractioninv}) is valid when the Lyapunov
~(\ref{eq:polar-ret-lyapunov}) has a unique SPD solution. If $\matY=R_{\matX}^{\text{polar}}(\xi_{\matX})$
for some $\xi_{\matX}$, then ~(\ref{eq:stiefelpolarretractioninv})
has an SPD solution $\matZ=(\matI_{p}+\xi_{\matX}^{\T}\matB\xi_{\matX})^{-1/2}$
(see Appendix \ref{subsec:Metric-Independent-Notions-1} for more
details). Let us now consider when ~(\ref{eq:stiefelpolarretractioninv})
has a unique solution. It has a unique solution if and only if $\matX^{\T}\matB\matY$
and $-\matY^{\T}\matB\matX$ do not share any eigenvalue \cite[Theorem 2.4.4.1]{horn2012matrix}.
Both $\matX^{\T}\matB\matY$ and $-\matY^{\T}\matB\matX$ are invertible
since they are products of full rank matrices, thus all eigenvalues
are not equal to zero. Next recall that $\matX^{\T}\matB\matX=\matI_{p}$,
and that eigenvalues of a matrix are a continuous function of the
matrix. Using the Bauer--Fike theorem \cite{bauer1960norms} for
a small enough perturbation of the matrix $\matX^{\T}\matB\matX$,
i.e., $\matX^{\T}\matB\matX+\delta\matX^{\T}\matB\matX=\matX^{\T}\matB\matY$,
the eigenvalues of $\matX^{\T}\matB\matY$ do not differ from the
eigenvalues of $\matX^{\T}\matB\matX$ more than the norm of the perturbation.
Thus, the real part of the eigenvalues of $\matX^{\T}\matB\matY$
remains strictly positive, leading to $\matX^{\T}\matB\matY$ and
$-\matY^{\T}\matB\matX$ not sharing any eigenvalue. The validity
of ~(\ref{eq:stiefelpolarretractioninv}) is the intersection of
the image of $R_{\matX}^{\text{polar}}(\cdot)$ and a neighborhood
of $\matX$ in which ~(\ref{eq:polar-ret-lyapunov}) has a unique
solution.

The second retraction mapping is based on the QR decomposition with
respect to the matrix $\matB$:
\begin{equation}
R_{\matX}^{\text{QR}}(\xi_{\matX})\coloneqq\qfm{\matX+\xi_{\matX}}{\matB}=\matB^{-\nicehalf}\qf{\matB^{\nicehalf}\left(\matX+\xi_{\matX}\right)},\label{eq:stiefelqrretraction}
\end{equation}
where $\xi_{\matX}\in T_{\matX}\text{St}_{\matB}(p,d)$~\cite{sato2019cholesky}.
One can show that if $\matR^{\T}\matR=(\matX+\xi_{\matX})^{\T}\matB(\matX+\xi_{\matX})$
is a Cholesky decomposition then $\qfm{\matX+\xi_{\matX}}{\matB}=(\matX+\xi_{\matX})\matR^{-1}$,
so once $\matB(\matX+\xi_{\matX})$ has been computed we can compute
$R_{\matX}^{\text{QR}}(\xi_{\matX})$ using $O(dp^{2})$ operations~\cite{sato2019cholesky}.

Given $\matY\in\text{St}_{\matB}(p,d)$ close enough to $\matX$ the
inverse of the QR-based retraction is
\begin{equation}
R_{\matX}^{\text{QR}^{-1}}(\matY)\coloneqq\matY\matR-\matX,\label{eq:stiefelqrretractioninv}
\end{equation}
where $\matR$ is the unique upper-triangular $p\times p$ matrix
with strictly positive elements on its main diagonal which is a solution
for the following Lyapunov-like equation:
\begin{equation}
2\matI_{p}=\matX^{\T}\matB\matY\matR+\matR^{\T}\matY^{\T}\matB\matX.\label{eq:qr-ret-lyapunov}
\end{equation}
Solving this equation takes $O(p^{4})$ operations \cite{kaneko2012empirical}.
Thus, once $\matB\matX$ is computed we can compute $R_{\matX}^{\text{QR}^{-1}}(\matY)$
using $O(p^{4}+dp^{2})$ operations. Note that this equation has a
solution for $\matY=\matX$, which is $\matR=\matI$. Then, by continuity
arguments, if $\matY$ is close enough to $\matX$, a solution exists.
To show uniqueness of solution, we use \cite[Eq. (14) and Algorithm 1]{kaneko2012empirical}.
According to Kaneko et al., using the constraint that $\matR$ is
upper-triangular we can reformulate  (\ref{eq:qr-ret-lyapunov}) as
an equivalent set of linear equations which has a unique solution
if and only if all the principal minors of $\matX^{\T}\matB\matY$
are non-singular (see Appendix (\ref{subsec:Metric-Independent-Notions-1})).
Similarly to the argument for inverse of the polar inverse retraction,
since $\matX^{\T}\matB\matX=\matI_{p}$, for a small enough perturbation,
i.e., $\matX^{\T}\matB\matX+\delta\matX^{\T}\matB\matX=\matX^{\T}\matB\matY$,
the real part of the eigenvalues of $\matX^{\T}\matB\matY$ remains
strictly positive, leading to non-singularity of $\matX^{\T}\matB\matY$.
Note that in order to have consistency, it is also required that the
diagonal elements of $\matR$ are strictly positive. Again, using
a similar continuity arguments we can achieve such a solution if $\matY$
is close enough to $\matX$. The validity of ~(\ref{eq:stiefelqrretractioninv})
is the intersection of the image of $R_{\matX}^{\text{QR}}(\cdot)$
and a neighborhood of $\matX$ in which ~(\ref{eq:qr-ret-lyapunov})
has a unique solution.

The third retraction mapping is based on the Cayley transform with
respect to the matrix $\matB$ (a generalization of the retraction
presented in \cite{kaneko2012empirical}):
\begin{equation}
R_{\matX}^{\text{Cayley}}(\xi_{\matX})\coloneqq\left(\matI_{d}-\frac{1}{2}\matW(\xi_{\matX})\right)^{-1}\left(\matI_{d}+\frac{1}{2}\matW(\xi_{\matX})\right)\matX,\label{eq:stiefelcayleyretraction}
\end{equation}
where 
\[
\matW(\xi_{\matX})\coloneqq(\matI_{d}-\frac{1}{2}\matX\matX^{\T}\matB)\xi_{\matX}\matX^{\T}\matB-\matX\xi_{\matX}^{\T}(\matI_{d}-\frac{1}{2}\matB\matX\matX^{\T})\matB.
\]
Once the multiplications with $\matB$ are computed, we need to compute
the inverse of a $d\times d$ matrix in order to find $R_{\matX}^{\text{Cayley}}(\xi_{\matX})$.
However, noticing that
\[
\matW(\xi_{\matX})=\left[\begin{array}{cc}
(\matI_{d}-\frac{1}{2}\matX\matX^{\T}\matB)\xi_{\matX} & \matX\end{array}\right]\left[\begin{array}{c}
\matX^{\T}\\
-\xi_{\matX}^{\T}(\matI_{d}-\frac{1}{2}\matX\matX^{\T}\matB)^{\T}
\end{array}\right]\matB,
\]
(i.e., a product of a $d\times2p$ matrix by a $2p\times d$ matrix),
we can use the Sherman-Morrison-Woodbury formula to only invert a
$2p\times2p$ matrix. A closed form for the inverse of this retraction
is only known when $d$ is even \cite{kaneko2012empirical}.

Similarly, there are several possible ways to compute a vector transport.
It is possible to define a metric independent vector transport, using
\cite[Equation 8.6]{AMS09} by differentiating a retraction mapping
\begin{eqnarray*}
\tau_{\eta_{\matX}}^{(\text{ind})}\xi_{X} & \coloneqq & \text{D}R_{\matX}(\eta_{\matX})[\xi_{\matX}].
\end{eqnarray*}
In Appendix~\ref{subsec:Metric-Related-Notions-1}, we derive concrete
formulas based on the polar and QR retractions, (\ref{eq:stiefelpolarretraction})
and~(\ref{eq:stiefelqrretraction}). A vector transport based on
the Cayley retraction is presented in \cite{zhu2017riemannian}. The
various vector transport have the same computational cost as computing
the corresponding retractions. Note that it is also possible to define
another vector transport that has this property by simply applying
the projection on the tangent space. However, this vector transport
is metric dependent, so we discuss it in the next subsection.

\subsection{\label{subsec:Metric-Related-Notions}Metric Related Notions}

This subsection is the main contribution of our paper. In this subsection
we derive explicit formulas for the orthogonal projection with respect
to the Riemannian metric, the Riemannian gradient and Hessian with
respect to the non-standard metric which allow the use of various
preconditioned Riemannian algorithms. Note that the formulas in this
subsection, unlike the previous one, \emph{cannot} be derived via
a change of variables $\hat{\matX}=\matB^{\nicehalf}\matX$ unless
a specific metric is used (corresponding to $\matM_{\matX}=\matB$
for all $\matX\in\stiefelB(p,d)$), since though this change of variables
makes $\hat{\matX}\in\stiefel(p,d)$, the induced metric on that manifold
is not the standard metric. Note that if indeed $\matM_{\matX}=\matB$
for all $\matX\in\stiefelB(p,d)$, then via the change of variables
$\hat{\matX}=\matB^{\nicehalf}\matX$ we have that $\hat{\matX}\in\stiefel(p,d)$
with the corresponding standard metric $\matM_{\hat{\matX}}=\matI$
for all $\hat{\matX}\in\stiefel(p,d)$. Unfortunately, the change
of variables $\hat{\matX}=\matB^{\nicehalf}\matX$ forces us to form
$\matB$ explicitly, which is prohibited in problems where computing
$\matB$ is as expensive as solving them with a direct method.

Specifically, we define a Riemannian metric on the ambient space $\R^{d\times p}$,
and this uniquely defines a metric on $\text{St}_{\matB}(p,d)$ that
makes it a Riemannian submanifold. The metric we define on $\R^{d\times p}$
is 
\[
\bar{g}_{\matX}(\bar{\xi}_{\matX},\bar{\eta}_{\matX})\coloneqq\dotprodM{\bar{\xi}_{\matX}}{\bar{\eta}_{\matX}}=\Trace{\bar{\xi}_{\matX}^{\T}\matM_{\matX}\bar{\eta}_{\matX}}
\]
where $\matX\mapsto\matM_{\matX}$ is a smooth mapping on $\R^{d\times p}$
(thus, the metric varies smoothly with $\matX$ making it a Riemannian
metric), and each $\matM_{\matX}$ is assumed to be an SPD matrix
so that we have a properly defined inner product on each tangent space,
and a Riemannian metric for $\R^{d\times p}$. Now, for any $\matX\in\text{St}_{\matB}(p,d),\ \xi_{\matX},\eta_{\matX}\in T_{\matX}\stiefel_{\matB}(p,d)$,
given in ambient space coordinates, the Riemannian metric on $\text{St}_{\matB}(p,d)$
is given by 
\begin{equation}
g_{\matX}(\xi_{\matX},\eta_{\matX})\coloneqq\dotprodM{\xi_{\matX}}{\eta_{\matX}}=\Trace{\xi_{\matX}^{\T}\matM_{\matX}\eta_{\matX}}.\label{eq:stiefelmetric}
\end{equation}
The cost of computing $g_{\matX}(\xi_{\matX},\eta_{\matX})$ is $O\left(T_{\matM}p+dp\right)$
where $T_{\matM}$ is the maximal cost (possibly after preprocessing)
of taking the product with $\matM_{\matX}$ with a vector for all
$\matX$.

The metric selection is how we propose to incorporate a preconditioner,
and so the mapping $\matX\mapsto\matM_{\matX}$ is termed a preconditioning
scheme. It should be chosen so that the Riemannian Hessian at the
optimum is well conditioned. We discuss this further in Subsection
\ref{subsec:Metric-selection-and-Hessian}. Classically, the metric
employed for the generalized Stiefel manifold corresponds to $\matM_{\matX}=\matB$
for all $\matX\in\stiefelB(p,d)$ \cite{EAS98}. In quite a few applications
this choice minimizes a-priori bounds on the condition number of the
Riemannian Hessian at the optimum~(see Subsections \ref{subsec:Metric-selection-and-Hessian}
and \ref{subsec:CCA}). However, as we shall see, various operations
required for Riemannian optimization require products with $\matM_{\matX}^{-1}$,
and in many applications this results in algorithms that are too expensive
when $\matM_{\matX}=\matB$ for some $\matX\in\stiefelB(p,d)$. In
such cases, there is a need to balance in the chosen $\matX\mapsto\matM_{\matX}$
between minimizing the condition number, and efficient products with
$\matM_{\matX}^{-1}$. This is a typical trade-off for preconditioning.

After defining the Riemannian metric we can derive the metric related
notions required for Riemannian optimization. Since $\text{St}_{\matB}(p,d)$
is an embedded submanifold of $\R^{d\times p}$, the orthogonal projection
on the tangent space with respect to the Riemannian metric is a key
component. We denote the orthogonal projection operator on $T_{\matX}\text{St}_{\matB}(p,d)$
by $\Pi_{\matX}\left(\cdot\right)$, and the orthogonal projection
operator (with respect to the metric defined by $\matX\mapsto\matM_{\matX}$)
on the normal space, $\left(T_{\matX}\text{St}_{\matB}(p,d)\right)^{\perp}$,
by $\Pi_{\matX}^{\perp}\left(\cdot\right)$.

In order to find analytic formulas for these operators, we first note
that the normal space is:
\begin{equation}
\left(T_{\matX}\stiefelB(p,d)\right)^{\perp}=\left\{ \matM_{\matX}^{-1}\matB\matX\matS\ :\ \matS\in{\cal {\cal S}_{\text{sym}}}(p)\right\} .\label{eq:normalstiefeMmetric}
\end{equation}
Indeed, recall that $\Trace{\matS^{\T}\Omega}=0$ for any symmetric
matrix $\matS$ and anti-symmetric matrix $\Omega$, thus by using
the representation in  (\ref{eq:tangentstiefel_st_B_2}) of tangent
vectors we get that any vector of the form $\matM_{\matX}^{-1}\matB\matX\matS$
where $\matS\in{\cal {\cal S}_{\text{sym}}}(p)$ is orthogonal to
the tangent space at $\matX\in\stiefelB(p,d)$. The dimension of the
normal space should be $p(p+1)/2$, thus since the set $\left\{ \matM_{\matX}^{-1}\matB\matX\matS\ :\ \matS\in{\cal {\cal S}_{\text{sym}}}(p)\right\} $
is $p(p+1)/2$ dimensional, it is indeed the normal space.

The following lemma gives a formula for the orthogonal projections
to the tangent and normal spaces.
\begin{lem}
\label{lem:The-orthogonal-projection}The orthogonal projections with
respect to $g_{\matX}(\cdot,\cdot)$ on $\left(T_{\matX}\stiefelB(p,d)\right)^{\perp}$
and on $T_{\matX}\stiefelB(p,d)$ (viewed as a subspace of $T_{\matX}\R^{d\times p}\simeq\R^{d\times p}$
and given in ambient coordinates) are:

\begin{equation}
\Pi_{\matX}^{\perp}\left(\xi_{\matX}\right)=\matM_{\matX}^{-1}\matB\matX\matS_{\xi_{\matX}}\label{eq:stiefelprojtonormal}
\end{equation}
 and

\begin{equation}
\Pi_{\matX}\left(\xi_{\matX}\right)=\left(\id_{T_{\matX}\R^{d\times p}}-\Pi_{\matX}^{\perp}\right)\left(\xi_{\matX}\right)=\xi_{\matX}-\matM_{\matX}^{-1}\matB\matX\matS_{\xi_{\matX}}\label{eq:stiefelprojtotangent}
\end{equation}
where $\xi_{\matX}\in T_{\matX}\R^{d\times p}$, $\id_{T_{\matX}\R^{d\times p}}$
denotes the identity mapping on $T_{\matX}\R^{d\times p}$, and $\matS_{\xi_{\matX}}\in\R^{p\times p}$
is the unique solution of the following Sylvester equation
\[
\left(\matX^{\T}\matB\matM_{\matX}^{-1}\matB\matX\right)\matS_{\xi_{\matX}}+\matS_{\xi_{\matX}}\left(\matX^{\T}\matB\matM_{\matX}^{-1}\matB\matX\right)=\matX^{\T}\matB\xi_{\matX}+\left(\matX^{\T}\matB\xi_{\matX}\right)^{\T}.
\]
The cost of computing (in ambient coordinates) $\Pi_{\matX}\left(\xi_{\matX}\right)$
for an arbitrary $\xi_{\matX}$ is $O(T_{\matB}p+T_{\matM^{-1}}p+dp^{2})$,
where $T_{\matB}$ and $T_{\matM^{-1}}$ are the cost of computing
the product of $\matB$ with a vector and the maximal cost of taking
the product with $\matM_{\matX}^{-1}$ with a vector for all $\matX\in\stiefelB(p,d)$.
\end{lem}

\begin{proof}
Note that $T_{\matX}\stiefel_{\matB}(p,d)\varoplus\left(T_{\matX}\stiefel_{\matB}(p,d)\right)^{\perp}=T_{\matX}\R^{d\times p}\simeq\R^{d\times p}$.
This implies that for any $\xi_{\matX}\in T_{\matX}\R^{d\times p}\simeq\R^{d\times p}$
there exists unique $\Omega_{\xi_{\matX}}\in{\cal {\cal S}_{\text{skew}}}(p)$,
$K_{\xi_{\matX}}\in\R^{(d-p)\times p}$ and $\matS_{\xi_{\matX}}\in{\cal {\cal S}_{\text{sym}}}(p)$
such that $\xi_{\matX}$ is decomposed to a unique component on the
tangent space of $\stiefel_{\matB}(p,d)$ and a unique component on
the normal space of $\stiefel_{\matB}(p,d)$:
\begin{equation}
\xi_{\matX}=\Pi_{\matX}\left(\xi_{\matX}\right)+\Pi_{\matX}^{\perp}\left(\xi_{\matX}\right)=\left(\matX\Omega_{\xi_{\matX}}+\matX_{\matB\perp}K_{\xi_{\matX}}\right)+\matM_{\matX}^{-1}\matB\matX\matS_{\xi_{\matX}}.\label{eq:generalvectorRnp}
\end{equation}

By left-multiplying  (\ref{eq:generalvectorRnp}) by $\matX^{\T}\matB$,
we get 
\[
\matX^{\T}\matB\xi_{\matX}=\Omega_{\xi_{\matX}}+\matX^{\T}\matB\matM_{\matX}^{-1}\matB\matX\matS_{\xi_{\matX}}\ .
\]
Summing $\matX^{\T}\matB\xi_{\matX}+\left(\matX^{\T}\matB\xi_{\matX}\right)^{\T}$,
and using the fact that $\Omega_{\xi_{\matX}}$ is skew-symmetric
so it vanishes in the sum, we get that $\matS_{\xi_{\matX}}$ solves
the following Sylvester equation (\cite[Subsection 2.4.4]{horn2012matrix}):
\begin{equation}
\matX^{\T}\matB\xi_{\matX}+\left(\matX^{\T}\matB\xi_{\matX}\right)^{\T}=\left(\matX^{\T}\matB\matM_{\matX}^{-1}\matB\matX\right)\matS_{\xi_{\matX}}+\matS_{\xi_{\matX}}\left(\matX^{\T}\matB\matM_{\matX}^{-1}\matB\matX\right)\ .\label{eq:sylvesterforproj}
\end{equation}
Indeed, according to \cite[Theorem 2.4.4.1]{horn2012matrix} there
is a unique solution to  (\ref{eq:sylvesterforproj}) for any $\matX^{\T}\matB\xi_{\matX}+\left(\matX^{\T}\matB\xi_{\matX}\right)^{\T}$,
since $\left(\matX^{\T}\matB\matM_{\matX}^{-1}\matB\matX\right)$
is positive definite ($\matX^{\T}\matB\matM_{\matX}^{-1}\matB\matX$
is a Gram matrix of $\matM_{\matX}^{-\nicehalf}\matB\matX$, which
consists of a product of three matrices, two invertible matrices $\matM_{\matX}^{-\nicehalf}$
and $\matB$, and one full-column rank matrix $\matX\in\stiefel_{\matB}(p,d)$)
and $-\left(\matX^{\T}\matB\matM_{\matX}^{-1}\matB\matX\right)$ is
negative definite, thus both matrices have no eigenvalues in common.
Solving  (\ref{eq:sylvesterforproj}) costs $O(p^{3})$ assuming we
already computed $\matX^{\T}\matB\matM_{\matX}^{-1}\matB\matX$. Furthermore,
as expected $\matS_{\xi_{\matX}}$ is symmetric since $\matS_{\xi_{\matX}}^{\T}$
again satisfies ~(\ref{eq:sylvesterforproj}), and the solution to
the equation is unique.

After obtaining $\matS_{\xi_{\matX}}$ by solving  (\ref{eq:sylvesterforproj}),
analytical expressions for the orthogonal projections on the normal
space and the tangent space are given by  (\ref{eq:stiefelprojtonormal})
and  (\ref{eq:stiefelprojtotangent}).

Note that the orthogonal projection on the normal space and the tangent
space satisfy the definition of an orthogonal projection with respect
to the inner product defined on $\R^{n\times p}$ with the matrix
$\matM_{\matX}$. Indeed, both projections satisfy the projection
property $\Pi_{\matX}^{2}\left(\cdot\right)=\Pi_{\matX}\left(\cdot\right)$
and $\left(\Pi_{\matX}^{\perp}\right)^{2}\left(\cdot\right)=\Pi_{\matX}^{\perp}\left(\cdot\right)$,
since $\matS_{\xi_{\Pi_{\matX}^{\perp}\left(\xi_{\matX}\right)}}$
and $\matS_{\xi_{\matX}}$ satisfy the same Sylvester equation. In
addition, both projections are orthogonal with respect to the inner
product defined on $\R^{n\times p}$ with the matrix $\matM_{\matX}$,
i.e.,
\begin{equation}
g_{\matX}(\Pi_{\matX}\left(\xi_{\matX}\right),\eta_{\matX})=g_{\matX}(\xi_{\matX},\Pi_{\matX}\left(\eta_{\matX}\right)),g_{\matX}(\Pi_{\matX}^{\perp}\left(\xi_{\matX}\right),\eta_{\matX})=g_{\matX}(\xi_{\matX},\Pi_{\matX}^{\perp}\left(\xi_{\matX}\right))\label{eq:projproperty2}
\end{equation}
 for all $\xi_{\matX},\eta_{\matX}\in\R^{n\times p}$, since by using
the properties of the trace operator.

The cost of computing (in ambient coordinates) $\Pi_{\matX}\left(\xi_{\matX}\right)$
for an arbitrary $\xi_{\matX}$ is $O(T_{\matB}p+T_{\matM^{-1}}p+dp^{2})$.
Indeed, after obtaining $\matS_{\xi_{\matX}}$ by solving a Sylvester
equation which costs $O(p^{3})$, we are left with taking product
of $\matB$ and $\matM_{\matX}^{-1}$ with matrices, and products
of matrices of the dimensions $p\times d$ by $d\times p$ , $d\times p$
by $p\times p$ and $p\times p$ by $p\times p$.
\end{proof}
In the special case where $\matM_{\matX}=\matB$ for all $\matX\in\stiefelB(p,d)$,
$\stiefelB(p,d)$ is isometric to $\stiefel(p,d)$ via the change
of variables $\hat{\matX}=\matB^{\nicehalf}\matX$. The orthogonal
projections on the normal space ~(\ref{eq:stiefelprojtonormal})
and on the tangent space ~(\ref{eq:stiefelprojtotangent}) are reduced
to a generalization of the orthogonal projection on the tangent space
of the Stiefel manifold \cite[Example 3.6.2]{AMS09}:
\begin{equation}
\Pi_{\matX}^{\perp}\left(\xi_{\matX}\right)=\matX\sym{\left(\matX^{\T}\matB\xi_{\matX}\right)}\label{eq:stiefelprojtonormalBmetric}
\end{equation}
 and

\begin{equation}
\Pi_{\matX}\left(\xi_{\matX}\right)=\left(\id_{T_{\matX}\R^{d\times p}}-\Pi_{\matX}^{\perp}\right)\left(\xi_{\matX}\right)=\left(\matI_{d}-\matX\matX^{\T}\matB\right)\xi_{\matX}+\matX\skew{\left(\matX^{\T}\matB\xi_{\matX}\right)}.\label{eq:stiefelprojtotangentBmetric}
\end{equation}
In such case, the cost of computing (in ambient coordinates) $\Pi_{\matX}\left(\xi_{\matX}\right)$
for an arbitrary $\xi_{\matX}$ is $O(T_{\matB}p+dp^{2})$. The cost
is evident from the formulas once we observe that none of the operations
require forming $\matB$, but instead require taking product of $\matB$
with a matrix of $p$ columns.

Using the orthogonal projection we can also propose a simple metric
dependent vector transport using the vector transport definition on
Riemannian submanifolds \cite[Subsection 8.1.3]{AMS09}:
\begin{eqnarray}
\tau_{\eta_{\matX}}^{(\text{dep})}\xi_{X} & \coloneqq & \Pi_{R_{\matX}(\eta_{\matX})}\left(\xi_{\matX}\right),\label{eq:vectransstiefeldep}
\end{eqnarray}
where $R_{\matX}(\cdot)$ is a retraction mapping of our choice, e.g.,
~(\ref{eq:stiefelpolarretraction}), ~(\ref{eq:stiefelqrretraction})
or  (\ref{eq:stiefelcayleyretraction}).

Let $f:\stiefelB(p,d)\to\R$ be a smooth function, and let $\bar{f}$
be a smooth extension of $f$ to $\R^{d\times p}$ (typically, $f$
is given in ambient coordinates, thereby making the extension $\bar{f}$
natural). We now develop first and second order Riemannian components
for $f$. The Riemannian gradient is an element of the tangent space,
and to derive an analytic formula for it we use \cite[Eq. 3.37]{AMS09}:
the Riemannian gradient can be computed by computing the Riemannian
gradient in $\R^{d\times p}$ of $\bar{f}$, and orthogonally projecting
it with respect to the Riemannnian metric to the tangent space of
$\stiefelB(p,d)$ using the orthogonal projection on the tangent space,
$\Pi_{\matX}\left(\cdot\right)$. In short, $\grad{f(\matX)=\Pi_{\matX}\left(\grad{\bar{f}(\matX)}\right)}$.
First, we consider $\grad{\bar{f}(\matX)}$. Note that it is not the
Euclidean gradient $\nabla\bar{f}(\matX)$, even though $\bar{f}$
is defined on $\R^{d\times p}$. The reason is that $\bar{f}$ is
defined on a $\R^{d\times p}$ endowed with a non-standard inner product.
According to \cite[Eq. 3.31]{AMS09}, we have 
\[
\Trace{\grad{\bar{f}(\matX)}^{\T}\matM_{\matX}\xi_{\matX}}=g_{\matX}(\grad{\bar{f}(\matX)},\xi_{\matX})=\text{D}\bar{f}(\matX)[\xi_{\matX}]=\Trace{\nabla\bar{f}(\matX)^{\T}\xi_{\matX}}
\]
for every $\bar{\xi}_{\matX}\in T_{\matX}\R^{d\times p}$ (in the
above, $\text{D}f(\matX)$ denotes the (Frechet) differential of $f$
at $\matX$), so $\grad{\bar{f}(\matX)}=\matM_{\matX}^{-1}\nabla\bar{f}(\matX)$.
Thus, we have 
\begin{equation}
\grad{f(\matX)}=\Pi_{\matX}\left(\matM_{\matX}^{-1}\nabla\bar{f}(\matX)\right)\,.\label{eq:Rgrad_st_B}
\end{equation}
The cost of computing the Riemannian gradient given the Euclidean
gradient of $\bar{f}$ is the cost of computing the orthogonal projection
on the tangent space, and taking the product of $\nabla\bar{f}(\matX)$
and $\matM_{\matX}^{-1}$.

The components developed so far, allow the application of any first
order Riemannian optimization algorithm, e.g., Riemannian gradient
and Riemannian conjugate-gradient. In order to apply second-order
methods, e.g., Riemannian Newton and Riemannian trust-region, the
Riemannian Hessian must also be derived. An expression for the Riemannian
Hessian is also useful for reasoning on the convergence rate by examining
the condition number of the Hessian at the optimum. However, any expression
for the Riemannian Hessian must depend on the specifics of the mapping
of $\matX$ to $\matM_{\matX}$. Thus, we focus on the simpler case
where $\matM_{\matX}=\matM$, i.e. $\matM_{\matX}$ is constant for
all $\matX\in\stiefelB(p,d)$). This is a reasonable choice for a
preconditioning metric since it still allows the use of different
cheap-to-invert constant approximations of $\matB$ (see Subsection
\ref{subsec:CCA} for an example).

Recall that in \cite[ Proposition 5.5.6]{AMS09}, it is shown that
at a critical point $\matX^{\star}$, i.e. $\grad{f(\matX^{\star})}=0$,
the Riemannian Hessian equals to the Riemannian Hessian of a composition
of the cost function with a retraction map (known in the literature
as the \emph{pullback} function). The pullback function is a function
from the tangent space which is a Euclidean space to $\R$, thus its
Riemannian Hessian is the Euclidean Hessian. In addition, retraction
maps typically do not depend on the choice of the Riemannian metric.
Therefore, the Euclidean Hessian of the pullback function only depends
on the Riemannian metric at a critical point through the directional
derivative of Riemannian gradient of the pullback function on the
tangent space at the critical point. Thus, the formula we derive for
the Riemannian Hessian in ambient coordinates is valid at a critical
point $\matX^{\star}$ when using any preconditioning scheme $\matX\mapsto\matM_{\matX}$
as well if we set $\matM=\matM_{\matX^{\star}}$. This property allows
the analysis of the condition number of the Riemannian Hessian at
the critical points with a preconditioning scheme $\matX\mapsto\matM_{\matX}$,
giving indication for the asymptotic convergence of Riemannian optimization
algorithms (e.g., \cite[Theorem 4.5.6, Theorem 7.4.11 and Eq. (7.50)]{AMS09}).

The Riemannian Hessian of $f$ at a point on the manifold is a linear
transformation from the tangent space to itself. When $\matM_{\matX}=\matM$
for all $\matX\in\stiefelB(p,d)$, we can compute the result of applying
the Riemannian Hessian to a tangent vector in ambient coordinates
via the formula~\cite{absil2013extrinsic}: 
\begin{equation}
\hess{f(\matX)[\eta_{\matX}]}=\Pi_{\matX}(\matM^{-1}\nabla^{2}\bar{f}(\matX)\eta_{\matX})+W_{\matX}(\eta_{\matX},\Pi_{\matX}^{\perp}(\matM^{-1}\nabla\bar{f}(\matX)))\label{eq:StiefelHessian}
\end{equation}
where $\nabla^{2}\bar{f}(\matX)$ is the Euclidean Hessian of $\bar{f}$
and $W_{\matX}$ is the Weingarten map on $\stiefelB(p,d)$. The Weingarten
map is an operator that takes as arguments a tangent vector $\eta_{\matX}\in T_{\matX}\stiefelB(p,d)$
and a normal vector $\matU_{\matX}\in\left(T_{\matX}\stiefelB(p,d)\right)^{\perp}$
and returns a tangent vector. An analytic formula for Weingarten map
on $\stiefelB(p,d)$, in ambient coordinates, is
\[
W_{\matX}\left(\eta_{\matX},\matU_{\matX}\right)=-\Pi_{\matX}\left(\matM^{-1}\matB\eta_{\matX}\left(\matX^{\T}\matM\matU_{\matX}\right)\right).
\]
The derivation of ~(\ref{eq:StiefelHessian}) is based on Lemma \ref{lem:modificationoftheorem}.
The complete derivation of the Riemannian connection, the Weingarten
map, and the Riemannian Hessian appears in Appendix \ref{subsec:Metric-Related-Notions-1}.
Based on these formulas, we have the following formula for the Riemannian
Hessian when $\matM_{\matX}=\matM$ for all $\matX$:
\begin{equation}
\hess{f(\matX)[\eta_{\matX}]}=\Pi_{\matX}\left(\matM^{-1}\nabla^{2}\bar{f}(\matX)\eta_{\matX}-\matM^{-1}\matB\eta_{\matX}\left(\matX^{\T}\nabla\bar{f}(\matX)-\matX^{\T}\matM\grad f(\matX)\right)\right).\label{eq:StiefelHessian-explicite}
\end{equation}
The cost of applying the Riemannian Hessian to a tangent vector given
the Euclidean Hessian of $\bar{f}$ is the cost of computing the orthogonal
projection on the tangent space, and taking the products with $\matB$,
$\matM$ and $\matM^{-1}$.

\paragraph{Exponential Map.}

An important metric related retraction map on a Riemannian manifold
is the exponential mapping. According to \cite[Proposition 5.4.1]{AMS09},
the exponential map induced by the Riemannian connection defined on
the manifold is a retraction map, termed the exponential retraction.
In particular, the exponential map is based on moving on geodesic
curves in the direction of a tangent vector. In the derivation of
the exponential map we assume $\matM_{\matX}=\matM$ for all $\matX\in\stiefelB(p,d)$.

First, let us recall the definition of a geodesic curve. A geodesic
$\gamma(t)$ on a manifold ${\cal M}$ endowed with a Riemannian connection
$\nabla$ is a curve with zero acceleration
\[
\frac{\text{D}^{2}}{\text{dt}^{2}}\gamma(t)=0,
\]
for all $t$ in the domain of $\gamma(t)$, where $\frac{\text{D}^{2}}{\text{dt}^{2}}\gamma(t)=\frac{\text{D}}{\text{dt}}\dot{\gamma}$~\cite[Section 5.4]{AMS09}
.

On the generalized Stiefel manifold, the function $\xi_{\matX}\longmapsto\frac{\text{D}}{\text{dt}}\xi_{\matX}$
from the set of all (smooth) vector fields on $\stiefelB(p,d)$ to
itself is $\frac{\text{D}}{\text{dt}}\left(\cdot\right)\coloneqq\Pi_{\matX(t)}\left(\frac{\text{d}}{\text{dt}}\left(\cdot\right)\right)$.
For every $\xi_{x}\in T_{\x}{\cal M}$, there exists an interval $I$
about $0$ and a unique geodesic $\gamma(t;\x,\xi):\,I\to{\cal M}$
such that $\gamma(0)=\x$ and $\dot{\gamma}(0)=$. Moreover, we have
the homogeneity property $\gamma(t;\x,a)=\gamma(at;\x,)$. The mapping
\[
\text{Exp}_{\x}:\:T_{\x}{\cal M}\to{\cal M}:\:\mapsto\text{Exp}_{\x}=\gamma(1;\x,),
\]
is called the exponential map at $\x$~\cite[Section 5.4]{AMS09}.

To find the exponential map on the Stiefel manifold $\stiefelB(p,d)$,
we need to find the geodesic given $\matX=\gamma(0)\in\stiefelB(p,d)$
and $\xi_{\matX}=\dot{\gamma}(0)\in T_{\matX}\stiefelB(p,d)$, i.e.,
we need to solve the differential equation
\begin{eqnarray}
\frac{\text{D}^{2}}{\text{dt}^{2}}\gamma(t) & = & 0\nonumber \\
\Pi_{\gamma(t)}\left(\frac{\text{d}}{\text{dt}}\left[\frac{\text{d}}{\text{dt}}\left(\gamma(t)\right)\right]\right) & = & 0\label{eq:exp_ode}\\
\Pi_{\gamma(t)}\left(\ddot{\gamma}(t)\right) & = & 0\nonumber \\
\ddot{\gamma}(t) & = & \matM^{-1}\matB\gamma(t)\matS_{\ddot{\gamma}(t)}\ .\nonumber 
\end{eqnarray}
where the matrix $\matS_{\ddot{\gamma}(t)}$ satisfies the following
Sylvester equation 
\[
\gamma(t)^{\T}\matB\ddot{\gamma}(t)+\left(\gamma(t)^{\T}\matB\ddot{\gamma}(t)\right)^{\T}=\left(\gamma(t)^{\T}\matB\matM^{-1}\matB\gamma(t)\right)\matS_{\ddot{\gamma}(t)}+\matS_{\ddot{\gamma}(t)}\left(\gamma(t)^{\T}\matB\matM^{-1}\matB\gamma(t)\right).
\]
Note that we can replace $\gamma(t)^{\T}\matB\ddot{\gamma}(t)+\left(\gamma(t)^{\T}\matB\ddot{\gamma}(t)\right)^{\T}$
by \textbf{$-2\dot{\gamma}(t)^{\T}\matB\dot{\gamma}(t)$ }since $\gamma(t)^{\T}\matB\gamma(t)=\matI_{p}$
when $\gamma(t)\in\stiefelB(p,d)$ and by differentiating two times
with respect to $t$ we get the equality.

Thus, in order to compute the exponential map, we simply need to solve
~(\ref{eq:exp_ode}). Unfortunately, in the general case we are unaware
of any analytical solution, and so the equation needs to be solved
numerically. However, in the special case where $\matM_{\matX}=\matB$
for all $\matX\in\stiefelB(p,d)$ such that $\stiefelB(p,d)$ is isometric
to $\stiefel(p,d)$ via the change of variables $\hat{\matX}=\matB^{\nicehalf}\matX$,
the equation can be solved analytically in a manner similar to \cite[Equation 5.26]{AMS09}.
For $\matM_{\matX}=\matB$ , the equation for the geodesic is reduced
to 
\[
\ddot{\gamma}(t)=-\gamma(t)\left(\dot{\gamma}(t)^{\T}\matB\dot{\gamma}(t)\right)\ .
\]
 We perform a small modification of the calculations given in \cite[Subsection 2.2.2]{EAS98}
(also developed by Ross Lippert). Denote
\[
\matC\coloneqq\gamma(t)^{\T}\matB\gamma(t),\ \matA\coloneqq\gamma(t)^{\T}\matB\dot{\gamma}(t),\ \matS\coloneqq\dot{\gamma}(t)^{\T}\matB\dot{\gamma}(t).
\]
By differentiating $\matC,\matA,\matS$ we get the following equations:
\begin{eqnarray*}
\dot{\matC} & = & \matA+\matA^{\T}\ ,\\
\dot{\matA} & = & \matS+\gamma(t)^{\T}\matB\ddot{\gamma}(t)=\matS-\matC\matS\ ,\\
\dot{\matS} & = & \ddot{\gamma}(t)^{\T}\matB\dot{\gamma}(t)+\dot{\gamma}(t)^{\T}\matB\ddot{\gamma}(t)=-\left[\matS\matA+\matA^{\T}\matS\right]\ .
\end{eqnarray*}
Recall that since $\gamma(t)\in\stiefelB(p,d)$ we get that $\matC=\matI_{p}$.
Thus, $\dot{\matC}=\mat 0_{p}$ so $\matA=-\matA^{\T}$ , i.e., $\matA$
is skew-symmetric. Moreover, $\dot{\matA}=\mat 0_{p}$ so that $\matA(t)=\matA(0)$.
In addition, the last equation can be rewritten as 
\begin{eqnarray*}
\dot{\matS} & = & \matA\matS-\matS\matA\ ,
\end{eqnarray*}
and it has a closed form (see \cite[Theorem 9.2]{bhatia1997and} for
a constant matrix $\matA$) solution of the form
\begin{eqnarray*}
\matS(t) & = & e^{\matA t}\matS(0)e^{-\matA t}\ .
\end{eqnarray*}
Finally, we can use the following equation
\[
\frac{\text{d}}{\text{dt}}\left[\gamma(t)e^{\matA t},\ \dot{\gamma}(t)e^{\matA t}\right]=\left[\gamma(t)e^{\matA t},\ \dot{\gamma}(t)e^{\matA t}\right]\left(\begin{array}{cc}
\matA & -\matS(0)\\
\matI_{p} & \matA
\end{array}\right)\ ,
\]
to find a closed form for the geodesic curve
\begin{equation}
\gamma(t)=\left[\matX,\ \xi\right]\exp\left(t\left(\begin{array}{cc}
\matA & -\matS(0)\\
\matI_{p} & \matA
\end{array}\right)\right)\left[\begin{array}{c}
\matI_{p}\\
\mat 0_{p}
\end{array}\right]e^{-\matA t}\ .\label{eq:geodesiconStiefelMeqB}
\end{equation}
Substituting $t=1$ into ~(\ref{eq:geodesiconStiefelMeqB}) gives
us the exponential mapping $\text{Exp}_{\matX}\xi_{\matX}$.

\subsection{\label{subsec:comp-costs}Computational Costs}

Table~\ref{tab:costs} summarizes the computational costs, measured
in terms of arithmetic operations, of computing the Riemannian components
on the generalized Stiefel manifold described in Subsections \ref{subsec:metric-independent}
and \ref{subsec:Metric-Related-Notions}. Note that all the costs
are for operations in ambient coordinates. In the table, we denote
by $T_{\matC}$ the cost of computing the product of $\matC$ with
a vector (potentially, after preprocessing $\matC$), for some matrix
$\matC$. Specifically, we use $T_{\matB},T_{\matB^{-\nicehalf}},T_{\matB^{\nicehalf}},,T_{\matM}$
and $T_{\matM^{-1}}$. In particular, $T_{\matM}$ and $T_{\matM^{-1}}$
denote the maximal cost (over $\matX\in\stiefelB(p,d)$) of taking
the product of $\matM_{\matX}$ and $\matM_{\matX}^{-1}$ (respectively)
with a vector. Also, we denote by $T_{\nabla\bar{f}}$ and by $T_{\nabla^{2}\bar{f}}$
the cost of computing the Euclidean gradient and the cost of applying
the Euclidean Hessian to a tangent vector.

Note that compared to the standard metric on $\stiefelB(p,d)$ (i.e.,
$\matM_{\matX}=\matB$ for all $\matX$), we replace products with
$\matB^{-1}$ by products with $\matM_{\matX}^{-1}$, and $\matB$
is accessed only through matrix-vector products.

\begin{table}[H]
\caption{\label{tab:costs}Summary of the cost of the Riemnnian components
on the generalized Stiefel manifold}

\centering{}{\scriptsize{}}%
\begin{tabular}{|>{\centering}p{8cm}|c|}
\hline 
\textbf{\scriptsize{}Operation} & \textbf{\scriptsize{}Cost}\tabularnewline
\hline 
\hline 
{\scriptsize{}Retraction maps (Eqs.~(\ref{eq:stiefelpolarretraction}),
(\ref{eq:stiefelqrretraction}), (\ref{eq:stiefelcayleyretraction}))} & {\scriptsize{}$O\left(T_{\matB}p+dp^{2}\right)$}\tabularnewline
\hline 
{\scriptsize{}Inverse of the polar-based retraction (Eq.~(\ref{eq:stiefelpolarretractioninv}))} & {\scriptsize{}$O\left(T_{\matB}p+dp^{2}\right)$}\tabularnewline
\hline 
{\scriptsize{}Inverse of the QR-based retraction (Eq.~(\ref{eq:stiefelqrretractioninv}))} & {\scriptsize{}$O\left(T_{\matB}p+dp^{2}+p^{4}\right)$}\tabularnewline
\hline 
{\scriptsize{}Vector Transport, associated with retractions (Eqs.~(\ref{eq:vectransportpolar}),
(\ref{eq:vectransportQR}), \cite[Eq. (16)]{zhu2017riemannian})} & {\scriptsize{}$O\left(T_{\matB}p+dp^{2}\right)$}\tabularnewline
\hline 
{\scriptsize{}Inner product on the tangent space (Eq.~(\ref{eq:stiefelmetric}))} & {\scriptsize{}$O\left(T_{\matM}p+dp\right)$}\tabularnewline
\hline 
{\scriptsize{}Orthogonal projections on the tangent/normal space,
$\matM_{\matX}$ metric (Eqs.~(\ref{eq:stiefelprojtotangent}),(\ref{eq:stiefelprojtonormal}))} & {\scriptsize{}$O\left(T_{\matB}p+T_{\matM^{-1}}p+dp^{2}\right)$}\tabularnewline
\hline 
{\scriptsize{}Orthogonal projections on the tangent/normal space,
$\matB$ metric (Eqs.~(\ref{eq:stiefelprojtotangentBmetric}),(\ref{eq:stiefelprojtonormalBmetric}))} & {\scriptsize{}$O\left(T_{\matB}p+dp^{2}\right)$}\tabularnewline
\hline 
{\scriptsize{}Vector Transport, based on the orthogonal projection
(Eq.~(\ref{eq:vectransstiefeldep}))} & {\scriptsize{}$O\left(T_{\matB}p+T_{\matM^{-1}}p+dp^{2}\right)$}\tabularnewline
\hline 
{\scriptsize{}Riemannian gradient computation (Eq.~(\ref{eq:Rgrad_st_B}))} & {\scriptsize{}$O\left(T_{\matB}p+T_{\matM^{-1}}p+dp^{2}+T_{\nabla\bar{f}}\right)$}\tabularnewline
\hline 
{\scriptsize{}Applying the Riemannian Hessian to a tangent vector
(Eq.~(\ref{eq:StiefelHessian-explicite}))} & {\scriptsize{}$O\left(T_{\matB}p+T_{\matM^{-1}}p+T_{\matM}p+dp^{2}+T_{\nabla\bar{f}}+T_{\nabla^{2}\bar{f}}\right)$}\tabularnewline
\hline 
\end{tabular}{\scriptsize\par}
\end{table}

\subsection{\label{subsec:Product-Manifold-of}Product Manifold of Generalized
Stiefel Manifolds}

In some cases it is desirable to solve optimization problems with
several sets of variables, in which each set of variables is constrained
to a different generalized Stiefel manifold. For example, the CCA
problem is formulated as an optimization problem with two generalized
orthogonality constraints. Such cases are easily addressed by using
the notion of product manifold \cite[ Section 3.1.6]{AMS09}. Here,
we briefly summarize how it applies to our settings.

The basic idea of the product manifold of generalized Stiefel manifolds
is to simply consider the Cartesian product of separately computed
Riemannian components on each of the manifolds in the product. In
particular, when the number of columns is equal for all the generalized
Stiefel manifolds in the product, then it is possible to simply stack
the component matrices on top of each other, and performing the operations
separably on each manifold.

Specifically, Let $\matB_{1},\dots,\matB_{k}$ be SPD matrices, where
the dimension of $\matB_{i}$ is $d_{i}\times d_{i}$, and denote
$d=d_{1}+\dots+d_{k}$. Suppose that the goal is to minimize $f(\matX_{1},\dots,\matX_{k})=f(\matX)$
with the constraint $\matX_{i}\in\stiefel_{\matB_{i}}(p,d_{i})$ for
$i=1,\dots,k$. The problem can be solved using Riemannian optimization
on the product manifold $\stiefel_{\matB_{1}}(p,d_{1})\times\stiefel_{\matB_{2}}(p,d_{2})\times\dots\times\stiefel_{\matB_{k}}(p,d_{k})$,
i.e., $\matX\in\stiefel_{\matB_{1}}(p,d_{1})\times\stiefel_{\matB_{2}}(p,d_{2})\times\dots\times\stiefel_{\matB_{k}}(p,d_{k})$.
Indeed, for the product manifold, there is a natural way to define
the differentiable structure so that manifold topology of $\stiefel_{\matB_{1}}(p,d_{1})\times\stiefel_{\matB_{2}}(p,d_{2})\times\dots\times\stiefel_{\matB_{k}}(p,d_{k})$
is the product topology. However, to employ Riemannian optimization
it is also necessary to define a metric on the product manifold.

Suppose that on each $\stiefel_{\matB_{k}}(p,d_{k})$ the metric is
defined by a smooth mapping $\matX_{i}\mapsto\matM_{\matX_{i}}^{(i)}$
such that $\matM_{\matX_{i}}^{(i)}$ is an SPD matrix (i.e., the metric
$g^{(i)}$ on $\stiefel_{\matB_{i}}(p,d_{i})$ is defined in ambient
coordinates by $g_{\matX}^{(i)}(\eta_{\matX},\xi_{\matX})=\Trace{\eta_{\matX}^{\T}\matM_{\matX_{i}}^{(i)}\xi_{\matX}}$).
The product manifold $\stiefel_{\matB_{1}}(p,d_{1})\times\stiefel_{\matB_{2}}(p,d_{2})\times\dots\times\stiefel_{\matB_{k}}(p,d_{k})$
is a Riemannian submanifold of $\R^{d_{1}\times p}\times\R^{d_{2}\times p}\times\dots\times\R^{d_{k}\times p}$
endowed with the product metric (sum of the metric values on each
product component). Since $\R^{d_{1}\times p}\times\R^{d_{2}\times p}\times\dots\times\R^{d_{k}\times p}$
is naturally isomorphic to $\R^{d\times p}$ by stacking the matrices
on top of each other, then $\stiefel_{\matB_{1}}(p,d_{1})\times\stiefel_{\matB_{2}}(p,d_{2})\times\dots\times\stiefel_{\matB_{k}}(p,d_{k})$
can be viewed as a Riemannian embedded submanifold of $\R^{d\times p}$
endowed with the metric defined by the $d\times d$ matrix $\matM_{\matX}\coloneqq\blockdiag{\matM_{\matX_{1}}^{(1)},\matM_{\matX_{2}}^{(2)},...,\matM_{\matX_{k}}^{(k)}}$,
and the mapping $\matX\mapsto\matM_{\matX}$ is smooth.

The various notions introduced previously now extend to the product
manifold in a straightforward way. Indeed, the tangent space of $\stiefel_{\matB_{1}}(p,d_{1})\times\stiefel_{\matB_{2}}(p,d_{2})\times\dots\times\stiefel_{\matB_{k}}(p,d_{k})$
is the Cartesian product of tangent spaces of each of the generalized
Stiefel manifolds. The retraction and vector transport, and orthogonal
projection on the tangent space is stacking the operations performed
separably on each manifold on top of each other. The Riemannian gradient
is computed using the orthogonal projection to the tangent space after
pre-multiplying by $\matM_{\matX}^{-1}$, i.e., $\grad{f(\matX)}=\Pi_{\matX}\left(\matM_{\matX}^{-1}\nabla\bar{f}(\matX)\right)$
for $\matX\in\stiefel_{\matB_{1}}(p,d_{1})\times\stiefel_{\matB_{2}}(p,d_{2})\times\dots\times\stiefel_{\matB_{k}}(p,d_{k})$,
where $\Pi_{\matX}\left(\cdot\right)$ is stacking the orthogonal
projections on the tangent space of each of the manifolds on top of
each other. The normal space is the product of the normal spaces of
each of the manifolds. Similarly to Subsection \ref{subsec:Metric-Related-Notions},
for the next components we assume $\matM_{\matX}$ is constant. The
Weingarten map is again obtained by stacking the Weingarten maps of
each of the manifolds 
\[
W_{\matX}\left(\xi_{\matX},\matU_{\matX}\right)=\left[\begin{array}{c}
W_{\matX_{1}}\left(\xi_{\matX_{1}},\matU_{\matX_{1}}\right)\\
\vdots\\
W_{\matX_{k}}\left(\xi_{\matX_{k}},\matU_{\matX_{k}}\right)
\end{array}\right]
\]
where $W_{\matX_{i}}\left(\xi_{\matX_{i}},\matU_{\matX_{i}}\right)$
is the Weingarten map on $\stiefel_{\matB_{i}}(p,d_{i})$. The Riemannian
connection on the product manifold is the classical directional derivative
on $\R^{d\times p}$ projected on the tangent space. Thus, the Riemannian
Hessian can be computed using the same formula for the Riemannian
Hessian on the generalized Stiefel manifold,  (\ref{eq:StiefelHessian}),
following similar reasoning as in Appendix \ref{subsec:Metric-Related-Notions-1}.

In the above, we assume the number of columns in each Stiefel component
is the same in all the manifolds in the product. One can also work
on the product manifold $\stiefel_{\matB_{1}}(p_{1},d_{1})\times\stiefel_{\matB_{2}}(p_{2},d_{2})\times\dots\times\stiefel_{\matB_{k}}(p_{k},d_{k})$
where the $p_{1},\dots,p_{k}$ are not necessarily equal. In this
case, we cannot simply stack the tangent vectors etc., but can still
work with Cartesian product of the different components, and operators
like $\matM_{\matX}$ and $\matB$ that operate on each component
separately. Logically, this is the same as we do above for $p_{1}=\dots=p_{k}$,
although the description is somewhat more complex, so we omit the
details.

\subsection{\label{subsec:Metric-selection-and-Hessian} Metric Selection and
Riemannian Hessian Conditioning}

In this subsection we discuss the effects of metric selection with
relation to the condition number of the Riemannian Hessian at the
optimum. Similarly to the unconstrained case, the condition number
of the Riemannian Hessian affects the asymptotic convergence of the
various optimization algorithms -- see \cite[Theorem 4.5.6, Theorem 7.4.11 and Eq. (7.50)]{AMS09}.
We remark that there are also (worst-case) global convergence results
which guarantee sublinear convergence to first and second order (approximate)
critical points (e.g., \cite{boumal2019global}). However, these guarantees
require additional assumptions, e.g., Lipschitz gradient for first-order
conditions and Lipschitz Hessian with second-order retraction for
second-order conditions. Moreover, these guarantees do not depend
on the condition number of the Riemannian Hessian. In practice, as
the iterations progress linear convergence is observed (see experiments
in Subsection \ref{subsec:CCA}) as guaranteed by \cite[Theorem 4.5.6]{AMS09},
and for smaller condition number the convergence is faster.

For simplicity of analysis, consider the case $p=1$, i.e., the generalized
Stiefel manifold in this case is an ellipsoid $\elp$. We also assume
that for all $\x\in\elp$ we have $\matM_{\x}=\matM$ for some fixed
SPD matrix $\matM$. In order to analyze the condition number of the
Riemannian Hessian at the optimum recall that the Riemannian Hessian
is self-adjoint with respect to the Riemannian metric (see \cite[Propositin 5.5.3]{AMS09}).
Thus, its condition number at the optimum, $\x^{\star}$, can be found
using the ratio between the maximal and minimal value of the Rayleigh
quotient 
\[
q(\xi_{\x^{\star}})\coloneqq\frac{g_{\x^{\star}}(\xi_{\x^{\star}},\hess{f(\x^{\star})}[\xi_{\x^{\star}}])}{g_{\x^{\star}}(\xi_{\x^{\star}},\xi_{\x^{\star}})}.
\]
Using  (\ref{eq:StiefelHessian-explicite}), the Riemannian Hessian
for $p=1$ is reduced to 
\begin{eqnarray*}
\hess{f(\x^{\star})[\eta_{\x^{\star}}]} & = & \Pi_{\x^{\star}}\left(\matM_{\x^{\star}}^{-1}\left[\nabla^{2}\bar{f}(\x^{\star})-\left(\left(\x^{\star}\right)^{\T}\nabla\bar{f}(\x^{\star})-g_{\x^{\star}}(\x^{\star},\grad{f(\x^{\star})})\right)\matB\right]\eta_{\x^{\star}}\right)\ .
\end{eqnarray*}

Recall that $\grad{f(\x^{\star})}=0$, also the projection on the
tangent space is self-adjoint with respect to the Riemannian metric,.
(\ref{eq:projproperty2}), and for any $\xi_{\x^{\star}}\in T_{\x^{\star}}\elp$
we have $\Pi_{\x^{\star}}\left(\xi_{\x^{\star}}\right)=\xi_{\x^{\star}}$,
we get:
\begin{eqnarray*}
q(\xi_{\x^{\star}}) & = & \frac{\xi_{\x^{\star}}^{\T}\matM_{\x^{\star}}\Pi_{\x^{\star}}\left(\matM_{\x^{\star}}^{-1}\left[\nabla^{2}\bar{f}(\x^{\star})-\left(\left(\x^{\star}\right)^{\T}\nabla\bar{f}(\x^{\star})\right)\matB\right]\xi_{\x^{\star}}\right)}{\xi_{\x^{\star}}^{\T}\matM_{\x^{\star}}\xi_{\x^{\star}}}\\
 & = & \frac{\left(\Pi_{\x^{\star}}\left(\xi_{\x^{\star}}\right)\right)^{\T}\left[\nabla^{2}\bar{f}(\x^{\star})-\left(\left(\x^{\star}\right)^{\T}\nabla\bar{f}(\x^{\star})\right)\matB\right]\xi_{\x^{\star}}}{\xi_{\x^{\star}}^{\T}\matM_{\x^{\star}}\xi_{\x^{\star}}}\\
 & = & \frac{\xi_{\x^{\star}}^{\T}\left[\nabla^{2}\bar{f}(\x^{\star})-\left(\left(\x^{\star}\right)^{\T}\nabla\bar{f}(\x^{\star})\right)\matB\right]\xi_{\x^{\star}}}{\xi_{\x^{\star}}^{\T}\matM_{\x^{\star}}\xi_{\x^{\star}}}\ .
\end{eqnarray*}
This is the Rayleigh quotient of the matrix pencil 
\[
\left(\nabla^{2}\bar{f}(\x^{\star})-\left(\left(\x^{\star}\right)^{\T}\nabla\bar{f}(\x^{\star})\right)\matB,\matM_{\x^{\star}}\right)
\]
on $T_{\x^{\star}}\elp$. So, if we want to bound the condition number
of the Riemannian Hessian at the optimum we need to look at the pencil
\begin{equation}
\left(\Pi_{\x^{\star}}\left(\nabla^{2}\bar{f}(\x^{\star})-\left(\left(\x^{\star}\right)^{\T}\nabla\bar{f}(\x^{\star})\right)\matB\right)\Pi_{\x^{\star}},\Pi_{\x^{\star}}\matM_{\x^{\star}}\Pi_{\x^{\star}}\right).\label{eq:pencil}
\end{equation}
Therefore, choosing a preconditioning scheme $\x\mapsto\matM_{\x}$
such that $\matM_{\x}$ is SPD for any $\x\in\elp$ and 
\begin{equation}
\matM_{\x^{\star}}\approx\nabla^{2}\bar{f}(\x^{\star})-\left(\left(\x^{\star}\right)^{\T}\nabla\bar{f}(\x^{\star})\right)\matB\label{eq:Lagrangian_metric}
\end{equation}
will precondition the Riemannian Hessian at the optimum. One such
example can be found in \cite{mor2020solving}. In addition, the preconditioners
proposed in \cite{MS16}, which are inspired by the Lagrangian, can
be viewed in such manner, thus, approximating the Riemannian Newton
method. For the generalized Stiefel manifold with $p>1$ such a choice
is less obvious, and we leave it for future work.

Recall that the standard choice for metric selection on the generalized
Stiefel manifold with $p=1$ is $\matM_{\x}=\matB$ for all $\x\in\elp$.
If $\nabla^{2}\bar{f}(\x^{\star})$ is well conditioned, it is often
the case that the pencil ~(\ref{eq:pencil}) is well conditioned
under certain assumptions. We demonstrate this in Section~\ref{sec:importance}
for the problem of finding the leading correlation in CCA. In such
cases, if we use a preconditioning scheme $\x\mapsto\matM_{\x}$ such
that $\matM_{\x^{\star}}\approx\matB$, the condition number grows
by at most $\kappa(\matB,\matM_{\x^{\star}})$, so if that quantity
is small (i.e., $\matM_{\x^{\star}}$ well approximates $\matB$)
we can expect fast convergence.

\section{\label{sec:importance} Theoretical and Numerical Illustrations}

\subsection{Simple Theoretical Examples}

Our proposed preconditioning strategy for orthogonality constrained
problems is based on using a preconditioning scheme to define the
Riemannian metric. In this section we illustrate this point using
a couple of simple examples. All examples correspond to the case $p=1$,
i.e., the ellipsoid.
\begin{example}
\textbf{Linear Objective.} Consider the following problem 
\[
\max_{\x\in\R^{d}}\b^{\T}\x\st\x^{\T}\matB\x=1
\]
for some vector $0\neq\b\in\R^{d}$, where $\matB\in\R^{d\times d}$.
It is easy to show that the solution is $\x^{\star}=\matB^{-1}\b/\BNorm{\matB^{-1}\b}.$
It is well known that solving a linear system is equivalent to an
unconstrained minimization of a quadratic objective. Here we can see
that solving a linear system is also equivalent to maximizing a linear
objective subject to a quadratic constraint. Note that this problem
is constrained on the ellipsoid manifold $\elp$. Let the inner product
on each tangent space (the Riemannian metric) be endowed from the
ambient space $\R^{d}$. Using $\elp$ with a metric selection $g_{\x}(,\eta_{\x})=\matM_{\x}\eta_{\x}$
(in ambient coordinates), where $\x\mapsto\matM_{\x}\in\R^{d\times d}$
is a smooth mapping that maps $\x\in\elp$ to an SPD matrix $\matM_{\x}$,
the Riemannian gradient is 
\[
\grad{f(\x)}=(\matI_{n}-(\x^{\T}\matB\mat M_{\x}^{-1}\matB\x)^{-1}\mat M_{\x}^{-1}\matB\x\x^{\T}\matB)\mat M_{\x}^{-1}\b
\]
since the Euclidean gradient is simply $\b$, independent of $\x$.
Thus, using Riemannian gradient ascent on $\elp$ with the polar based
retraction, ~(\ref{eq:stiefelpolarretraction}), we get the iteration
\begin{eqnarray*}
\y_{k+1} & = & \x_{k}+\alpha_{k}\left(\mat M_{\x_{k}}^{-1}\b-\frac{\x_{k}^{\T}\matB\mat M_{\x_{k}}^{-1}\b}{\x_{k}^{\T}\matB\mat M_{\x_{k}}^{-1}\matB\x_{k}}\mat M_{\x_{k}}^{-1}\matB\x_{k}\right)\\
\x_{k+1} & = & \frac{\y_{k+1}}{\BNorm{\y_{k+1}}}\,.
\end{eqnarray*}
We see, as expected, that the iterations depend on the choice of the
Riemannian metric defined by the matrix $\mat M_{\x}$. If we impose
the metric $\matM_{\x}=\matB$ for all $\x\in\elp$, and take step
size $\alpha_{0}=1/\x_{k}^{\T}\b$, then the iterations reduce to
$\x_{1}=\matB^{-1}\b/\BNorm{\matB^{-1}\b}$, and the problem is solved
in a single iteration.

As expected, with $\matM_{\x}=\matB$ for all $\x$, the Riemannian
Hessian at $\x^{\star}$ is well conditioned. Indeed, we have 
\[
\hess{f(\x^{\star})}=-\Pi_{\x^{\star}}\left(\left(\x^{\star}{}^{\T}\b\right)\matI_{d}\right),
\]
and its corresponding Rayleigh quotient is
\[
q(\xi_{\x^{\star}})=\frac{\xi_{\x^{\star}}^{\T}\matB\left[-\Pi_{\x^{\star}}\left(\left(\x^{\star}{}^{\T}\b\right)\matI_{d}\right)\right]\xi_{\x^{\star}}}{\xi_{\x^{\star}}^{\T}\matB\xi_{\x^{\star}}}=-\left(\x^{\star}{}^{\T}\b\right)=-\BNorm{\matB^{-1}\b}\ ,
\]
which is constant so the condition number equals $1$. Note that the
metric selection $\matM_{\x}=\matB$ also satisfies  (\ref{eq:Lagrangian_metric}).
\end{example}

\begin{example}
\textbf{Inverse Power Iteration.} Consider the following problem 
\[
\max_{\x\in\R^{d}}\frac{1}{2}\x^{\T}\x\st\x^{\T}\matB\x=1
\]
where $\matB\in\R^{d\times d}$ is an SPD matrix. The solution is
an eigenvector corresponding the smallest eigenvalue of $\matB$,
$\lambda_{d}(\matB)$, (which is also the eigenvector corresponding
to the maximal eigenvalue of $\matB^{-1}$), since this problem is
equivalent to maximizing the Rayleigh quotient $\ensuremath{\x^{\T}\x/\x^{\T}\matB\x}$.
Note that this problem is constrained on the ellipsoid manifold $\elp$.
Using $\elp$ with metric selection $g_{\x}(,\eta_{\x})=\mat M_{\x}\eta_{\x}$
(in ambient coordinates), where $\mat M_{\x}\in\R^{d\times d}$ is
an SPD matrix for any $\x\in\elp$, the Riemannian gradient is 
\[
\grad{f(\x)}=(\matI_{d}-(\x^{\T}\matB\mat M_{\x}^{-1}\matB\x)^{-1}\mat M_{\x}^{-1}\matB\x\x^{\T}\matB)\mat M_{\x}^{-1}\x
\]
since the Euclidean gradient is $\x$. Thus, using Riemannian gradient
ascent on $\elp$ with the polar based retraction, ~(\ref{eq:stiefelpolarretraction}),
we get the iteration
\begin{eqnarray*}
\y_{k+1} & = & \x_{k}+\alpha_{k}\left(\mat M_{\x}^{-1}\x_{k}-\frac{\x_{k}^{\T}\matB\mat M_{\x}^{-1}\x_{k}}{\x_{k}^{\T}\matB\mat M_{\x}^{-1}\matB\x_{k}}\mat M_{\x}^{-1}\matB\x_{k}\right)\\
\x_{k+1} & = & \frac{\y_{k+1}}{\BNorm{\y_{k+1}}}\,.
\end{eqnarray*}
If we impose the metric $\mat M_{\x}=\matB$ for all $\x\in\elp$,
and take step sizes $\alpha_{k}=(\x_{k}^{\T}\x_{k})^{-1},$ then the
iterations reduce to $\x_{k+1}=\matB^{-1}\x_{k}/\TNorm{\matB^{-1}\x_{k}}$,
i.e., the inverse power method, which is well known for its good convergence
properties for eigenvalues near zero.

Let us examine the Riemannian Hessian at the optimal point $\x^{\star}$
(i.e $\left(\x^{\star}\right)^{\T}\x^{\star}=1/\lambda_{\min}(\matB)=1/\lambda_{d}(\matB)$):
\[
\hess{f(\x^{\star})}=\Pi_{\x^{\star}}\left(\matB^{-1}\left[\matI_{d}-\left(\left(\x^{\star}\right)^{\T}\x^{\star}\right)\matB\right]\right)=\Pi_{\x^{\star}}\left(\matB^{-1}\left[\matI_{d}-\left(1/\lambda_{d}(\matB)\right)\matB\right]\right).
\]
The corresponding Rayleigh quotient is reduced to the following form
using similar reasoning as in Subsection \ref{subsec:Metric-selection-and-Hessian}:
\[
q(\xi_{\x^{\star}})=\frac{\xi_{\x^{\star}}^{\T}\matB\left[\Pi_{\x^{\star}}\left(\matB^{-1}\left[\matI_{d}-\left(1/\lambda_{d}(\matB))\right)\matB\right]\right)\right]\xi_{\x^{\star}}}{\xi_{\x^{\star}}^{\T}\matB\xi_{\x^{\star}}}=\frac{\xi_{\x^{\star}}^{\T}\left[\matI_{d}-\left(1/\lambda_{d}(\matB)\right)\matB\right]\xi_{\x^{\star}}}{\xi_{\x^{\star}}^{\T}\matB\xi_{\x^{\star}}}.
\]
Thus, the eigenvalues of the Riemannian Hessian at $\x^{\star}$ correspond
to the generalized eigenvalues of the matrix pencil $\left(\matI_{d}-\left(1/\lambda_{d}(\matB)\right)\matB,\matB\right)$
on $T_{\x^{\star}}\elp$, i.e., the eigenvalues of $\matB^{-1}$ deflated
by $-1/\lambda_{d}(\matB)$ on $T_{\x^{\star}}\elp$. Moreover, since
$\xi_{\x^{\star}}\in T_{\x^{\star}}\elp$, we have $\xi_{\x^{\star}}^{\T}\matB\x^{\star}=0$,
thus $\xi_{\x^{\star}}$ is constrained not to correspond to $1/\lambda_{d}(\matB)$.
Assume that $\lambda_{d-1}(\matB)>\lambda_{d}(\matB)$, then the condition
number is bounded by 
\[
\frac{1/\lambda_{d}(\matB)-1/\lambda_{1}(\matB)}{1/\lambda_{d}(\matB)-1/\lambda_{d-1}(\matB)},
\]
which for $\lambda_{d}(\matB)$ that is close to $0$, and $\lambda_{d-1}(\matB)\gg0$
is close to $1$.

Note that if we try to impose the metric $\mat M_{\x}=-(\matI_{d}-\left(1/\lambda_{d}(\matB)\right)\matB)$
for all $\x\in\elp$ (following  (\ref{eq:Lagrangian_metric})), we
have that $\mat M_{\x}$ is singular since it has a zero eigenvalue
(corresponding to the eigenvector $\x^{\star}$), thus it cannot be
a Riemannian metric inherited from the ambient space $\R^{d}$.
\end{example}

\subsection{\label{subsec:CCA}Canonical Correlation Analysis: Theory and Experiment}

In this subsection we illustrate our approach on the problem of finding
the top correlation between two datasets. This problem can be written
as optimization problem whose constraint set is the product of two
ellipsoids.

CCA, originally introduced by \cite{hotelling1936relations}, is a
well-established method in statistical learning with numerous applications
(e.g., \cite{sun2010scalable,chaudhuri2009multi,dhillon2011multi,dhillon2012two,su2012discriminant,kim2007discriminative}).
In CCA the relation between a pair of datasets in matrix form is analyzed,
where the goal is to find the directions of maximal correlation between
a pair of observed variables. In the language of linear algebra, CCA
measures the similarities between two subspaces spanned by the columns
of of the two matrices. Here, we consider a regularized version of
CCA defined below:
\begin{defn}
Let $\matX\in\R^{n\times d_{\x}}$ and $\matY\in\R^{n\times d_{\y}}$
be two data matrices, and $\lambda_{\x},\lambda_{\y}\geq0$ be two
regularization parameter. Let 
\[
q=\max\left(\rank{\matX^{\T}\matX+\lambda_{\x}\matI_{d_{\x}}},\rank{\matY^{\T}\matY+\lambda_{\y}\matI_{d_{\y}}}\right).
\]
The $(\lambda_{\x},\lambda_{\y})$ canonical correlations $\sigma_{1}\geq\dots\geq\sigma_{q}$
and the $(\lambda_{\x},\lambda_{\y})$ canonical weights $\u_{1},\dots,\u_{q}\in\R^{d_{x}}$,
$\v_{1},\dots,\v_{q}\in\R^{d_{y}}$, are the ones that maximize 
\[
\Trace{\matU^{\T}\matX^{\T}\matY\matV}
\]
subject to 
\[
\matU^{\T}(\matX^{\T}\matX+\lambda_{\x}\matI_{d_{\x}})\matU=\matI_{d_{\x}},\quad\matV^{\T}(\matY^{\T}\matY+\lambda_{\y}\matI_{d_{\y}})\matV=\matI_{d_{\y}}
\]
where $\matU^{\T}\matX^{\T}\matY\matV=\diag{\sigma_{1},\dots,\sigma_{q}}$,
$\matU=\left[\begin{array}{ccc}
\u_{1} & \dots & \u_{q}\end{array}\right]\in\R^{d_{x}\times q}$ and $\matV=\left[\begin{array}{ccc}
\v_{1} & \dots & \v_{q}\end{array}\right]\in\R^{d_{y}\times q}$.
\end{defn}

In this paper, we focus on finding the top correlation, i.e., finding
$\sigma_{1},\u_{1}$ and $\v_{1}$. It is useful to introduce the
following notations: 
\[
\Sigma_{\x\x}=\matX^{\T}\matX+\lambda_{\x}\matI_{d_{\x}},\Sigma_{\y\y}=\matY^{\T}\matY+\lambda_{\y}\matI_{d_{\y}},\Sigma_{\x\y}=\matX^{\T}\matY\,.
\]
Restricting to finding the top correlation, the optimization problem
becomes:
\begin{equation}
\max\u^{\T}\Sigma_{\x\y}\v\st\u\in\elpsigx,\v\in\elpsigy\label{eq:CCAdemo}
\end{equation}
It is well known (\cite{bjorck1973numerical}) that the optimal solution
of Problem (\ref{eq:CCAdemo}) is (up to the sign of the vectors)
\begin{equation}
\u_{1}\coloneqq\Sigma_{\x\x}^{-\nicehalf}\phi\quad\v_{1}\coloneqq\Sigma_{\y\y}^{-\nicehalf}\psi\label{eq:CCAoptsol}
\end{equation}
where $\phi\in\R^{d_{x}}$ and $\psi\in\R^{d_{y}}$ are the left and
right unit-length singular vector corresponding to the largest singular
value $\sigma_{1}$ of the matrix 
\begin{equation}
\matR\coloneqq\Sigma_{\x\x}^{-\nicehalf}\Sigma_{\x\y}\Sigma_{\y\y}^{-\nicehalf}\:.\label{eq:Tmatrix}
\end{equation}

In order to conveniently use the Riemannian optimization framework,
we also denote $d=d_{\x}+d_{\y}$, and $\z=[\u^{\T},\v^{\T}]^{\T}\in\R^{d}$
where $\u\in\R^{d_{\x}}$ and $\v\in\R^{d_{\y}}$. Then the constraint
set is a product manifold of two ellipsoids $\z\in\elpCCA\coloneqq\elpsigx\times\elpsigy$.
The objective function to be minimized is then
\begin{equation}
f(\z)=-\frac{1}{2}\z^{\T}\left[\begin{array}{cc}
0 & \Sigma_{\x\y}\\
\Sigma_{\x\y}^{\T} & 0
\end{array}\right]\z\ .\label{eq:costriemanniancca}
\end{equation}
We endow the manifold $\elpsigx$ and $\elpsigy$ with a metric defined
by two preconditioning schemes $\u\mapsto\matM_{\u}^{(\x\x)}$ and
$\v\mapsto\matM_{\v}^{(\y\y)}$. The metric on the product manifold
$\elpCCA$ is defined by $\z\mapsto\matM_{\z}=\blockdiag{\matM_{\u}^{(\x\x)},\matM_{\v}^{(\y\y)}}$
as explained in Section~\ref{subsec:Product-Manifold-of}. Using
the formulas in Section~\ref{subsec:Metric-Related-Notions} we find
that the Riemannian gradient and the Riemannian Hessian (at the critical
points or if $\matM_{\z}\coloneqq\matM=\blockdiag{\matM^{(\x\x)},\matM^{(\y\y)}}$)
are given by: 
\[
\grad f(\z)=\Pi_{\z}\left(\matM_{\z}^{-1}\nabla\bar{f}(\z)\right)=-\left[\begin{array}{c}
\Pi_{\u}\left(\left(\matM_{\u}^{(\x\x)}\right)^{-1}\Sigma_{\x\y}\v\right)\\
\Pi_{\v}\left(\left(\matM_{\v}^{(\y\y)}\right)^{-1}\Sigma_{\x\y}^{\T}\u\right)
\end{array}\right]\ ,
\]
\[
\hess{f(\z)}[\eta_{\z}]=\Pi_{\z}\left(\matM_{\z}^{-1}\left[\begin{array}{cc}
(\u^{\T}\matM^{(\x\x)}\Pi_{\u}^{\perp}\left(\left(\matM^{(\x\x)}\right)^{-1}\Sigma_{\x\y}\v\right))\cdot\Sigma_{\x\x} & -\Sigma_{\x\y}\\
-\Sigma_{\x\y}^{\T} & \left(\v^{\T}\matM^{(\y\y)}\Pi_{\v}^{\perp}\left(\left(\matM^{(\y\y)}\right)^{-1}\Sigma_{\x\y}^{\T}\u\right)\right)\cdot\Sigma_{\y\y}
\end{array}\right]\eta_{\z}\right)\ .
\]
Along with formulas for the retraction and vector transport (see
Subsection \ref{subsec:metric-independent}), various Riemannian optimization
algorithms can be applied to solve Problem (\ref{eq:CCAdemo}).

As expected, at the optimal solution $\z^{\star}=[\u_{1}^{\T},\v_{1}^{\T}]^{\T}$
(see  (\ref{eq:CCAoptsol})) the Riemannian gradient vanishes: $\grad f(\z^{\star})=0$.
Moreover, the Riemannian Hessian at the optimum becomes
\begin{equation}
\hess{f(\z^{\star})}=\Pi_{\z^{\star}}\left(\matM_{\z^{\star}}^{-1}\left[\begin{array}{cc}
\sigma_{1}\cdot\Sigma_{\x\x} & -\Sigma_{\x\y}\\
-\Sigma_{\x\y}^{\T} & \sigma_{1}\cdot\Sigma_{\y\y}
\end{array}\right]\right)\ .\label{eq:CCAHess_opt}
\end{equation}

\begin{figure}[t]
\begin{centering}
\begin{tabular}{ccc}
\includegraphics[width=0.38\columnwidth]{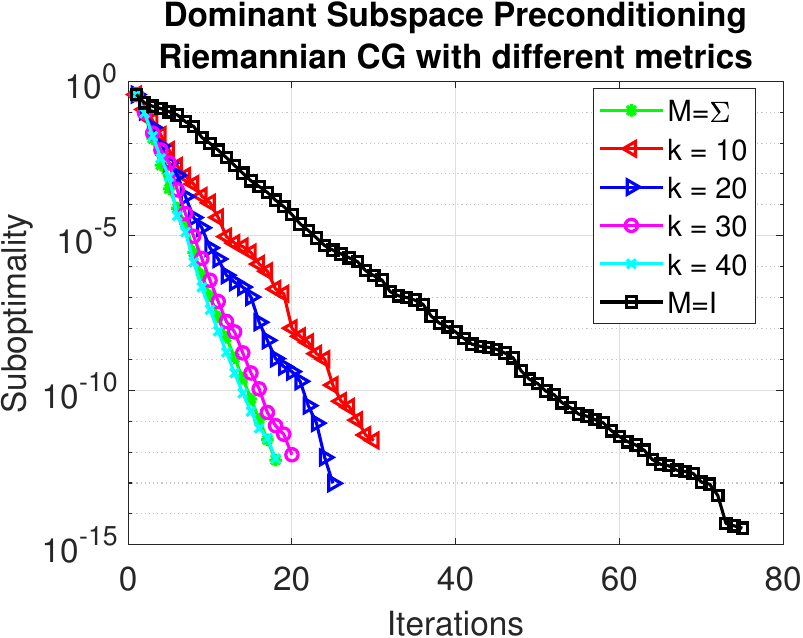} & ~ & \includegraphics[width=0.38\columnwidth]{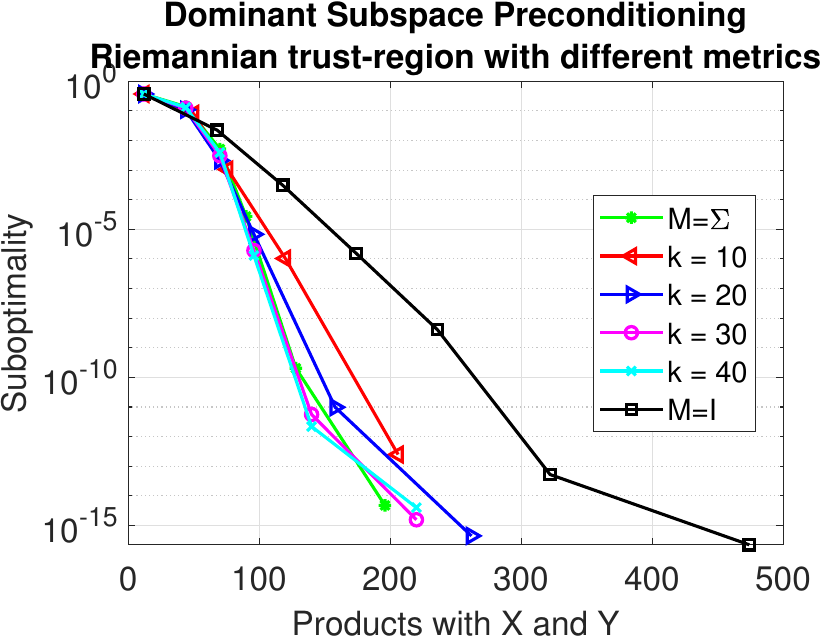}\tabularnewline
\end{tabular}
\par\end{centering}
\caption{\label{fig:Results-for-CCA}Results for CCA with Riemannian conjugate-gradient
(left - suboptimality vs. \#iterations) and Riemannian trust-region
(right - suboptimality vs. products with the data matrices) with various
choices of metrics for $p=1$. The number of leading eigenvalues used
to form the Dominant Subspace Preconditioner is denoted by $k$.}
\end{figure}

Next, we demonstrate the effect of preconditioning on the condition
number of the Riemannian Hessian at $\z^{\star}$. We show that if
the leading correlation is strictly larger than the second largest
one, and we select a smooth preconditioning scheme $\z\mapsto\matM_{\z}$
such that $\matM_{\z^{\star}}=\Sigma\coloneqq\blockdiag{\Sigma_{\x\x},\Sigma_{\y\y}}$,
the condition number of the Riemannian Hessian at the optimum is equal
to $(\sigma_{1}+\sigma_{2})/(\sigma_{1}-\sigma_{2})$. Thus, if the
leading correlation gap $\sigma_{1}-\sigma_{2}$ is $O(\sigma_{1})$
then the condition number at the optimum is $O(1)$, and we can expect
fast convergence (dependence on the gap between the correlations is
expected). Furthermore, if we select a smooth preconditioning scheme
$\z\mapsto\matM_{\z}$ such that $\matM_{\z^{\star}}\approx\Sigma$
(see for example Fig. \ref{fig:Results-for-CCA}) the condition number
bound grows by at most a small factor: $\kappa\left(\matB,\matM_{\z^{\star}}\right)$.
\begin{lem}
\label{lem:CCA_cond_num}Assuming $\sigma_{1}-\sigma_{2}>0$ and that
$\Sigma$ is an SPD matrix, if $\elpCCA$ is equipped with a metric
defined by a smooth preconditioning scheme $\z\mapsto\matM_{\z}$
such that $\matM_{\z^{\star}}=\Sigma$, then the condition number
of Riemannian Hessian on $\elpCCA$ of  (\ref{eq:costriemanniancca})
at $\z^{\star}$ is equal to $\frac{\sigma_{1}+\sigma_{2}}{\sigma_{1}-\sigma_{2}}$.
Additionally, if $\matM_{\z^{\star}}\approx\Sigma$ then the condition
number is bounded by $\frac{\sigma_{1}+\sigma_{2}}{\sigma_{1}-\sigma_{2}}\cdot\kappa\left(\matB,\matM_{\z^{\star}}\right)$.
\end{lem}

\begin{proof}
In order to bound the condition number of Riemannian Hessian on $\elpCCA$
of  (\ref{eq:costriemanniancca}) at $\z^{\star}$ we use the Courant-Fischer
Theorem for the compact self-adjoint linear operator $\hess{f(\z^{\star})[\cdot]}:T_{\z^{\star}}\elpCCA\to T_{\z^{\star}}\elpCCA$
over the finite dimensional vector space $T_{\matZ}\elpCCA$:
\begin{eqnarray*}
\lambda_{k}(\hess{f(\z^{\star})}) & = & \min_{U,\dim(U)=k-1}\max_{\mat 0\neq\xi_{\z^{\star}}\in U^{\perp}}q(\xi_{\z^{\star}}),
\end{eqnarray*}
\begin{eqnarray*}
\lambda_{k}(\hess{f(\z^{\star})}) & = & \max_{U,\dim(U)=k}\min_{\mat 0\neq\xi_{\z^{\star}}\in U}q(\xi_{\z^{\star}}),
\end{eqnarray*}
where 
\[
q(\xi_{\z^{\star}})\coloneqq\frac{g_{\z^{\star}}(\xi_{\z^{\star}},\hess{f(\z^{\star})[\xi_{\z^{\star}}]})}{g_{\z^{\star}}(\xi_{\z^{\star}},\xi_{\z^{\star}})},
\]
is the Rayleigh quotient. In the above, $\lambda_{k}(\hess{f(\z^{\star})})$
is the $k$-th largest eigenvalue (i.e., eigenvalues are ordered in
a descending order) of $\hess{f(\z^{\star})}$, and $U$ is a linear
subspace of $T_{\z^{\star}}\elpCCA$. In particular, the maximal and
minimal eigenvalues are given by the formulas
\[
\lambda_{\max}(\hess{f(\z^{\star})})=\max_{\mat 0\neq\xi_{\z^{\star}}\in T_{\z^{\star}}\elpCCA}q(\xi_{\z^{\star}}),
\]
 
\[
\lambda_{\min}(\hess{f(\z^{\star})})=\min_{\mat 0\neq\xi_{\z^{\star}}\in T_{\z^{\star}}\elpCCA}q(\xi_{\z^{\star}}),
\]
and the condition number of the Riemannian Hessian at $\z^{\star}$
is the ratio of these two eigenvalues.
\[
\kappa(\hess{f(\z^{\star})})=\frac{\lambda_{\max}(\hess{f(\z^{\star})})}{\lambda_{\min}(\hess{f(\z^{\star})})}\ .
\]

We begin by simplifying the quotient $q(\xi_{\z^{\star}})$. At the
optimum, $\z^{\star}$, we have $f(\z^{\star})=-\u_{1}^{\T}\Sigma_{\x\y}\v_{1}=-\v_{1}^{\T}\Sigma_{\x\y}^{\T}\u_{1}=-\sigma_{1}$.
The formula for the Riemannian Hessian, $\hess{f(\z^{\star})}$, is
given by  (\ref{eq:CCAHess_opt}). Using the following notation for
the Euclidean Hessian of 
\begin{equation}
\nabla^{2}\bar{f}(\z^{\star})\coloneqq\left[\begin{array}{cc}
0 & -\Sigma_{\x\y}\\
-\Sigma_{\x\y}^{\T} & 0
\end{array}\right],\label{eq:CCAeucHess}
\end{equation}
and $\Sigma$ we can compactly write  (\ref{eq:CCAHess_opt}):
\[
\hess{f(\z^{\star})}=\Pi_{\z^{\star}}\left(\matM_{\z^{\star}}^{-1}\left(\nabla^{2}\bar{f}(\z^{\star})+\sigma_{1}\Sigma\right)\right).
\]
Next, as in Subsection \ref{subsec:Metric-selection-and-Hessian},
recall that $\Pi_{\z^{\star}}$ is self-adjoint with respect to the
Riemannian metric,  (\ref{eq:projproperty2}), and that for any $\xi_{\z^{\star}}\in T_{\z^{\star}}\elpCCA$
we have $\Pi_{\z^{\star}}\left(\xi_{\z^{\star}}\right)=\xi_{\z^{\star}}$,
we get:
\[
q(\xi_{\z^{\star}})=\frac{\xi_{\z^{\star}}^{\T}\left(\nabla^{2}\bar{f}(\z^{\star})+\sigma_{1}\cdot\Sigma\right)\xi_{\z^{\star}}}{\xi_{\z^{\star}}^{\T}\matM_{\z^{\star}}\xi_{\z^{\star}}}=\frac{\xi_{\z^{\star}}^{\T}\left(\nabla^{2}\bar{f}(\z^{\star})+\sigma_{1}\cdot\Sigma\right)\xi_{\z^{\star}}}{\xi_{\z^{\star}}^{\T}\Sigma\xi_{\z^{\star}}}\cdot\frac{\xi_{\z^{\star}}^{\T}\Sigma\xi_{\z^{\star}}}{\xi_{\z^{\star}}^{\T}\matM_{\z^{\star}}\xi_{\z^{\star}}}\ ,
\]
where we use the fact that $\Sigma$ is not singular. Note that the
quotient 
\[
\frac{\xi_{\z^{\star}}^{\T}\left(\nabla^{2}\bar{f}(\z^{\star})+\sigma_{1}\cdot\Sigma\right)\xi_{\z^{\star}}}{\xi_{\z^{\star}}^{\T}\Sigma\xi_{\z^{\star}}},
\]
corresponds to the Rayleigh quotient of the Riemannian Hessian at
$\z^{\star}$ if $\matM_{\z^{\star}}=\Sigma$.

Let us first find the eigenvalues of the Riemannian Hessian for the
case $\matM_{\z^{\star}}=\Sigma$. We perform the following invertible
change of variables $\tilde{\xi}_{\z^{\star}}\coloneqq\Sigma^{\nicehalf}\xi_{\z^{\star}}$,
to find that 
\[
q(\xi_{\z^{\star}})=\frac{\tilde{\xi}_{\z^{\star}}^{\T}\left(\Sigma^{-\nicehalf}\nabla^{2}\bar{f}(\z^{\star})\Sigma^{-\nicehalf}+\sigma_{1}\cdot\matI_{d}\right)\tilde{\xi}_{\z^{\star}}}{\tilde{\xi}_{\z^{\star}}^{\T}\tilde{\xi}_{\z^{\star}}}\coloneqq\tilde{q}(\tilde{\xi}_{\z^{\star}})\ .
\]
Denote the space of vectors $\tilde{\xi}_{\z^{\star}}$ such that
$\Sigma^{-\nicehalf}\tilde{\xi}_{\z^{\star}}\in T_{\z^{\star}}\elpCCA$
by $\Sigma^{\nicehalf}T_{\z^{\star}}\elpCCA$, and the orthogonal
space to it by $(\Sigma^{\nicehalf}T_{\z^{\star}}\elpCCA)^{\perp}$.
The above expression, $\tilde{q}(\tilde{\xi}_{\z^{\star}})$, is the
Rayleigh quotient for the symmetric matrix $\Sigma^{-\nicehalf}\nabla^{2}\bar{f}(\z^{\star})\Sigma^{-\nicehalf}+\sigma_{1}\cdot\matI_{d}$.
Thus, applying the Courant-Fischer theorem for $\tilde{q}(\tilde{\xi}_{\z^{\star}})$,
where $\tilde{\xi}_{\z^{\star}}\in\Sigma^{\nicehalf}T_{\z^{\star}}\elpCCA$,
the minimal and the maximal values of $R(\xi_{\z^{\star}})$, where
$\xi_{\z^{\star}}\in T_{\z^{\star}}\elpCCA$, are the minimal and
the maximal eigenvalues of the matrix $\Sigma^{-\nicehalf}\nabla^{2}\bar{f}(\z^{\star})\Sigma^{-\nicehalf}+\sigma_{1}\cdot\matI_{d}$
in the space $\Sigma^{\nicehalf}T_{\z^{\star}}\elpCCA$.

To find the eigenvalues of the matrix $\Sigma^{-\nicehalf}\nabla^{2}\bar{f}(\z^{\star})\Sigma^{-\nicehalf}+\sigma_{1}\cdot\matI_{d}$
in the space $\Sigma^{\nicehalf}T_{\z^{\star}}\elpCCA$, we first
note that all the eigenvalues of $\Sigma^{-\nicehalf}\nabla^{2}\bar{f}(\z^{\star})\Sigma^{-\nicehalf}$
are $-\sigma_{1}<-\sigma_{2}\leq...\leq-\sigma_{q}\leq0\leq...\leq0\leq\sigma_{q}\leq...\leq\sigma_{2}<\sigma_{1}$
(see \cite{golub1995canonical}). So, all the eigenvalue of $\Sigma^{-\nicehalf}\nabla^{2}\bar{f}(\z^{\star})\Sigma^{-\nicehalf}+\sigma_{1}\cdot\matI_{d}$
are $0<\sigma_{1}-\sigma_{2}\leq\cdots\leq\sigma_{1}-\sigma_{q}\leq\sigma_{1}\leq...\leq\sigma_{q}+\sigma_{1}\leq...\leq\sigma_{2}+\sigma_{1}<2\sigma_{1}$.
Next, note that the eigenspaces of $\Sigma^{-\nicehalf}\nabla^{2}\bar{f}(\z^{\star})\Sigma^{-\nicehalf}+\sigma_{1}\cdot\matI_{d}$
corresponding to the eigenvalues $0$ and $2\sigma_{1}$ is exactly
the two dimensional space $(\Sigma^{\frac{1}{2}}T_{\z^{\star}}\elpCCA)^{\perp}$.
Indeed, according to  (\ref{eq:normalstiefeMmetric}) and Subsection
\ref{subsec:Product-Manifold-of}:
\[
(\Sigma^{\nicehalf}T_{\z^{\star}}\elpCCA)^{\perp}=\span{\left\{ \Sigma^{\nicehalf}\left[\begin{array}{c}
\u_{1}\\
\v_{1}
\end{array}\right],\Sigma^{\nicehalf}\left[\begin{array}{c}
\u_{1}\\
-\v_{1}
\end{array}\right]\right\} \ ,}
\]
where using  (\ref{eq:CCAoptsol})
\[
\Sigma^{\nicehalf}\left[\begin{array}{c}
\u_{1}\\
\v_{1}
\end{array}\right]=\left[\begin{array}{c}
\phi\\
\psi
\end{array}\right]\quad\textrm{and}\quad\Sigma^{\nicehalf}\left[\begin{array}{c}
\u_{1}\\
-\v_{1}
\end{array}\right]=\left[\begin{array}{c}
\phi\\
-\psi
\end{array}\right]\ .
\]
Recall that the normal space $(T_{\z^{\star}}\elpCCA)^{\perp}$ is
the Cartesian product of the normal spaces $(T_{\u_{1}}\elpsigx)^{\perp}$
and $(T_{\v_{1}}\elpsigy)^{\perp}$ which are spanned by $\u_{1}$
and $\v_{1}$ correspondingly when $\matM_{\z^{\star}}=\Sigma$. Thus,
the Cartesian product $(T_{\z^{\star}}\elpCCA)^{\perp}$ can be spanned
by $[\u_{1}^{\T},\v_{1}^{\T}]^{\T}$ and $[\u_{1}^{\T},-\v_{1}^{\T}]^{\T}$.

Then, using  (\ref{eq:Tmatrix}) and  (\ref{eq:CCAeucHess}) we have
\[
(\Sigma^{-\nicehalf}\nabla^{2}\bar{f}(\z^{\star})\Sigma^{-\nicehalf}+\sigma_{1}\cdot\matI_{d})\Sigma^{\nicehalf}\left[\begin{array}{c}
\u_{1}\\
\v_{1}
\end{array}\right]=\left(\left[\begin{array}{cc}
 & -\matR\\
-\matR^{\T}
\end{array}\right]+\sigma_{1}\matI_{d}\right)\left[\begin{array}{c}
\phi\\
\psi
\end{array}\right]=0\ ,
\]
where the last equality follows from the fact that $\left[\begin{array}{cc}
 & \matR\\
\matR^{\T}
\end{array}\right]$ is the augmented matrix associated with $\matR$, so $\left[\begin{array}{c}
\phi\\
\psi
\end{array}\right]$, which has the dominant left and right singular vectors stacked,
is the eigenvalue corresponding to the largest eigenvalue $\sigma_{1}$
of the augmented matrix. Similarly, since the vector $\left[\begin{array}{c}
\phi\\
-\psi
\end{array}\right]$ is the eigenvector corresponding to the smallest eigenvalue $-\sigma_{1}$
of the augmented matrix, then
\[
(\Sigma^{-\nicehalf}\nabla^{2}\bar{f}(\z^{\star})\Sigma^{-\nicehalf}+\sigma_{1}\cdot\matI_{d})\Sigma^{\nicehalf}\left[\begin{array}{c}
\u_{1}\\
-\v_{1}
\end{array}\right]=\left(\left[\begin{array}{cc}
 & -\matR\\
-\matR^{\T}
\end{array}\right]+\sigma_{1}\matI_{d}\right)\left[\begin{array}{c}
\phi\\
-\psi
\end{array}\right]=2\sigma_{1}\left[\begin{array}{c}
\phi\\
-\psi
\end{array}\right]\ .
\]

Finally, the minimal and the maximal eigenvalues of the matrix $\Sigma^{-\nicehalf}\nabla^{2}\bar{f}(\z^{\star})\Sigma^{-\nicehalf}+\sigma_{1}\cdot\matI_{d}$
in the space of vectors $\xi$ such that $\Sigma^{-\nicehalf}\xi\in T_{\z^{\star}}\elpCCA$
are $\sigma_{1}-\sigma_{2}$ and $\sigma_{1}+\sigma_{2}$ correspondingly.
Thus, 
\[
\lambda_{\max}(\hess{f(\z^{\star})})=\max_{\mat 0\neq\xi_{\z^{\star}}\in T_{\z^{\star}}\elpCCA}q(\xi_{\z^{\star}})=\sigma_{1}+\sigma_{2}>0\ ,
\]
and,
\[
\lambda_{\min}(\hess{f(\z^{\star})})=\min_{\mat 0\neq\xi_{\z^{\star}}\in T_{\z^{\star}}\elpCCA}q(\xi_{\z^{\star}})=\sigma_{1}-\sigma_{2}>0\ ,
\]
The condition number for the case $\matM_{\z^{\star}}=\Sigma$ is
obtained by dividing the last two quantities.

If $\matM_{\z^{\star}}\approx\Sigma$, we can bound the smallest and
largest eigenvalues of the Riemannian Hessian at $\z^{\star}$ by
\begin{eqnarray*}
\lambda_{\min}(\hess{f(\z^{\star})}) & \geq & \min_{0\neq\eta_{\z^{\star}}\in T_{\z^{\star}}\elpCCA}\frac{\eta_{\z^{\star}}^{\T}\left(\nabla^{2}\bar{f}(\z^{\star})+\sigma_{1}\cdot\Sigma\right)\eta_{\z^{\star}}}{\eta_{\z^{\star}}^{\T}\Sigma\eta_{\z^{\star}}}\cdot\min_{\eta_{\z^{\star}}\neq0}\frac{\eta_{\z^{\star}}^{\T}\Sigma\eta_{\z^{\star}}}{\eta_{\z^{\star}}^{\T}\matM_{\z^{\star}}\eta_{\z^{\star}}}\\
 & = & \lambda_{\min}(\Sigma,\matM_{\z^{\star}})\cdot(\sigma_{1}-\sigma_{2})\ ,
\end{eqnarray*}
and
\begin{eqnarray*}
\lambda_{\max}(\hess{f(\z^{\star})}) & \leq & \max_{0\neq\eta_{\z^{\star}}\in T_{\z^{\star}}\elpCCA}\frac{\eta_{\z^{\star}}^{\T}\left(\nabla^{2}\bar{f}(\z^{\star})+\sigma_{1}\cdot\Sigma\right)\eta_{\z^{\star}}}{\eta_{\z^{\star}}^{\T}\Sigma\eta_{\z^{\star}}}\cdot\max_{\eta_{\z^{\star}}\neq0}\frac{\eta_{\z^{\star}}^{\T}\Sigma\eta_{\z^{\star}}}{\eta_{\z^{\star}}^{\T}\matM_{\z^{\star}}\eta_{\z^{\star}}}\\
 & = & \lambda_{\max}(\Sigma,\matM_{\z^{\star}})\cdot(\sigma_{1}+\sigma_{2})\ .
\end{eqnarray*}
Finally, we get
\[
\kappa(\hess{f(\z^{\star})})=\frac{\lambda_{\max}(\hess{f(\z^{\star})})}{\lambda_{\min}(\hess{f(\z^{\star})})}\leq\frac{\sigma_{1}+\sigma_{2}}{\sigma_{1}-\sigma_{2}}\cdot\kappa\left(\matB,\matM_{\z^{\star}}\right)\ .
\]
\end{proof}
We now illustrate the effect of the preconditioning scheme $\z\mapsto\matM_{\z}$
numerically. In our experiments, we use six metric choices with constant
matrices, i.e., $\matM_{\z}\coloneqq\matM$ independent of $\z\in\elpCCA$:
the trivial choice of a unit matrix $\matM=\mat I_{d}$, the standard
but expensive choice $\matM=\Sigma$ which achieves the optimal bound
according to Lemma \ref{lem:CCA_cond_num}, and four approximations
of $\Sigma$ via the (exact) sketched preconditioning strategy described
by Gonen et al. \cite{GOS16}, which we term as \emph{Dominant Subspace
Preconditioning}.

Dominant Subspace Preconditioning was originally designed for ridge
regression to speed up Stochastic Variance Reduced Gradient via an
approximation of the empirical correlation matrix. In our experiments
we use this preconditioning strategy to approximate $\Sigma_{\x\x}$
and $\Sigma_{\y\y}$. The approximation is done as follows: suppose
$\matA=\hat{\matX}\hat{\matX}^{\T}\in\R^{d\times d}$ be some positive
semi-definite matrix, and let $\hat{\matX}=\matU\mat{\Lambda}^{\nicehalf}\matV^{\T}$
be an SVD decomposition of $\hat{\matX}$ such that $\matA=\matU\Lambda\matU^{\T}$
is an eigendecomposition, with the diagonal entries in $\Lambda$
sorted in descending order. Given $k$, let us denote by $\matU_{k}$
the first $k$ columns of $\matU$, $\Lambda_{k}$ denote the leading
$k\times k$ minor of $\Lambda$, and $\lambda_{k}$ the $k$-th largest
eigenvalue of $\matA$. The $k$-dominant subspace preconditioner
of $\matA+\lambda\matI_{d}$ is $\matU_{k}(\Lambda_{k}-\lambda_{k}\matI)\matU_{k}^{\T}+(\lambda_{k}+\lambda)\matI_{d}$.
The dominant subspace can be found using a sparse SVD solver (we use
\noun{MATLAB}'s svds). Moreover, its inverse can be easily computed
using the formula
\[
\matU_{k}(\Lambda_{k}+\lambda\matI)^{-1}\matU_{k}^{\T}+\frac{1}{\lambda_{k}+\lambda}(\matI_{d}-\matU_{k}\matU_{k}^{\T}).
\]

The experiments are performed with the \noun{MEDIANILL}\footnote{Datasets were downloaded for libsvm's website: https://www.csie.ntu.edu.tw/\textasciitilde cjlin/libsvmtools/datasets/}
dataset where the dimensions are $n=43907$, $d_{\x}=120$, and $d_{\y}=101$.
The implementation uses \noun{Manopt} which is a \noun{MATLAB} library
that implements some Riemannian optimization algorithms \cite{manopt}.
In Fig. \ref{fig:Results-for-CCA} the left graph presents suboptimality
vs. iteration count for Riemannian CG, and the right graph presents
suboptimality vs. products with the data matrices for Riemannian trust-region.
Note that in Riemannian trust-region, different iterations do a variable
amount of passes over the data, thus, this is the dominant cost of
the trust-region method. The graphs in Fig. \ref{fig:Results-for-CCA}
demonstrate that the choice $\matM=\Sigma$ leads to the lowest iteration
count. This observation is also supported by the condition number
of the Riemannian Hessian at the optimum. We evaluated it using \noun{Manopt},
and indeed, the lowest condition number, $4.03$, is achieved when
$\matM=\Sigma$, and the highest, $60.2$, for $\matM=\mat I_{d}$.

\section{Conclusions}

In this paper, we developed the preconditioned geometric components
for optimization on the generalized Stiefel manifold. The main mechanism
for introducing a preconditioner is via the Riemannian metric. The
technique can be used to precondition any underlying Riemannian optimization
method. Our method can also be applied to constraints which are described
by the product of two or more generalized Stiefel manifolds. We demonstrated
our method both theoretically and numerically on the problem of computing
the dominant canonical correlation. As part of developing the related
geometrical components of the generalized Stiefel manifold equipped
with a non standard Riemannian metric, we evaluate the costs of computing
these components and relate the preconditioner to asymptotic convergence
via the condition number of the Riemannian Hessian at the optimum.

In a sense, this paper presents only part of the picture. While it
presents a methodology for building preconditioned algorithms for
optimization with generalized orthogonality constraints, it does not
explains how to build effective preconditioners to be used in conjunction
with those algorithms, and we leave it for future work. Additional
research directions include addressing other constraints using similar
ideas, e.g., fixed-rank matrices, products of different types of manifolds,
quotient manifolds, etc.
\begin{acknowledgement*}
The authors thank Bart Vandereycken for useful discussions. The authors
would also like to thank the referees for their valuable comments.
This research was supported by the Israel Science Foundation (grant
no. 1272/17).
\end{acknowledgement*}
\bibliographystyle{elsarticle-num}
\bibliography{ellipsoid}

\begin{thebibliography}{10}
\expandafter\ifx\csname url\endcsname\relax
  \def\url#1{\texttt{#1}}\fi
\expandafter\ifx\csname urlprefix\endcsname\relax\def\urlprefix{URL }\fi
\expandafter\ifx\csname href\endcsname\relax
  \def\href#1#2{#2} \def\path#1{#1}\fi

\bibitem{hotelling1936relations}
H.~Hotelling, {Relations between two sets of variates}, Biometrika 28~(3/4)
  (1936) 321--377.

\bibitem{fisher1936use}
R.~A. Fisher, {The use of multiple measurements in taxonomic problems}, Ann.
  Eugen. 7~(2) (1936) 179--188.

\bibitem{bjorck1973numerical}
A.~Bj\"{o}rck, G.~H. Golub, {Numerical methods for computing angles between
  linear subspaces}, Math. Comput. 27~(123) (1973) 579--594.

\bibitem{EAS98}
A.~Edelman, T.~Arias, S.~Smith,
  \href{https://doi.org/10.1137/S0895479895290954}{{The Geometry of Algorithms
  with Orthogonality Constraints}}, SIAM J. Matrix Anal. Appl. 20~(2) (1998)
  303--353.
\newblock \href
  {http://arxiv.org/abs/https://doi.org/10.1137/S0895479895290954}
  {\path{arXiv:https://doi.org/10.1137/S0895479895290954}}, \href
  {https://doi.org/10.1137/S0895479895290954}
  {\path{doi:10.1137/S0895479895290954}}.
\newline\urlprefix\url{https://doi.org/10.1137/S0895479895290954}

\bibitem{AMS09}
P.~A. Absil, R.~Mahony, R.~Sepulchre, {Optimization Algorithms on Matrix
  Manifolds}, Princeton University Press, 2009.

\bibitem{boumal2020intromanifolds}
N.~Boumal, \href{http://www.nicolasboumal.net/book}{{An introduction to
  optimization on smooth manifolds}}, Available online (May 2020).
\newline\urlprefix\url{http://www.nicolasboumal.net/book}

\bibitem{manopt}
N.~Boumal, B.~Mishra, P.~A. Absil, R.~Sepulchre,
  \href{http://www.manopt.org}{{{Manopt}, a {M}atlab Toolbox for Optimization
  on Manifolds}}, J Mach Learn Res. 15 (2014) 1455--1459.
\newline\urlprefix\url{http://www.manopt.org}

\bibitem{MS16}
B.~Mishra, R.~Sepulchre, \href{http://dx.doi.org/10.1137/140970860}{{Riemannian
  Preconditioning}}, SIAM J. Optim. 26~(1) (2016) 635--660.
\newblock \href {http://arxiv.org/abs/http://dx.doi.org/10.1137/140970860}
  {\path{arXiv:http://dx.doi.org/10.1137/140970860}}, \href
  {https://doi.org/10.1137/140970860} {\path{doi:10.1137/140970860}}.
\newline\urlprefix\url{http://dx.doi.org/10.1137/140970860}

\bibitem{luenberger1972gradient}
D.~G. Luenberger, {The gradient projection method along geodesics}, Manage Sci.
  18~(11) (1972) 620--631.

\bibitem{gabay1982minimizing}
D.~Gabay, {Minimizing a differentiable function over a differential manifold},
  J Optim Theory Appl. 37~(2) (1982) 177--219.

\bibitem{smith1994optimization}
S.~T. Smith, {Optimization Techniques on Riemannian Manifolds}, Fields
  institute communications 3~(3) (1994) 113--135.

\bibitem{wen2013feasible}
Z.~Wen, W.~Yin, {A feasible method for optimization with orthogonality
  constraints}, Math Program 142~(1) (2013) 397--434.

\bibitem{zhu2017riemannian}
X.~Zhu, {A {Riemannian} conjugate gradient method for optimization on the
  {Stiefel} manifold}, Comput Optim Appl. 67~(1) (2017) 73--110.

\bibitem{li2020efficient}
J.~Li, F.~Li, S.~Todorovic,
  \href{https://openreview.net/forum?id=HJxV-ANKDH}{{Efficient Riemannian
  Optimization on the Stiefel Manifold via the Cayley Transform}}, in:
  International Conference on Learning Representations, 2020.
\newline\urlprefix\url{https://openreview.net/forum?id=HJxV-ANKDH}

\bibitem{sato2019cholesky}
H.~Sato, K.~Aihara, {Cholesky QR-based retraction on the generalized Stiefel
  manifold}, Comput Optim Appl. 72~(2) (2019) 293--308.

\bibitem{kaneko2012empirical}
T.~Kaneko, S.~Fiori, T.~Tanaka, {Empirical arithmetic averaging over the
  compact {Stiefel} manifold}, IEEE Trans. Signal 61~(4) (2012) 883--894.

\bibitem{beck2009fast}
A.~Beck, M.~Teboulle, {A fast iterative shrinkage-thresholding algorithm for
  linear inverse problems}, SIAM J. Imaging Sci. 2~(1) (2009) 183--202.

\bibitem{chen2018proximal}
S.~Chen, S.~Ma, A.~Man-Cho~So, T.~Zhang,
  \href{https://doi.org/10.1137/18M122457X}{{Proximal Gradient Method for
  Nonsmooth Optimization over the {Stiefel} Manifold}}, SIAM J. Optim. 30~(1)
  (2020) 210--239.
\newblock \href {http://arxiv.org/abs/https://doi.org/10.1137/18M122457X}
  {\path{arXiv:https://doi.org/10.1137/18M122457X}}, \href
  {https://doi.org/10.1137/18M122457X} {\path{doi:10.1137/18M122457X}}.
\newline\urlprefix\url{https://doi.org/10.1137/18M122457X}

\bibitem{chen2019alternating}
S.~Chen, S.~Ma, L.~Xue, H.~Zou,
  \href{https://doi.org/10.1287/ijoo.2019.0032}{{An Alternating Manifold
  Proximal Gradient Method for Sparse Principal Component Analysis and Sparse
  Canonical Correlation Analysis}}, INFORMS J Optim. 2~(3) (2020) 192--208.
\newblock \href {http://arxiv.org/abs/https://doi.org/10.1287/ijoo.2019.0032}
  {\path{arXiv:https://doi.org/10.1287/ijoo.2019.0032}}, \href
  {https://doi.org/10.1287/ijoo.2019.0032} {\path{doi:10.1287/ijoo.2019.0032}}.
\newline\urlprefix\url{https://doi.org/10.1287/ijoo.2019.0032}

\bibitem{huang2019extending}
W.~Huang, K.~Wei, {An Extension of {FISTA} to {Riemannian} Optimization for
  Sparse {PCA}}, arXiv preprint arXiv:1909.05485 (2019).

\bibitem{yger2012adaptive}
F.~Yger, M.~Berar, G.~Gasso, A.~Rakotomamonjy,
  \href{http://dl.acm.org/citation.cfm?id=3042573.3042615}{{Adaptive
  {Canonical} {Correlation} {Analysis} Based on Matrix Manifolds}}, in:
  Proceedings of the 29th International Coference on International Conference
  on Machine Learning, ICML'12, Omnipress, USA, 2012, pp. 299--306.
\newline\urlprefix\url{http://dl.acm.org/citation.cfm?id=3042573.3042615}

\bibitem{pechen2008control}
A.~Pechen, D.~Prokhorenko, R.~Wu, H.~Rabitz, {Control landscapes for two-level
  open quantum systems}, J. Phys. A Math. 41~(4) (2008) 045205.

\bibitem{oza2009optimization}
A.~Oza, A.~Pechen, J.~Dominy, V.~Beltrani, K.~Moore, H.~Rabitz, {Optimization
  search effort over the control landscapes for open quantum systems with
  Kraus-map evolution}, J. Phys. A Math. 42~(20) (2009) 205305.

\bibitem{sato2013complex}
H.~Sato, T.~Iwai, {A complex singular value decomposition algorithm based on
  the {R}iemannian {N}ewton method}, in: Decision and Control ({CDC}), 2013
  {IEEE} 52nd Annual Conference on, IEEE, 2013, pp. 2972--2978.

\bibitem{sato2014Riemannian}
H.~Sato, {{R}iemannian conjugate gradient method for complex singular value
  decomposition problem}, in: Decision and Control ({CDC}), 2014 {IEEE} 53rd
  Annual Conference on, IEEE, 2014, pp. 5849--5854.

\bibitem{udriste2013convex}
C.~Udriste, {Convex functions and optimization methods on {Riemannian}
  manifolds}, Vol. 297, Springer Science \& Business Media, 2013.

\bibitem{ngo2012scaled}
T.~Ngo, Y.~Saad, {Scaled gradients on Grassmann manifolds for matrix
  completion}, in: Adv Neural Inf Process Syst., 2012, pp. 1412--1420.

\bibitem{mishra2014r3mc}
B.~Mishra, R.~Sepulchre, {R3MC: A Riemannian three-factor algorithm for
  low-rank matrix completion}, in: Decision and Control (CDC), 2014 IEEE 53rd
  Annual Conference on, IEEE, 2014, pp. 1137--1142.

\bibitem{shi2016low}
Y.~Shi, J.~Zhang, K.~B. Letaief, {Low-rank matrix completion for topological
  interference management by Riemannian pursuit}, IEEE Trans. Wirel. 15~(7)
  (2016) 4703--4717.

\bibitem{zhou2016riemannian}
T.~Zhou, H.~Qian, Z.~Shen, C.~Zhang, C.~Xu,
  \href{https://doi.org/10.24963/ijcai.2017/495}{{Tensor Completion with Side
  Information: A {Riemannian} Manifold Approach}}, in: Proceedings of the
  Twenty-Sixth International Joint Conference on Artificial Intelligence,
  {IJCAI-17}, 2017, pp. 3539--3545.
\newblock \href {https://doi.org/10.24963/ijcai.2017/495}
  {\path{doi:10.24963/ijcai.2017/495}}.
\newline\urlprefix\url{https://doi.org/10.24963/ijcai.2017/495}

\bibitem{vandereycken2010riemannian}
B.~Vandereycken, S.~Vandewalle, {A Riemannian optimization approach for
  computing low-rank solutions of Lyapunov equations}, SIAM J. Matrix Anal.
  Appl. 31~(5) (2010) 2553--2579.

\bibitem{mor2020solving}
U.~Mor, H.~Avron, {Solving Trust Region Subproblems Using {Riemannian}
  Optimization}, arXiv preprint arXiv:2010.07547 (2020).

\bibitem{kressner2016preconditioned}
D.~Kressner, M.~Steinlechner, B.~Vandereycken, {Preconditioned low-rank
  Riemannian optimization for linear systems with tensor product structure},
  SIAM J Sci Comput. 38~(4) (2016) A2018--A2044.

\bibitem{bock1987randwertproblemmethoden}
H.-G. Bock, {Randwertproblemmethoden zur Parameteridentifizierung in systemen
  nichtlinearer Differentialgleichungen}, Bonner mathematische Schriften,
  Rheinische Friedrich-Wilhelm Universit{\"a}t, Academic Dissertation 16795956,
  1987 (in German).

\bibitem{Golub:1996:MC:248979}
G.~H. Golub, C.~F. Van~Loan, {Matrix Computations (4rd Ed.)}, Johns Hopkins
  University Press, Baltimore, MD, USA, 2013.

\bibitem{trefethen1997numerical}
L.~N. Trefethen, D.~Bau~III, {Numerical linear algebra}, Vol.~50, SIAM, 1997.

\bibitem{lezcano2019trivializations}
M.~Lezcano~Casado, {Trivializations for gradient-based optimization on
  manifolds}, Adv Neural Inf Process Syst. 32 (2019).

\bibitem{criscitiello2020accelerated}
C.~Criscitiello, N.~Boumal, {An accelerated first-order method for non-convex
  optimization on manifolds}, arXiv preprint arXiv:2008.02252 (2020).

\bibitem{bento2017iteration}
G.~C. Bento, O.~P. Ferreira, J.~G. Melo, {Iteration-complexity of gradient,
  subgradient and proximal point methods on Riemannian manifolds}, J Optim
  Theory Appl. 173~(2) (2017) 548--562.

\bibitem{ferreira2002proximal}
O.~P. Ferreira, P.~R. Oliveira, {Proximal point algorithm on {Riemannian}
  manifolds}, Optim. 51~(2) (2002) 257--270.

\bibitem{zhu2020riemannian}
X.~Zhu, H.~Sato, {Riemannian conjugate gradient methods with inverse
  retraction}, Comput Optim Appl. 77~(3) (2020) 779--810.

\bibitem{horn2012matrix}
R.~A. Horn, C.~R. Johnson, {Matrix Analysis, 2nd Ed.}, Cambridge University
  Press, 2012.

\bibitem{bauer1960norms}
F.~L. Bauer, C.~T. Fike, {Norms and exclusion theorems}, Numer Math (Heidelb)
  2~(1) (1960) 137--141.

\bibitem{absil2013extrinsic}
P.~A. Absil, R.~Mahony, J.~Trumpf, {An Extrinsic Look at the Riemannian
  Hessian}, in: Geometric Science of Information, Springer, 2013, pp. 361--368.

\bibitem{bhatia1997and}
R.~Bhatia, P.~Rosenthal, {How and why to solve the operator equation
  $AX-XB=Y$}, Bull. London Math. Soc. 29~(1) (1997) 1--21.

\bibitem{boumal2019global}
N.~Boumal, P.~A. Absil, C.~Cartis, {Global rates of convergence for nonconvex
  optimization on manifolds}, IMA J. Numer. Anal. 39~(1) (2019) 1--33.

\bibitem{sun2010scalable}
L.~Sun, B.~Ceran, J.~Ye, {A scalable two-stage approach for a class of
  dimensionality reduction techniques}, in: Proceedings of the 16th ACM SIGKDD
  international conference on Knowledge discovery and data mining, ACM, 2010,
  pp. 313--322.

\bibitem{chaudhuri2009multi}
K.~Chaudhuri, S.~M. Kakade, K.~Livescu, K.~Sridharan, {Multi-view clustering
  via canonical correlation analysis}, in: Proceedings of the 26th {A}nnual
  International Conference on Machine Learning (ICML), ACM, 2009, pp. 129--136.

\bibitem{dhillon2011multi}
P.~Dhillon, D.~P. Foster, L.~H. Ungar, {Multi-view learning of word embeddings
  via {CCA}}, in: Adv Neural Inf Process Syst., 2011, pp. 199--207.

\bibitem{dhillon2012two}
P.~S. Dhillon, J.~Rodu, D.~P. Foster, L.~H. Ungar,
  \href{http://dl.acm.org/citation.cfm?id=3042573.3042586}{{Two Step {CCA}: A
  New Spectral Method for Estimating Vector Models of Words}}, in: Proceedings
  of the 29th International Conference on International Conference on Machine
  Learning, ICML'12, Omnipress, USA, 2012, pp. 67--74.
\newline\urlprefix\url{http://dl.acm.org/citation.cfm?id=3042573.3042586}

\bibitem{su2012discriminant}
Y.~Su, Y.~Fu, X.~Gao, Q.~Tian, {Discriminant learning through multiple
  principal angles for visual recognition}, IEEE Trans Image Process. 21~(3)
  (2012) 1381--1390.

\bibitem{kim2007discriminative}
T.-K. Kim, J.~Kittler, R.~Cipolla, {Discriminative learning and recognition of
  image set classes using canonical correlations}, IEEE Trans Pattern Anal Mach
  Intell. 29~(6) (2007) 1005--1018.

\bibitem{golub1995canonical}
G.~H. Golub, H.~Zha, {The canonical correlations of matrix pairs and their
  numerical computation}, in: Linear Algebra for Signal Processing, Springer,
  1995, pp. 27--49.

\bibitem{GOS16}
A.~Gonen, F.~Orabona, S.~Shalev-Shwartz,
  \href{http://dl.acm.org/citation.cfm?id=3045390.3045538}{{Solving Ridge
  Regression Using Sketched Preconditioned {SVRG}}}, in: Proceedings of the
  33rd International Conference on International Conference on Machine Learning
  - Volume 48, ICML'16, JMLR.org, 2016, pp. 1397--1405.
\newline\urlprefix\url{http://dl.acm.org/citation.cfm?id=3045390.3045538}

\bibitem{bartels1972solution}
R.~H. Bartels, G.~W. Stewart,
  \href{https://doi.org/10.1145/361573.361582}{{Solution of the Matrix Equation
  $AX + XB = C$}}, Commun. ACM 15~(9) (1972) 820--826.
\newblock \href {https://doi.org/10.1145/361573.361582}
  {\path{doi:10.1145/361573.361582}}.
\newline\urlprefix\url{https://doi.org/10.1145/361573.361582}

\bibitem{absil2009all}
P.~A. Absil, J.~Trumpf, R.~Mahony, B.~Andrews, {All roads lead to {N}ewton:
  Feasible second-order methods for equality-constrained optimization},
  Technical Report UCL-INMA-2009.024 (2009).

\bibitem{von1937some}
J.~Von~Neumann, {Some matrix-inequalities and metrization of matric-space.
  Tomsk Univ. Rev. 1 (1937) 286--300}.

\end{thebibliography}

\appendix

\section{Further Details on the Preconditioned Geometric Components}

In this section we elaborate on the derivations of the Riemannian
components that appear in Section \ref{sec:precond-geometry}. Our
main contribution is the metric dependent components in Subsection
\ref{subsec:Metric-Related-Notions-1}. The metric independent components
are included for completeness.

\subsection{\label{subsec:Metric-Independent-Notions-1}Metric Independent Notions}

We begin with the metric independent notions that appear in Subsection
\ref{subsec:metric-independent}. Recall that the tangent space has
two common characterizations. The first characterization
\begin{equation}
T_{\matX}\text{St}_{\matB}(p,d)=\left\{ \matZ\in\R^{d\times p}\,:\,\matZ^{\T}\matB\matX+\matX^{\T}\matB\matZ=0_{p}\right\} ,\label{eq:tangent_st_B_appendix}
\end{equation}
is based on the Submersion Theorem \cite[Proposition 3.3.3]{AMS09}.
$\text{St}_{\matB}(p,d)$ is the kernel of the mapping $F(\matX)=\matX^{\T}\matB\matX-\matI_{p}$,
i.e., $\text{St}_{\matB}(p,d)=F^{-1}(0_{p})$. This mapping is a submersion
since the rank of $F$ is $p(p+1)/2$ (i.e., $F$ is full rank); indeed,
the rank of $F$ is determined by the range of $DF(\matX)[\cdot]:\R^{d\times p}\to{\cal {\cal S}_{\text{sym}}}(p)$.
For every $\hat{\mat Z}\in{\cal {\cal S}_{\text{sym}}}(p)$, the matrix
$\matZ=\frac{1}{2}\matX\hat{Z}\in\R^{d\times p}$ satisfies $DF(\matX)[\matZ]=\hat{Z}$.
According to \cite[Proposition 3.3.3]{AMS09} then $\text{St}_{\matB}(p,d)$
is an embedded submanifold of $\R^{d\times p}$, and its dimension
is $dp-\frac{p(p+1)}{2}$.

The second characterization is:
\begin{equation}
T_{\matX}\text{St}_{\matB}(p,d)=\left\{ \matZ=\matX\Omega+\matX_{\matB\perp}K\in\R^{d\times p}\,:\,\Omega\in{\cal {\cal S}_{\text{skew}}}(p),\ K\in\R^{(d-p)\times p}\right\} ,\label{eq:tangentstiefel_st_B_2_appendix}
\end{equation}
where $\Omega$ is a skew-symmetric matrix (i.e., $\,\Omega^{\T}=-\Omega$),
$K$ is arbitrary, and $\matX_{\matB\perp}\in\R^{d\times(d-p)}$ satisfies
that its columns are an orthonormal basis for the orthogonal complement
of the column space of $\matX$ with respect to the matrix $\matB$,
i.e., $\matX_{\matB\perp}^{\T}\matB\matX_{\matB\perp}=\matI_{d-p}$,
and $\matX_{\matB\perp}^{\T}\matB\matX=0_{(d-p)\times p}$. The dimension
of the space defined in ~(\ref{eq:tangentstiefel_st_B_2_appendix})
is $p(p-1)/2+p(d-p)=dp-p(p+1)/2$. Both characterizations of $T_{\matX}\text{St}_{\matB}(p,d)$,
~(\ref{eq:tangent_st_B_appendix}) and ~(\ref{eq:tangentstiefel_st_B_2_appendix}),
are equal. Indeed, every $\matZ\in\R^{d\times p}$ can be represented
by $\matX\Omega+\matX_{\matB\perp}K$ for arbitrary $\Omega\in\R^{p\times p}$
and $K\in\R^{(d-p)\times p}$ ($dp$ degrees of freedom), where the
columns of $\mat X$ and $\matX_{\matB\perp}$ are linearly independent,
thus each of the columns of $\matZ$ can be any vector in $\R^{d}$,
and $\matZ$ any matrix in $\R^{d\times p}$. Suppose $\matZ$ satisfies
~(\ref{eq:tangent_st_B_appendix}), then $\Omega^{\T}=-\Omega$,
so that $\matZ$ belongs to the set defined in ~(\ref{eq:tangentstiefel_st_B_2_appendix}).
Thus, the set defined in ~(\ref{eq:tangent_st_B_appendix}) is a
subset (subspace) of the set defined in ~(\ref{eq:tangentstiefel_st_B_2_appendix}).
Finally, since both the sets defined in ~(\ref{eq:tangent_st_B_appendix})
and ~(\ref{eq:tangentstiefel_st_B_2_appendix}) are subspaces of
$T_{\matX}\R^{d\times p}\simeq\R^{d\times p}$, and both are with
the same dimension we get that ~(\ref{eq:tangent_st_B_appendix})
and ~(\ref{eq:tangentstiefel_st_B_2_appendix}) are equal.

In this article, we consider the use of three retractions mappings:
\begin{equation}
R_{\matX}^{\text{polar}}(\xi_{\matX})\coloneqq(\matX+\xi_{\matX})(\mat I_{p}+\xi_{\matX}^{\T}B\xi_{\matX})^{-\nicehalf}\label{eq:stiefelpolarretraction_appendix}
\end{equation}
\begin{equation}
R_{\matX}^{\text{QR}}(\xi_{\matX})\coloneqq\qfm{\matX+\xi_{\matX}}{\matB}=\matB^{-\nicehalf}\qf{\matB^{\nicehalf}\left(\matX+\xi_{\matX}\right)}\label{eq:stiefelqrretraction_appendix}
\end{equation}
\begin{equation}
R_{\matX}^{\text{Cayley}}(\xi_{\matX})\coloneqq(\matI_{d}-\frac{1}{2}\matW(\xi_{\matX}))^{-1}(\matI_{d}+\frac{1}{2}\matW(\xi_{\matX}))\matX\label{eq:stiefelcayleyretraction_appendix}
\end{equation}
where 
\[
\matW(\xi_{\matX})\coloneqq(\matI_{d}-\frac{1}{2}\matX\matX^{\T}\matB)\xi_{\matX}\matX^{\T}\matB-\matX\xi_{\matX}^{\T}(\matI_{d}-\frac{1}{2}\matB\matX\matX^{\T})\matB.
\]
The cost of computing the polar-based retraction,  (\ref{eq:stiefelpolarretraction_appendix}),
is $O\left(T_{\matB}p+dp^{2}\right)$ where $T_{\matB}$ is the cost
of computing the product of $\matB$ with a vector. This is evident
from the formulas since none of the operations require forming $\matB$,
but instead require taking product of $\matB$ with matrices, finding
the inverse of a square root of a $p\times p$ matrix, multiplying
a $d\times p$ matrix by a $p\times p$ matrix, and multiplying a
$p\times d$ matrix by a $d\times p$ matrix. This is also mentioned
in \cite[Section 3.2]{sato2019cholesky}. The cost of computing the
QR-based retraction, (\ref{eq:stiefelqrretraction_appendix}), is
also $O\left(T_{\matB}p+dp^{2}\right)$. This is shown in \cite[Section 3.2]{sato2019cholesky}.
Though, in \cite{sato2019cholesky}, it is claimed that for large
$p$ ($p\leq d)$ the QR-based retraction has an advantage in computational
costs compared to the polar-based retraction, since the eigenvalue
decomposition of $(\matX+\xi_{\matX})^{\T}\matB(\matX+\xi_{\matX})$
(or SVD decomposition of $\matX+\xi_{\matX}$) can be replaced with
a Cholesky decomposition of the same matrix. The cost of computing
the Cayley transform based retraction, (\ref{eq:stiefelcayleyretraction_appendix}),
is $O\left(T_{\matB}p+dp^{2}\right)$ which follows using the Sherman-Morrison-Woodbury
formula as described in Subsection \ref{subsec:metric-independent}.
Another approach suggested in \cite{li2020efficient} by Li et al.
is to use a fixed point method to approximate the retraction.

The retraction in ~(\ref{eq:stiefelqrretraction_appendix}) is proven
to be indeed a retraction mapping in \cite[Theorem 3.1]{sato2019cholesky}.
For the retraction in ~(\ref{eq:stiefelpolarretraction_appendix}),
though we found the equation in the literature, we could not find
a formal argument that it is a retraction. Therefore, we show this
by showing that it meets the conditions in \cite[Definition 4.1.1]{AMS09}.
The first condition of \cite[Definition 4.1.1]{AMS09} is that $R_{\matX}(\mat 0_{\matX})=\matX$,
and it indeed holds since $R_{\matX}^{\text{polar}}(\mat 0_{\matX})=(\matX+\mat 0_{\matX})(\mat I_{p}+\mat 0_{\matX}^{\T}B\mat 0_{\matX})^{-\nicehalf}=\matX$.
The second condition of \cite[Definition 4.1.1]{AMS09} is that $\text{D}R_{\x}(0_{\x})=\text{id}_{T_{\matX}\text{St}_{\matB}(p,d)}$,
where $\text{id}_{T_{\matX}\text{St}_{\matB}(p,d)}$ denotes the identity
mapping on $T_{\matX}\text{St}_{\matB}(p,d)$. This condition is equivalent
to the condition that for every vector $\xi_{\matX}\in T_{\matX}\text{St}_{\matB}(p,d)$
we have $\left.\frac{\text{d}}{\text{dt}}R_{\mat X}(t\xi_{\matX})\right|_{t=0}=\xi_{\matX}$.
Denote by $\lambda_{1},...,\lambda_{p}\geq0$ the eigenvalues of $\xi_{\matX}^{\T}B\xi_{\matX}$,
then 
\[
(\mat I_{p}+t^{2}\xi_{\matX}^{\T}B\xi_{\matX})^{-\nicehalf}=\matQ\left(\begin{array}{ccc}
\frac{1}{\sqrt{1+t^{2}\lambda_{1}}}\\
 & \ddots\\
 &  & \frac{1}{\sqrt{1+t^{2}\lambda_{p}}}
\end{array}\right)\matQ^{\T},
\]
where $\matQ$ is an orthogonal matrix that diagonalizes $\xi_{\matX}^{\T}B\xi_{\matX}$.
Then,
\begin{eqnarray*}
\left.\frac{\text{d}}{\text{dt}}R_{\matX}^{\text{polar}}(t\xi_{\matX})\right|_{t=0} & = & \left.\frac{\text{d}}{\text{dt}}\left[(\matX+t\xi_{\matX})(\mat I_{p}+t^{2}\xi_{\matX}^{\T}B\xi_{\matX})^{\nicehalf}\right]\right|_{t=0}=
\end{eqnarray*}
\begin{eqnarray*}
 & = & \left.\frac{\text{d}}{\text{dt}}\left[(\matX+t\xi_{\matX})\matQ\left(\begin{array}{ccc}
\frac{1}{\sqrt{1+t^{2}\lambda_{1}}}\\
 & \ddots\\
 &  & \frac{1}{\sqrt{1+t^{2}\lambda_{p}}}
\end{array}\right)\matQ^{\T}\right]\right|_{t=0}
\end{eqnarray*}
\begin{eqnarray*}
=\left.\xi_{\matX}\matQ\left(\begin{array}{ccc}
\frac{1}{\sqrt{1+t^{2}\lambda_{1}}}\\
 & \ddots\\
 &  & \frac{1}{\sqrt{1+t^{2}\lambda_{p}}}
\end{array}\right)\matQ^{\T}-(\matX+t\xi_{\matX})\matQ\left(\begin{array}{ccc}
\frac{t\lambda_{1}}{\left(1+t^{2}\lambda_{1}\right)^{1.5}}\\
 & \ddots\\
 &  & \frac{t\lambda_{p}}{\left(1+t^{2}\lambda_{p}\right)^{1.5}}
\end{array}\right)\matQ^{\T}\right|_{t=0} & = & \xi_{\matX}.
\end{eqnarray*}

Similarly, the retraction in ~(\ref{eq:stiefelcayleyretraction_appendix})
is also proven to be a retraction mapping in \cite[Eq. (4)]{li2020efficient}
for the Stiefel manifold. In order to generalize it to the generalized
Stiefel manifold, we show it meets the conditions in \cite[Definition 4.1.1]{AMS09}.
For the first condition we have $\matW(\mat 0_{\matX})=\mat 0_{d}$,
thus
\[
R_{\matX}^{\text{Cayley}}(\mat 0_{\matX})=(\matI_{d}-\mat 0_{d})^{-1}(\matI_{d}+\mat 0_{d})\matX=\matX\ .
\]
For the second condition, we have
\[
\left.\frac{\text{d}}{\text{dt}}R_{\matX}^{\text{Cayley}}(t\xi_{\matX})\right|_{t=0}=\matW(\xi_{\matX})\matX=\xi_{\matX}\ ,
\]
where we used $\matX^{\T}\matB\matX=\matI_{p}$ and the definition
of tangent vectors on $\text{St}_{\matB}(p,d)$, (\ref{eq:tangent_st_B_appendix}),
i.e., $\xi_{\matX}^{\T}\matB\matX+\matX^{\T}\matB\xi_{\matX}=\mat 0_{p}$.

Let us consider the inverse of the polar retraction. Suppose that
$\matY=R_{\matX}^{\text{polar}}(\xi_{\matX}).$ Using the definition
of the polar retraction, and reordering the equation we find that
\begin{equation}
\xi_{\matX}=\matY(\mat I_{p}+\xi_{\matX}^{\T}B\xi_{\matX})^{\nicehalf}-\matX\ .\label{eq:polar-xi-appendix}
\end{equation}
Left multiply by $\matX^{\T}\matB$, and recall that $\matX^{\T}\matB\matX=\matI_{p}$,
to find that 
\[
\matX^{\T}\matB\xi_{\matX}=\matX^{\T}\matB\matY(\mat I_{p}+\xi_{\matX}^{\T}B\xi_{\matX})^{\nicehalf}-\matI_{p}\ .
\]
Now using the fact that $\matX^{\T}\matB\xi_{\matX}+\xi_{\matB}^{\T}\matB\matX=\mat 0_{p}$
(since $\xi_{\matX}$ is a tangent vector), we find that 
\[
\matX^{\T}\matB\matY(\mat I_{p}+\xi_{\matX}^{\T}B\xi_{\matX})^{\nicehalf}+(\mat I_{p}+\xi_{\matX}^{\T}B\xi_{\matX})^{\nicehalf}\matY^{\T}\matB\matX-2\matI_{p}=\mat 0_{p}\,.
\]
Thus $\matZ=(\mat I_{p}+\xi_{\matX}^{\T}B\xi_{\matX})^{\nicehalf}$
is SPD solution to ~(\ref{eq:polar-ret-lyapunov}). If we can uniquely
recover $(\mat I_{p}+\xi_{\matX}^{\T}B\xi_{\matX})^{\nicehalf}$ by
solving ~(\ref{eq:polar-ret-lyapunov}) (something we can do in a
small neighborhood of $\matX$ that intersects with the image of the
polar retraction), we can use ~(\ref{eq:polar-xi-appendix}) to invert
the polar retraction.

The derivation of the inverse of the QR retraction is similar. Suppose
that $\matY=R_{\matX}^{\text{QR}}(\xi_{\matX}).$ Using the definition
of the QR-based retraction, and reordering the equation we find that
\begin{equation}
\xi_{\matX}=\matY\mat R-\matX,\label{eq:qr-xi-appendix}
\end{equation}
where $\matR$ is an upper-triangular matrix with strictly positive
elements on its main diagonal such that 
\[
\qf{\matB^{\nicehalf}\left(\matX+\xi_{\matX}\right)}\matR=\matB^{\nicehalf}\left(\matX+\xi_{\matX}\right).
\]
To find $\matR$, left multiply by $\matX^{\T}\matB$, and recall
that $\matX^{\T}\matB\matX=\matI_{p}$ to find that 
\[
\matX^{\T}\matB\xi_{\matX}=\matX^{\T}\matB\matY\matR-\matI_{p}\ .
\]
Now using the fact that $\matX^{\T}\matB\xi_{\matX}+\xi_{\matB}^{\T}\matB\matX=\mat 0_{p}$
(since $\xi_{\matX}$ is a tangent vector), we find that 
\[
\matX^{\T}\matB\matY\matR+\matR^{\T}\matY^{\T}\matB\matX-2\matI_{p}=\mat 0_{p}\,.
\]
Thus, $\matR$ is an upper-triangular matrix with strictly positive
elements on its main diagonal solving to ~(\ref{eq:qr-ret-lyapunov}).
If we can uniquely recover $\matR$ by solving ~(\ref{eq:qr-ret-lyapunov})
(something we can do in a small neighborhood of $\matX$ that intersects
with the image of the QR-based retraction), we can use ~(\ref{eq:qr-xi-appendix})
to invert the QR-based retraction.

We remind here the conditions for a unique solution for  (\ref{eq:qr-ret-lyapunov}).
According to \cite[Eq. (14) and Algorithm 1]{kaneko2012empirical},
 (\ref{eq:qr-ret-lyapunov}) is equivalent to the set of the following
$p$ linear equations
\[
\tilde{\matM}_{i}\tilde{r}_{i}=\b_{i},\ i=1,...,p,
\]
where $\tilde{\matM}_{i}$ is the $i$-th principal minor extracted
from the matrix $\matX^{\T}\matB\matY$, $\tilde{r}_{i}$ is the column-vector
formed by the first $i$ elements of the $i$-th column of the matrix
$\matR$, and $\b_{i}$ is the column-vector whose first $i-1$ elements
are the product 
\[
-[m_{i1},...m_{ij}]\tilde{r}_{j}\ ,
\]
where $j=1,...,i-1$, $m_{ik}$ are elements of the $i$-th row of
$\tilde{\matM}_{i}$, and the $i$-th element of $\b_{i}$ equals
$1$. Thus, this set of linear equations has a unique solution if
and only if all the principal minors of $\matX^{\T}\matB\matY$ are
non-singular. In addition, we also demand that the diagonal elements
of $\matR$ are strictly positive. Note that, since for $\matY$ close
enough to $\matX$ the eigenvalues of $\matX^{\T}\matB\matY$ are
strictly positive, thus $\det(\tilde{\matM}_{i})>0$. Moreover, using
Cramer's rule for $r_{ii}$ the denominator is positive and the nominator
is also positive for $\matY$ close enough to $\matX$ which satisfies
the second constraint on $\matR$.

We also show the derivation of retraction based vector transports
using equations (\ref{eq:stiefelpolarretraction_appendix}) and (\ref{eq:stiefelqrretraction_appendix})
similarly to \cite{zhu2017riemannian}. For ~(\ref{eq:stiefelpolarretraction_appendix})
denote
\[
\matA(t)\coloneqq\mat I_{p}+\left(\eta_{\matX}+t\xi_{\matX}\right)^{\T}B\left(\eta_{\matX}+t\xi_{\matX}\right),
\]
then,
\begin{eqnarray}
\tau_{\eta_{\matX}}^{(\text{polar})} & \coloneqq & \text{D}R_{\matX}^{\text{polar}}(\eta_{\matX})[\xi_{\matX}]\label{eq:vectransportpolar}\\
 & = & \left.\frac{\text{d}}{\text{dt}}R_{\matX}^{\text{polar}}(\eta_{\matX}+t\xi_{\matX})\right|_{t=0}\nonumber \\
 & = & \left.\frac{\text{d}}{\text{dt}}\left[(\matX+\eta_{\matX}+t\xi_{\matX})\left(\matA(t)\right)^{-\nicehalf}\right]\right|_{t=0}\nonumber \\
 & = & \xi_{\matX}\left(\matA(0)\right)^{-\nicehalf}+(\matX+\eta_{\matX})\left.\frac{\text{d}}{\text{dt}}\left(\matA(t)\right)^{-\nicehalf}\right|_{t=0}\nonumber \\
 & = & \xi_{\matX}\left(\matA(0)\right)^{-\nicehalf}-(\matX+\eta_{\matX})\left(\matA(0)\right)^{-\nicehalf}\left.\frac{\text{d}}{\text{dt}}\left(\matA(t)\right)^{\nicehalf}\right|_{t=0}\left(\matA(0)\right)^{-\nicehalf},\nonumber 
\end{eqnarray}
where the last equality is due to the differentiation of the following
two identities
\[
\matI=\left(\matA(t)\right)^{-\nicehalf}\left(\matA(t)\right)^{\nicehalf},
\]
\[
\matA(t)=\left(\matA(t)\right)^{\nicehalf}\left(\matA(t)\right)^{\nicehalf},
\]
which leads to
\[
\mat 0=\frac{\text{d}}{\text{dt}}\left(\matA(t)\right)^{-\nicehalf}\left(\matA(t)\right)^{\nicehalf}+\left(\matA(t)\right)^{-\nicehalf}\frac{\text{d}}{\text{dt}}\left(\matA(t)\right)^{\nicehalf},
\]
\[
\frac{\text{d}}{\text{dt}}\left(\matA(t)\right)^{-\nicehalf}=-\left(\matA(t)\right)^{-\nicehalf}\frac{\text{d}}{\text{dt}}\left(\matA(t)\right)^{\nicehalf}\left(\matA(t)\right)^{-\nicehalf},
\]
and
\[
\frac{\text{d}}{\text{dt}}\matA(t)=\frac{\text{d}}{\text{dt}}\left(\matA(t)\right)^{\nicehalf}\left(\matA(t)\right)^{\nicehalf}+\left(\matA(t)\right)^{\nicehalf}\frac{\text{d}}{\text{dt}}\left(\matA(t)\right)^{\nicehalf},
\]
\[
\xi_{\matX}^{\T}B\eta_{\matX}+\eta_{\matX}^{\T}B\xi_{\matX}+2t\xi_{\matX}^{\T}B\xi_{\matX}=\frac{\text{d}}{\text{dt}}\left(\matA(t)\right)^{\nicehalf}\left(\matA(t)\right)^{\nicehalf}+\left(\matA(t)\right)^{\nicehalf}\frac{\text{d}}{\text{dt}}\left(\matA(t)\right)^{\nicehalf}.
\]
Thus, $\left.\frac{\text{d}}{\text{dt}}\left(\matA(t)\right)^{\nicehalf}\right|_{t=0}$
is a $p\times p$ matrix which is the solution of the following Sylvester
equation:
\begin{eqnarray}
\left.\frac{\text{d}}{\text{dt}}\left(\matA(t)\right)^{\nicehalf}\right|_{t=0}\left(\matA(0)\right)^{\nicehalf}+\left(\matA(0)\right)^{\nicehalf}\left.\frac{\text{d}}{\text{dt}}\left(\matA(t)\right)^{\nicehalf}\right|_{t=0} & = & \xi_{\matX}^{\T}B\eta_{\matX}+\eta_{\matX}^{\T}B\xi_{\matX}.\label{eq:Sylvesterpolarvectrans}
\end{eqnarray}
According to \cite[Theorem 2.4.4.1]{horn2012matrix}, there is a unique
solution to  (\ref{eq:Sylvesterpolarvectrans}) for any $\xi_{\matX}^{\T}B\eta_{\matX}+\eta_{\matX}^{\T}B\xi_{\matX}$,
since $\left(\matA(0)\right)^{\nicehalf}=\left(\mat I_{p}+\eta_{\matX}^{\T}B\eta_{\matX}\right)^{\nicehalf}$
is positive definite ($\eta_{\matX}^{\T}B\eta_{\matX}$ is a symmetric
positive semi-definite matrix) and $-\left(\mat I_{p}+\eta_{\matX}^{\T}B\eta_{\matX}\right)^{\nicehalf}$
is negative definite, thus they have no eigenvalues in common. Solving
 (\ref{eq:Sylvesterpolarvectrans}) costs $O(p^{3})$ (e.g., using
the Bartels--Stewart algorithm \cite{bartels1972solution}). In addition,
to compute this vector transport we need to find the square root of
a $p\times p$ matrix and its inverse which also costs $O(p^{3})$,
compute the product of $\matB$ with matrices, compute the matrix
multiplication of $d\times p$ matrices by $p\times p$ matrices,
of $p\times d$ matrices by $d\times p$ matrices and of $p\times p$
matrices by $p\times p$ matrices. Thus, the total computational cost
of using the vector transport given in  (\ref{eq:vectransportpolar})
is $O(T_{\matB}p+dp^{2})$.

For ~(\ref{eq:stiefelqrretraction}):
\begin{eqnarray}
\tau_{\eta_{\matX}}^{(\text{QR})} & \coloneqq & \text{D}R_{\matX}^{\text{QR}}(\eta_{\matX})[\xi_{\matX}]=\text{D}\text{qf}_{\matB}(\matX+\eta_{\matX})[\xi_{\matX}]\label{eq:vectransportQR}\\
 & = & \matB^{-\nicehalf}\text{D}\text{qf}\left(\matB^{\nicehalf}\left(\matX+\xi_{\matX}\right)\right)[\matB^{\nicehalf}\xi_{\matX}]\nonumber \\
 & = & \matB^{-\nicehalf}\left[\text{qf}\left(\matB^{\nicehalf}\left(\matX+\xi_{\matX}\right)\right)\rho_{\text{skew}}\left(\text{qf}\left(\matB^{\nicehalf}\left(\matX+\xi_{\matX}\right)\right)^{\T}\matB^{\nicehalf}\xi_{\matX}\left(\text{qf}\left(\matB^{\nicehalf}\left(\matX+\xi_{\matX}\right)\right)^{\T}\matB^{\nicehalf}\left(\matX+\xi_{\matX}\right)\right)^{-1}\right)+\right.\nonumber \\
 &  & \left.+\left(\matI_{n}-\text{qf}\left(\matB^{\nicehalf}\left(\matX+\xi_{\matX}\right)\right)\text{qf}\left(\matB^{\nicehalf}\left(\matX+\xi_{\matX}\right)\right)^{\T}\right)\matB^{\nicehalf}\xi_{\matX}\left(\text{qf}\left(\matB^{\nicehalf}\left(\matX+\xi_{\matX}\right)\right)^{\T}\matB^{\nicehalf}\left(\matX+\xi_{\matX}\right)\right)^{-1}\right],\nonumber 
\end{eqnarray}
where the last equality is due to \cite[Example 8.1.5]{AMS09}:
\begin{eqnarray*}
\text{D}\text{qf}(\mat Y)[\mat U] & = & \text{qf}(\mat Y)\rho_{\text{skew}}\left(\text{qf}(\mat Y)^{\T}\mat U\left(\text{qf}(\mat Y)^{\T}\mat Y\right)^{-1}\right)+\left(\matI_{n}-\text{qf}(\mat Y)\text{qf}(\mat Y)^{\T}\right)\mat U\left(\text{qf}(\mat Y)^{\T}\mat Y\right)^{-1},
\end{eqnarray*}
and $\rho_{\text{skew}}(\cdot)$ is the skew-symmetric term of the
decomposition of a square matrix $\matA$ into the sum of a skew-symmetric
term and an upper triangular term, i.e,
\begin{eqnarray*}
\left(\rho_{\text{skew}}(\matA)\right)_{i,j} & = & \begin{cases}
\matA_{i,j} & i>j\\
0 & i=j\\
-\matA_{j,i} & i<j
\end{cases}\ .
\end{eqnarray*}
Computing (\ref{eq:vectransportQR}) can be done in the following
way. First, computing $\matB^{-\nicehalf}\qf{\matB^{\nicehalf}\left(\matX+\xi_{\matX}\right)}$
costs $O\left(T_{\matB}p+dp^{2}\right)$ (see computational cost of
 (\ref{eq:stiefelqrretraction})). Also, computing $\text{qf}\left(\matB^{\nicehalf}\left(\matX+\xi_{\matX}\right)\right)^{\T}\matB^{\nicehalf}$
has the same cost since it is equivalent to computing $\matR^{-\T}\left(\matX+\xi_{\matX}\right)^{\T}\matB$,
where $\matR$ is the $\matR$ matrix of the thin $\matI$-QR decomposition
of $\matB^{\nicehalf}\left(\matX+\xi_{\matX}\right)$, and it can
be found using the Cholesky decomposition of $(\matX+\xi_{\matX})^{\T}\matB(\matX+\xi_{\matX})$.
Applying $\rho_{\text{skew}}(\cdot)$ takes $O(1)$. Finally, all
other computations evolve products of matrices which cost at most
$O(dp^{2})$ and computing the inverse of a $p\times p$ matrix. Thus,
the total computational cost of  (\ref{eq:vectransportQR}) is $O(T_{\matB}p+dp^{2})$.

Both forms of vector transport (\ref{eq:vectransportpolar}) and (\ref{eq:vectransportQR})
satisfy \cite[Definition 8.1.1]{AMS09}. The vector transport based
on the Cayley transform is derived in \cite[Eq. (16)]{zhu2017riemannian}.
It features the same computational complexity as computing the retraction
(\ref{eq:stiefelcayleyretraction_appendix}).

\subsection{\label{subsec:Metric-Related-Notions-1}Metric Related Notions}

We detail the derivation of the Riemannian Hessian that led to ~(\ref{eq:StiefelHessian})
stated in Subsection \ref{subsec:Metric-Related-Notions}. For the
derivation of the Riemannian Hessian we assume that the preconditioning
scheme defining the Riemannian metric is constant, i.e., $\matM_{\matX}\coloneqq\matM$
for all $\matX\in\stiefel_{\matB}(p,d)$. We remark again that ~(\ref{eq:StiefelHessian})
holds also with a non-constant $\matM_{\matX}$ at the critical points.

We use \cite[Definition 5.5.1]{AMS09} of the Riemannian Hessian:
For a real-valued function $f$ on $\text{St}_{\matB}(p,d)$, at a
point $\matX\in\text{St}_{\matB}(p,d)$ the Riemannian Hessian $\hess{f(\matX)}$
is a linear mapping of $T_{\matX}\stiefel_{\matB}(p,d)$ into itself
such that 
\[
\hess{f(\matX)}[\eta_{\matX}]=\nabla_{\eta_{\matX}}\grad{f(\matX)},
\]
for all $\eta_{\matX}\in T_{\matX}\stiefel_{\matB}(p,d)$. In the
previous equation, $\nabla$ is the Riemannian connection, which should
not be confused with the Euclidean gradient.

First, we find the Riemannian connection on $\stiefel_{\matB}(p,d)$
and show that it is the classical directional derivative of vector
fields projected on the tangent space. We can find the Riemannian
connection in a similar manner to the gradient computation performed
in Section \ref{subsec:Metric-Related-Notions} by using \cite[Proposition 5.3.2]{AMS09}:
composing the connection in the ambient space with the projection
on the tangent space. Let $\bar{\nabla}$ be the Levi-Civita connection
on $\R^{d\times p}$ endowed with the metric $\bar{g}$. Let $(e_{1},...,e_{dp})=(\mat E_{11},\mat E_{21},...,\mat E_{d1},\mat E_{12},...,\mat E_{d2},...,\mat E_{dp})$
be the canonical basis of $\mathbb{R}^{d\times p}$, that is matrices
$\mat E_{ij}\in\R^{d\times p}$ such that their only non-zero element
is in the $ij$-th position and its value is $1$. The matrices are
ordered by columns, i.e., for $i=kd+r$ where $k,r\in\mathbb{N}\cup\{0\}$
and $0\leq k\leq p,\ 0\leq r<d$ we have that 
\[
e_{i}=\begin{cases}
E_{r,(k+1)} & r\neq0\\
E_{d,k} & r=0
\end{cases}
\]
(first only the matrices with $1$ in their first column appear, then
in the second column, as so on). Then we have 
\[
\overline{\nabla}_{\overline{\eta}(\cdot)}\overline{\xi}(\cdot)=\sum_{i,j}\left(\overline{\eta^{i}}(\cdot)\,\overline{\xi^{j}}(\cdot)\overline{\nabla}_{e_{i}(\cdot)}e_{j}(\cdot)+\overline{\eta^{i}}(\cdot)\partial_{i}\overline{\xi^{j}}(\cdot)e_{j}(\cdot)\right),
\]
where $\overline{\eta}(\cdot),\overline{\xi}(\cdot),e_{i}(\cdot),\overline{\nabla}_{\overline{\eta}(\cdot)}\overline{\xi}(\cdot),\overline{\nabla}_{e_{i}(\cdot)}e_{j}(\cdot)$
are all vector fields on $\R^{d\times p}$ (i.e., given a point $\matX\in\R^{d\times p}$
the vector field assigns a tangent vector in $T_{\matX}\R^{d\times p}\cong\R^{d\times p}$,
e.g., $\overline{\eta}(\matX)=\overline{\eta}_{\matX}$). In particular,
$\overline{\eta}(\cdot)$ and $\overline{\xi}(\cdot)$ are smooth
local extensions of the vector fields $\eta(\cdot)$ and $\xi(\cdot)$
on $\stiefel_{\matB}(p,d)$ in a neighborhood of $\matX\in\stiefel_{\matB}(p,d)$
in $\R^{d\times p}$, in the sense that for $\matX\in\R^{d\times p}$
the vector fields $\overline{\eta}(\cdot)$ and $\overline{\xi}(\cdot)$
assign the same tangent vectors as $\eta(\cdot)$ and $\xi(\cdot)$.
Note that the vector field $\overline{\nabla}_{\overline{\eta}(\cdot)}\overline{\xi}(\cdot)$
at $\matX$ depends on $\overline{\eta}(\matX)=\overline{\eta}_{\matX}$
(see \cite[Proposition 5.18.]{boumal2020intromanifolds}). Thus, we
can write $\overline{\nabla}_{\overline{\eta}(\cdot)}\overline{\xi}(\cdot)$
at $\matX$ in the following way $\overline{\nabla}_{\overline{\eta}_{\matX}}\overline{\xi}(\matX)$.
In addition, given $\eta_{\matX}\in T_{\matX}\stiefel_{\matB}(p,d)$
and a vector field $\xi(\matX)$ on $\text{St}_{\matB}(p,d)$, the
connection $\overline{\nabla}_{\eta_{\matX}}\xi(\matX)$ is defined
by $\overline{\nabla}_{\overline{\eta}_{\matX}}\overline{\xi}(\matX)$
according to \cite[Equation 5.13]{AMS09} and it does not depend on
the local extension of $\xi(\matX)$. Recall that $\left(\overline{\nabla}_{e_{i}(\cdot)}e_{j}(\cdot)\right)_{k}=\Gamma_{i,j}^{k}$
($k$-th coordinate of $\overline{\nabla}_{e_{i}}e_{j}$) are the
Christoffel symbols. These symbols determine the connection $\overline{\nabla}$
uniquely, using the Fundamental Theorem of Riemannian Geometry for
the Levi-Civita connection. The Christoffel symbols can be calculated
using
\[
\Gamma_{i,j}^{k}=\frac{1}{2}\sum_{l=1}^{dp}g^{kl}(\partial_{i}g_{lj}+\partial_{j}g_{li}-\partial_{l}g_{ij}),
\]
where $g^{kl}$ is the $(k,l)$th entry of the inverse of the matrix
$dp\times dp$ matrix $\mat G$ which is defined by
\[
\left(\mat G\right)_{kl}\coloneqq g_{kl}=g_{\matX}(e_{k},e_{l})=g_{\matX}(E_{ij},E_{hm})=\Trace{E_{ij}^{\T}\matM E_{hm}}=\begin{cases}
0 & ,\ j\neq m\\
\matM_{ih} & ,\ j=m
\end{cases}\ .
\]
Since the components of the matrix $\matM$ do not depend on $\matX$
and on $(e_{1},...,e_{dp})$ (it is a constant matrix) we have $\forall i,j,k:\quad\Gamma_{i,j}^{k}=0$.
Therefore, $\overline{\nabla}$ is reduced to the classical directional
derivative in $\mathbb{R}^{d\times p}$
\[
\overline{\nabla}_{\overline{\eta}_{\matX}}\overline{\xi}(\matX)=\sum_{j=1}^{dp}\sum_{i=1}^{dp}\left(\overline{\eta_{\matX}^{i}}\partial_{i}\overline{\xi^{j}}(\matX)e_{j}\right)=\matJ_{\overline{\xi}(\matX)}\overline{\eta}_{\matX},
\]
where $\matJ_{\overline{\xi}(\matX)}\overline{\eta}_{\matX}$ denotes
the Jacobian matrix of $\overline{\xi}(\matX)$ at $\matX$ in the
direction $\overline{\eta}_{\matX}$. Now that we have the connection
on the ambient space $\R^{d\times p}$, which is a Riemannian manifold,
we can compute the connection on the submanifold $\stiefel_{\matB}(p,d)$.
Given $\eta_{\matX}\in T_{\matX}\stiefel_{\matB}(p,d)$ and a vector
field $\xi(\matX)$ on $\stiefel_{\matB}(p,d)$, the Riemannian connection
is (written, as usual, in terms of ambient coordinates):
\begin{equation}
\nabla_{\eta_{\matX}}\xi(\matX)=\Pi_{\matX}\left(\overline{\nabla}_{\overline{\eta}_{\matX}}\overline{\xi}(\matX)\right)=\Pi_{\matX}\left(\matJ_{\overline{\xi}(\matX)}\eta_{\matX}\right)\label{eq:connection}
\end{equation}
 where $\eta_{\matX}=\overline{\eta}_{\matX}$ and $\overline{\xi}(\cdot)$
is any smooth local extension of $\xi(\cdot)$ in a neighborhood of
$\matX\in\stiefel_{\matB}(p,d)$ in $\R^{d\times p}$.

Next, we can find the Riemannian Hessian using  (\ref{eq:connection}),
the product rule for derivation and according to \cite{absil2013extrinsic}:
\begin{eqnarray}
\hess{f(\matX)[\eta_{\matX}]} & = & \nabla_{\eta_{\matX}}\grad{f(\matX)}\nonumber \\
 & = & \Pi_{\matX}\left(\matJ_{h(\matX)}\eta_{\matX}\right)\nonumber \\
 & = & \Pi_{\matX}\left[\matP_{\matX}\matM^{-1}\nabla^{2}\bar{f}(\matX)\eta_{\matX}+(\text{D\ensuremath{\Pi_{\matX}}})[\eta_{\matX}]\matM^{-1}\nabla\bar{f}(\matX)\right]\nonumber \\
 & = & \Pi_{\matX}\left(\matM^{-1}\nabla^{2}\bar{f}(\matX)\eta_{\matX}\right)+\Pi_{\matX}\left((\text{D\ensuremath{\Pi_{\matX}}})(\matX)[\eta_{\matX}]\matM^{-1}\nabla\bar{f}(\matX)\right)\label{eq:RHess-Stiefel}
\end{eqnarray}
where $\nabla\bar{f}(\matX)$ and $\nabla^{2}\bar{f}(\matX)$ are
the Euclidean gradient and Hessian (respectively) of $\bar{f}$ and
\[
h:\R^{d\times p}\to\R^{d\times p},\,\,h(\matX)=\Pi_{\matX}\left(\matM^{-1}\nabla\bar{f}(\matX)\right).
\]
Note that for $\matX\in\stiefel_{\matB}(p,d)$ we have $h(\matX)=\grad f(\matX)$
so $h$ is a smooth local extension of the vector field $\grad f$
to $\R^{d\times p}$, and its Jacobian is calculated as follows
\[
\matJ_{h(\matX)}\eta_{\matX}=(\text{D\ensuremath{\Pi_{\matX}}})[\eta_{\matX}]\matM^{-1}\nabla\bar{f}(\matX)+\Pi_{\matX}\left(\matM^{-1}\nabla^{2}\bar{f}(\matX)\eta_{\matX}\right)\ ,
\]
where $(\text{D\ensuremath{\Pi_{\matX}}})[\eta_{\matX}]$ (here and
in  (\ref{eq:RHess-Stiefel})) is the derivative at $\matX$ along
$\eta_{\matX}$ of the function that maps $\matX$ to $\Pi_{\matX}$.

The main challenge in computing the Riemannian Hessian from  (\ref{eq:RHess-Stiefel})
is in computing $(\text{D\ensuremath{\Pi_{\matX}}})[\eta_{\matX}]$.
In order to circumvent this issue, we use a simple modification of
a result found in \cite{absil2013extrinsic} to the case in which
the Riemannian metric induced from $\R^{d\times p}$ on any Riemannian
submanifold of $\R^{d\times p}$ is of the form $g_{\matX}(\xi_{\matX},\eta_{\matX})=\Trace{\xi_{\matX}^{\T}\matM\eta_{\matX}}$
where $\matM\in\R^{d\times d}$ is any constant, SPD matrix. In order
to so, first we introduce the notion of the \emph{Weingarten map}.
\begin{defn}
\label{def:Weingarten-map}(\cite[Section 6.1]{absil2009all}, \cite[Definition 1]{absil2013extrinsic})
Given a Riemannian manifold ${\cal M}$, a point $\x\in{\cal M}$
on the manifold, a tangent vector $\eta_{\x}\in T_{\x}{\cal M}$ at
$\x$, and a normal vector $\u_{\x}\in(T_{\x}{\cal M})^{\perp}$,
we define the Weingarten map by 
\begin{equation}
W_{\x}\left(\eta_{\x},\u_{\x}\right)\coloneqq-\Pi_{\x}(\text{D}\u(\x)[\eta_{\x}])\label{eq:weingartendef}
\end{equation}
 where $\u(\cdot)$ is a smooth normal vector field on ${\cal M}$
which satisfies $\u(\x)=\u_{\x}$.
\end{defn}

For the manifold $\stiefel_{\matB}(p,d)$, viewed as an embedded submanifold
of $\R^{d\times p}$,  (\ref{eq:weingartendef}) reduces to 
\begin{eqnarray*}
W_{\matX}\left(\eta_{\matX},\matU(\matX)\right) & = & -\Pi_{\matX}\left(\matJ_{\bar{\mat U}(\matX)}\eta_{\matX}\right),
\end{eqnarray*}
where $\bar{\mat U}(\cdot)$ is any smooth local extension of the
normal vector field $\matU(\cdot)$ such that $\matU(\matX)=\mat U_{\matX}$
on $\stiefel_{\matB}(p,d)$. Now, that at a point $\matX\in\stiefel_{\matB}(p,d)$
any normal vector is of the form $\matU_{\matX}=\matM^{-1}\matB\matX\matS_{\matX}$
for some $\matS_{\matX}\in{\cal {\cal S}_{\text{sym}}}(p)$. Left
multiplying by $\matX^{\T}\matM$ we get $\matX^{\T}\matM\matU_{\matX}=\matX^{\T}\matM\matM^{-1}\matB\matX\mat S_{\matX}=\mat S_{\matX}$.
Now we can define a normal field on $\stiefel_{\matB}(p,d)$ by the
formula $\matU(\matX)=\matM^{-1}\matB\matX\mat S_{\matX}$ such that
$\matU(\matX)=\matU_{\matX}$ with $\mat S_{\matX}=\matX^{\T}\matM\matU_{\matX}$
such that $\mat S_{\matX}\in{\cal {\cal S}_{\text{sym}}}(p)$. The
vector field can be extended to $\R^{d\times p}$ by the same formula
such that $\bar{\mat U}(\cdot)$ and $\matU(\cdot)$ coincide on $\stiefel_{\matB}(p,d)$.
Next, we calculate the Jacobian of $\bar{\mat U}(\matX)$ at the direction
$\eta_{\matX}$: 
\[
\matJ_{\bar{\mat U}(\matX)}\eta_{\matX}=\matM^{-1}\matB\eta_{\matX}\mat S_{\matX}\ ,
\]
Therefore the Weingarten map for $\stiefel_{\matB}(p,d)$ is 
\begin{eqnarray*}
W_{\matX}\left(\eta_{\matX},\matU_{\matX}\right) & = & -\Pi_{\matX}(\matM^{-1}\matB\eta_{\matX}\mat S_{\matX})\\
 & = & -\Pi_{\matX}\left(\matM^{-1}\matB\eta_{\matX}\left(\matX^{\T}\matM\matU_{\matX}\right)\right)\ .
\end{eqnarray*}

The following lemma is a simple modification of \cite[Theorem 1]{absil2013extrinsic}.
Although the proof is almost identical, we include it here for completeness.
\begin{lem}
\label{lem:modificationoftheorem} For the Riemannian submanifold
$\stiefel_{\matB}(p,d)$ of $\R^{d\times p}$ endowed with $\bar{g}_{\matX}(\bar{\xi}_{\matX},\bar{\eta}_{\matX})=\Trace{\bar{\xi}_{\matX}^{\T}\matM\bar{\eta}_{\matX}}$
we have 
\begin{eqnarray*}
W_{\matX}\left(\eta_{\matX},\Pi_{\matX}^{\perp}\left(\matM^{-1}\matU\right)\right) & = & \Pi_{\matX}\left((\text{D}\Pi_{\matX})[\eta_{\matX}]\left(\matM^{-1}\matU\right)\right)\\
 & = & \Pi_{\matX}\left((\text{D}\Pi_{\matX})[\eta_{\matX}]\left(\Pi_{\matX}^{\perp}\left(\matM^{-1}\matU\right)\right)\right)\ ,
\end{eqnarray*}
 for all $\matX\in\stiefel_{\matB}(p,d)$, $\eta_{\matX}\in T_{\matX}\stiefel_{\matB}(p,d)$
and $\matU\in\R^{d\times p}$.
\end{lem}

\begin{proof}
First, we show that 
\begin{equation}
\Pi_{\matX}\left((\text{D}\Pi_{\matX})[\eta_{\matX}]\left(\cdot\right)\right)=\Pi_{\matX}\left((\text{D}\Pi_{\matX})[\eta_{\matX}]\left(\Pi_{\matX}^{\perp}\left(\cdot\right)\right)\right)\,\label{eq:prpositionproof1}
\end{equation}
holds. Then applying both sided on $\matM^{-1}\matU$ gives us the
equality 
\[
\Pi_{\matX}\left((\text{D}\Pi_{\matX})[\eta_{\matX}]\matM^{-1}\matU\right)=\Pi_{\matX}\left((\text{D}\Pi_{\matX})[\eta_{\matX}]\Pi_{\matX}^{\perp}\left(\matM^{-1}\matU\right)\right).
\]
To show this, we take the directional derivative of the equality $\Pi_{\matX}\left(\Pi_{\matX}^{\perp}\left(\cdot\right)\right)=\mat 0$
in the direction $\eta_{\matX}$, and we use $\Pi_{\matX}^{\perp}\left(\cdot\right)=\left(\id_{T_{\matX}\R^{d\times p}}-\Pi_{\matX}\right)\left(\cdot\right)$
to get
\begin{eqnarray*}
\mat 0 & = & (\text{D}\Pi_{\matX})[\eta_{\matX}]\Pi_{\matX}^{\perp}\left(\cdot\right)+\Pi_{\matX}\left((\text{D}\Pi_{\matX}^{\perp})[\eta_{\matX}]\left(\cdot\right)\right)\\
 & = & (\text{D}\Pi_{\matX})[\eta_{\matX}]\Pi_{\matX}^{\perp}\left(\cdot\right)-\Pi_{\matX}\left((\text{D}\Pi_{\matX})[\eta_{\matX}]\left(\cdot\right)\right).
\end{eqnarray*}
Substituting any tangent vector in both sides of the equation nullifies
the term $(\text{D}\Pi_{\matX})[\eta_{\matX}]\Pi_{\matX}^{\perp}\left(\cdot\right)$.
Thus, we substitute $\Pi_{\matX}$ and use $\Pi_{\matX}^{\perp}\left(\Pi_{\matX}\left(\cdot\right)\right)=\mat 0$
\begin{eqnarray*}
\mat 0 & = & \Pi_{\matX}\left((\text{D}\Pi_{\matX})[\eta_{\matX}]\left(\Pi_{\matX}\left(\cdot\right)\right)\right).
\end{eqnarray*}
Finally, to get (\ref{eq:prpositionproof1}), we use $\id_{T_{\matX}\R^{d\times p}}\left(\cdot\right)=\left(\Pi_{\matX}+\Pi_{\matX}^{\perp}\right)\left(\cdot\right)$
and get
\begin{eqnarray*}
\Pi_{\matX}\left((\text{D}\Pi_{\matX})[\eta_{\matX}]\left(\cdot\right)\right) & = & \Pi_{\matX}\left((\text{D}\Pi_{\matX})[\eta_{\matX}]\left(\left(\Pi_{\matX}+\Pi_{\matX}^{\perp}\right)\left(\cdot\right)\right)\right)\\
 & = & \Pi_{\matX}\left((\text{D}\Pi_{\matX})[\eta_{\matX}]\left(\Pi_{\matX}^{\perp}\left(\cdot\right)\right)\right).
\end{eqnarray*}

To conclude the proof we show that $W_{\matX}\left(\eta_{\matX},\Pi_{\matX}^{\perp}\left(\matM^{-1}\matU\right)\right)=\Pi_{\matX}\left((\text{D}\Pi_{\matX})[\eta_{\matX}]\left(\matM^{-1}\matU\right)\right)$.
Note that for embedded submanifolds of $\R^{d\times p}$ with a metric
derived from $\matM$, the Weingarten map reduces to $W_{\matX}\left(\eta_{\matX},\matU_{\matX}\right)=-\Pi_{\matX}\left(\matJ_{\matU(\matX)}\eta_{\matX}\right)$.
Using Definition \ref{def:Weingarten-map} along this observation,
we have 
\begin{eqnarray*}
W_{\matX}\left(\eta_{\matX},\Pi_{\matX}^{\perp}\left(\matM^{-1}\matU\right)\right) & = & -\Pi_{\matX}\left(\matJ_{\Pi_{\matX}^{\perp}\left(\matM^{-1}\matU(\matX)\right)}\eta_{\matX}\right)\\
 & = & -\Pi_{\matX}\left((\text{D}\Pi_{\matX}^{\perp})[\eta_{\matX}]\left(\matM^{-1}\matU_{\matX}\right)\right)-\Pi_{\matX}\left(\Pi_{\matX}^{\perp}\left(\matJ_{\matM^{-1}\matU(\matX)}\eta_{\matX}\right)\right)\\
 & = & \Pi_{\matX}\left((\text{D}\Pi_{\matX})[\eta_{\matX}]\left(\matM^{-1}\matU_{\matX}\right)\right),
\end{eqnarray*}
where in the last equality we used $\Pi_{\matX}\left(\Pi_{\matX}^{\perp}\left(\cdot\right)\right)=0$
and $\Pi_{\matX}^{\perp}\left(\cdot\right)=\left(\id_{T_{\matX}\R^{d\times p}}-\Pi_{\matX}\right)\left(\cdot\right)$.
\end{proof}
As a consequence of Lemma \ref{lem:modificationoftheorem}, we can
replace $\Pi_{\matX}\left((\text{D}\Pi_{\matX})[\eta_{\matX}]\left(\matM^{-1}\nabla\bar{f}(\matX)\right)\right)$
by $W_{\matX}\left(\eta_{\matX},\Pi_{\matX}^{\perp}\left(\matM^{-1}\nabla\bar{f}(\matX)\right)\right)$
in  (\ref{eq:RHess-Stiefel}). Therefore the expression for the Riemannian
Hessian becomes 
\[
\hess{f(\matX)[\eta_{\matX}]}=\Pi_{\matX}\left(\matM^{-1}\nabla^{2}\bar{f}(\matX)\eta_{\matX}\right)+W_{\matX}\left(\eta_{\matX},\Pi_{\matX}^{\perp}\left(\matM^{-1}\nabla\bar{f}(\matX)\right)\right).
\]
In particular, the Riemannian Hessian on $\stiefel_{\matB}(p,d)$
is 
\begin{eqnarray*}
\hess{f(\matX)[\eta_{\matX}]} & = & \Pi_{\matX}\left(\matM^{-1}\nabla^{2}\bar{f}(\matX)\eta_{\matX}\right)-\Pi_{\matX}\left(\matM^{-1}\matB\eta_{\matX}\left(\matX^{T}\matM\left(\Pi_{\matX}^{\perp}\left(\matM^{-1}\nabla\bar{f}(\matX)\right)\right)\right)\right)
\end{eqnarray*}
Note that some simplification of these expressions can be made by
using $\Pi_{\matX}^{\perp}=\text{id}_{T_{\matX}\text{St}_{\matB}(p,d)}-\Pi_{\matX}$:
\begin{eqnarray*}
\hess{f(\matX)[\eta_{\matX}]} & = & \Pi_{\matX}\left(\matM^{-1}\nabla^{2}\bar{f}(\matX)\eta_{\matX}\right)-\Pi_{\matX}\left(\matM^{-1}\matB\eta_{\matX}\left(\matX^{T}\nabla\bar{f}(\matX)-\matX^{\T}\matM\left(\Pi_{\matX}\left(\matM^{-1}\nabla\bar{f}(\matX)\right)\right)\right)\right)\\
 & = & \Pi_{\matX}\left(\matM^{-1}\nabla^{2}\bar{f}(\matX)\eta_{\matX}\right)-\Pi_{\matX}\left(\matM^{-1}\matB\eta_{\matX}\left(\matX^{T}\nabla\bar{f}(\matX)-\matX^{\T}\matM\grad f(\matX)\right)\right)\ .
\end{eqnarray*}

\section{Experiments With $p=2$}

Similarly to the experiments in Subsection \ref{subsec:CCA}, we perform
experiments with the \noun{MEDIANILL }dataset to demonstrate CCA for
$p=2$. We use the same choices for Riemannian metric: the trivial
choice of a unit matrix $\matM=\mat I_{d}$, the standard but expensive
choice $\matM=\Sigma$, and four approximations of $\Sigma$ via the
(exact) sketched preconditioning strategy. Finding the top two correlations
requires the \emph{von Neumann cost function} \cite{von1937some}
formulation:
\[
\begin{array}{c}
\max\Trace{\matU^{\T}\Sigma_{\x\y}\matV\matN}\\
\text{subject to}\\
\matU^{\T}\Sigma_{\x\x}\matU=\matI_{p}\\
\matV^{\T}\Sigma_{\y\y}\matV=\matI_{p}
\end{array}
\]
where $\matN=\diag{\mu_{1},\mu_{2}}$ and any $\mu_{1}>\mu_{2}>0$
(here we take $\mu_{1}=5$ and $\mu_{2}=1$). The corresponding Riemannian
components are constructed in a similar manner to Subsection \ref{subsec:CCA}.

\begin{figure}[t]
\begin{centering}
\begin{tabular}{ccc}
\includegraphics[width=0.38\columnwidth]{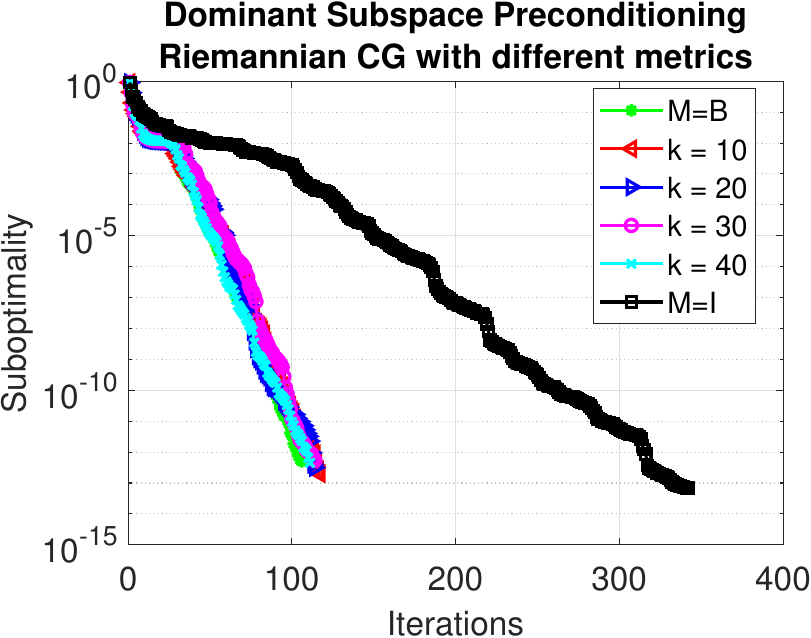} & ~ & \includegraphics[width=0.38\columnwidth]{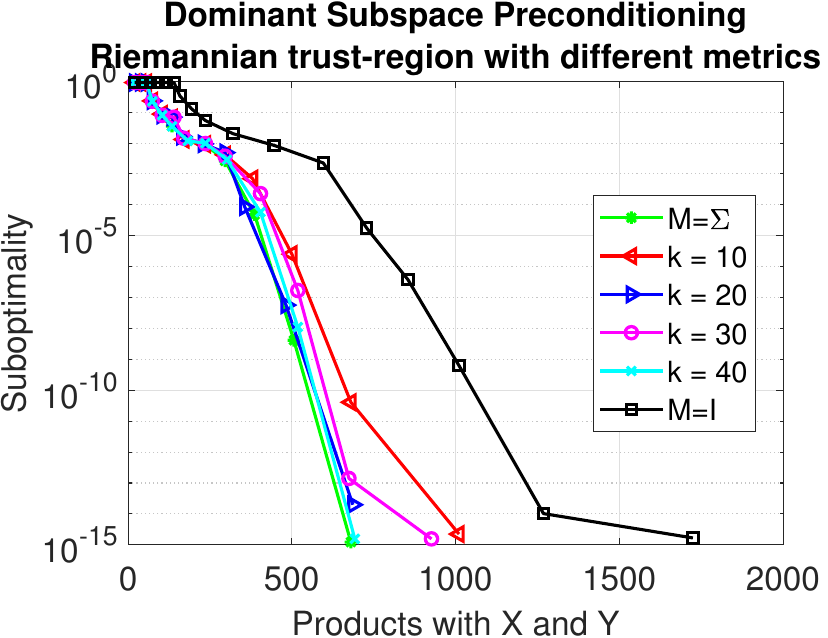}\tabularnewline
\end{tabular}
\par\end{centering}
\caption{\label{fig:Results-for-CCA-1}Results for CCA with Riemannian conjugate-gradient
(left - suboptimality vs. \#iterations) and Riemannain trust-region
(right - suboptimality vs. products with the data matrices) with various
choices of metrics for $p=2$. The number of leading eigenvalues used
to form the Dominant Subspace Preconditioner is denoted by $k$.}
\end{figure}

The graphs in Fig. \ref{fig:Results-for-CCA-1} demonstrate that the
choice $\matM=\Sigma$ leads to the lowest iteration count. This observation
is also supported by the condition number of the Riemannian Hessian
at the optimum, which is evaluated using \noun{Manopt}: the lowest
condition number, $115.68$, is achieved when $\matM=\Sigma$, and
the highest, $805.2$, for $\matM=\mat I_{d}$.

\end{document}